\documentclass[12pt,a4paper]{article}
\usepackage{epsfig,latexsym,amsfonts,amssymb,amsmath,amscd,
graphics,theorem,epic,epstopdf}

\makeatletter
\def\author#1{\gdef\autrun{\def\and{\unskip, }#1}\gdef\@author{#1}}
\def\address#1{{\def\and{\\\hspace*{18pt}}\renewcommand{\thefootnote}{}%
\footnote {#1}}%
}
\makeatother
\def\email#1{e-mail: #1}

\oddsidemargin36pt \evensidemargin74pt  \theoremstyle{plain}
\newtheorem{lemma}{Lemma}[section]
\newtheorem{proposition}{Proposition}[section]
\newtheorem{remark}{Remark}[section]

\newtheorem{theorem}{Theorem}[section]

\newtheorem{corollary}{Corollary}[section]
{\theorembodyfont{\rmfamily} 
  
\begin{document}
\baselineskip=15pt
\newcommand{\pperp}{\hbox{$\perp\hskip-6pt\perp$}}
\newcommand{\ssim}{\hbox{$\hskip-2pt\sim$}}
\newcommand{\N}{{\mathbb N}}\newcommand{\Delp}{{\Pi}}
\newcommand{\SSS}{{\mathbb S}}
\newcommand{\Z}{{\mathbb Z}}
\newcommand{\R}{{\mathbb R}}\newcommand{\F}{{\mathbb F}}
\newcommand{\mm}{{\mathfrak m}}
\newcommand{\C}{{\mathbb C}}\newcommand{\B}{{\mathbb B}}
\newcommand{\Q}{{\mathbb Q}}\newcommand{\T}{{\mathbb T}}\newcommand{\K}{{\mathbb
K}}\newcommand{\D}{{\mathcal D}}
\newcommand{\PP}{{\mathbb P}}
\newcommand{\st}{{*}}\newcommand{\coker}{{\operatorname{coker}}}\newcommand{\Prec}{{\operatorname{Prec}}}
\newcommand{\ev}{{\operatorname{ev}}}
\newcommand{\Id}{{\operatorname{Id}}}\newcommand{\irr}{{\operatorname{irr}}}
\newcommand{\Sym}{{\operatorname{Sym}}}
\newcommand{\fix}{{\operatorname{fix}}}\newcommand{\re}{{\operatorname{re}}}
\newcommand{\ima}{{\operatorname{im}}}
\newcommand{\nod}{{\operatorname{nod}}}\newcommand{\lcm}{{\operatorname{lcm}}}
\newcommand{\oeps}{{\overline\eps}}\newcommand{\Area}{{\operatorname{Area}}}
\newcommand{\codim}{{\operatorname{codim}}}\newcommand{\sol}{{\operatorname{sol}}}
\newcommand{\oDel}{{\widetilde\Del}}
\newcommand{\odd}{{\operatorname{odd}}}\newcommand{\even}{{\operatorname{even}}}
\newcommand{\real}{{\operatorname{Re}}}\newcommand{\ind}{{\operatorname{ind}}}
\newcommand{\conv}{{\operatorname{conv}}}\newcommand{\Leaf}{{\operatorname{Lf}}}
\newcommand{\Span}{{\operatorname{Span}}}
\newcommand{\Ker}{{\operatorname{Ker}}}\newcommand{\cheb}{{\operatorname{cheb}}}
\newcommand{\Fix}{{\operatorname{Fix}}}\newcommand{\ord}{{\operatorname{ord}}}
\newcommand{\sign}{{\operatorname{sign}}}
\newcommand{\Log}{{\operatorname{Log}}}
\newcommand{\Ima}{{\operatorname{Im}}}
\newcommand{\oi}{{\overline i}}
\newcommand{\oj}{{\overline j}}
\newcommand{\ob}{{\overline b}}
\newcommand{\os}{{\overline s}}
\newcommand{\oa}{{\overline a}}
\newcommand{\oy}{{\overline y}}
\newcommand{\ow}{{\overline w}}
\newcommand{\ou}{{\overline u}}
\newcommand{\ot}{{\overline t}}
\newcommand{\oz}{{\overline z}}
\newcommand{\newi}{i}\newcommand{\bd}{{\boldsymbol d}}
\newcommand{\newj}{j}
\newcommand{\newm}{m}
\newcommand{\newl}{{\ell}}
\newcommand{\bw}{{\boldsymbol w}}\newcommand{\bi}{{\boldsymbol i}}
\newcommand{\bx}{{\boldsymbol p}}\newcommand{\bp}{{\boldsymbol p}}
\newcommand{\bpp}{{\boldsymbol P}}\newcommand{\bq}{{\boldsymbol q}}
\newcommand{\by}{{\boldsymbol q}}\newcommand{\bn}{{\boldsymbol n}}
\newcommand{\bz}{{\boldsymbol z}}
\newcommand{\eps}{{\varepsilon}}
\newcommand{\proofend}{\hfill$\Box$\bigskip}
\newcommand{\Int}{{\operatorname{Int}}}
\newcommand{\pr}{{\operatorname{pr}}}
\newcommand{\grad}{{\operatorname{grad}}}
\newcommand{\rk}{{\operatorname{rk}}}
\newcommand{\im}{{\operatorname{Im}}}
\newcommand{\sk}{{\operatorname{sk}}}

\newcommand{\const}{{\operatorname{const}}}
\newcommand{\Sing}{{\operatorname{Sing}}}
\newcommand{\conj}{{\operatorname{Conj}}}
\newcommand{\cconj}{{\operatorname{conj}}}
\newcommand{\Pic}{{\operatorname{Pic}}}
\newcommand{\PicP}{{\operatorname{Pic_+}}}
\newcommand{\PicPR}{{\operatorname{Pic^\R_+}}}
\newcommand{\PicPP}{{\operatorname{Pic_{++}}}}
\newcommand{\PicPPR}{{\operatorname{Pic^\R_{++}}}}
\newcommand{\Crit}{{\operatorname{Crit}}}
\newcommand{\Ch}{{\operatorname{Ch}}}
\newcommand{\discr}{{\operatorname{discr}}}
\newcommand{\card}{{\operatorname{card}}}
\newcommand{\Tor}{{\operatorname{Tor}}}
\newcommand{\Conj}{{\operatorname{Conj}}}
\newcommand{\val}{{\operatorname{val}}}
\newcommand{\Val}{{\operatorname{Val}}}
\newcommand{\res}{{\operatorname{res}}}
\newcommand{\add}{{\operatorname{add}}}
\newcommand{\BH}{{\operatorname{bh}}}
\newcommand{\tmu}{{\C\mu}}
\newcommand{\ov}{{\overline v}}\newcommand{\on}{{\overline n}}
\newcommand{\ox}{{\overline{x}}}
\newcommand{\tet}{{\theta}}
\newcommand{\Del}{{\Delta}}
\newcommand{\bet}{{\beta}}
\newcommand{\kap}{{\kappa}}
\newcommand{\del}{{\delta}}
\newcommand{\sig}{{\sigma}}
\newcommand{\alp}{{\alpha}}
\newcommand{\Sig}{{\Sigma}}
\newcommand{\Gam}{{\Gamma}}
\newcommand{\gam}{{\gamma}}
\newcommand{\Lam}{{\Lambda}}
\newcommand{\lam}{{\lambda}}
\newcommand{\SC}{{SC}}
\newcommand{\MC}{{MC}}
\newcommand{\nek}{{,...,}}
\newcommand{\cim}{{c_{\mbox{\rm im}}}}
\newcommand{\op}{{\overline p}}
\newcommand{\CL}{{\mathcal L}}
\newcommand{\CH}{{\mathcal CH}}
\newcommand{\Sp}{{\mathbb S}}
\newcommand{\CX}{{\mathcal X}}
\newcommand{\CV}{{\mathcal V}}

\newcommand{\WW}{{\mathbb W}}

\title{Welschinger invariants \\
of real del Pezzo surfaces of degree $\ge2$}
\author{Ilia Itenberg \and Viatcheslav Kharlamov \and Eugenii Shustin}

\date{}
\maketitle

\address{I. Itenberg: Universit\'e Pierre et Marie Curie,
Institut de Math\'ematiques de
Jussieu - Paris Rive Gauche, 4 place Jussieu, 75252 Paris Cedex 5, France
\hskip5pt and \hskip5pt
D\'epartement de Math\'emati\-ques et Applications,
Ecole Normale Sup\'erieure,
45 rue d'Ulm, 75230 Paris Cedex 5, France;
\email{ilia.itenberg@imj-prg.fr}; fax: +33-144276366 \and V. Kharlamov: Universit\'{e}
de Strasbourg and IRMA, 7 rue Ren\'{e}-Descartes, 67084 Strasbourg
Cedex, France; \email{viatcheslav.kharlamov@math.unistra.fr}; fax: +33-390240328 \and E.
Shustin: School of Mathematical Sciences, Raymond and Beverly
Sackler Faculty of Exact Sciences, Tel Aviv University, Ramat Aviv,
69978 Tel Aviv, Israel; \email{shustin@post.tau.ac.il}; fax: +972-3-6409357}

\begin{abstract}
We compute the purely real Welschinger invariants, both original
and
modified,
for all real del Pezzo surfaces
of degree $\ge2$. We show that
under some
conditions,
for
such a surface $X$
and
a
real nef and big divisor class
$D\in\Pic(X)$, through any generic
collection
of $-DK_X-1$ real
points lying on a
connected component of the real part $\R X$ of $X$
one can
trace a real rational curve $C\in|D|$.
This is derived from the
positivity of
appropriate Welschinger invariants. We furthermore show that these
invariants are asymptotically equivalent, in the logarithmic scale,
to the corresponding genus
zero Gromov-Witten invariants.
Our approach
consists in a conversion of
Shoval-Shustin
recursive formulas counting complex
curves on
the plane blown up at seven points
and of
Vakil's extension
of the Abramovich-Bertram formula for Gromov-Witten invariants
into formulas computing real
enumerative invariants.

\medskip\noindent {\bf MSC2010}:
Primary 14N10. Secondary 14P05, 14N35.

\medskip\noindent {\bf Keywords}: real rational curves,
enumerative geometry, Welschinger invariants, Caporaso-Harris
formula, Abramovich-Bertram-Vakil formula.
\end{abstract}

\tableofcontents

\bigskip
\bigskip

{\hskip3in {\it The cashier gave us a sad smile, took a small hammer out of her mouth, and moving her nose slightly back and forth, she said:

\hskip3in - In my opinion, a seven comes after an eight, only if an eight comes after a seven.

}}

{\hskip3in (Daniil Kharms "A sonnet"
)}

\section{Introduction}\label{intro}

Welschinger invariants can be regarded as real analogues of genus zero Gromov-Witten invariants.
They were introduced in \cite{W1}, \cite{W2}
and count, with appropriate signs, the real rational pseudo-holomorphic
curves which pass through given real collections of points in a given real rational symplectic four-fold.
In the case of real del Pezzo surfaces, the Welschinger count is equivalent to the
enumeration of
real rational algebraic curves. In the present paper, we continue the study of purely real
Welschinger invariants (that is, Welschinger invariants in the situation when all
the point constraints are real)
of del Pezzo surfaces. These invariants, as well as their modifications
introduced in \cite{IKS7}, can be used to prove the existence of interpolating real rational
curves.

As we proved in \cite{IKS, IKS2,  IKS4, IKS6,IKS7}, if $X$ is
either the plane blown up at
$a$ real points and
$b$ pairs of complex conjugate points, where
$a+2b\le6$, $b\le1$, or  a minimal two-component real conic bundle over
$\PP^1$,
or a two-component
real cubic surface, then the (modified) Welschinger invariants
of $X$ are positive and are asymptotically equivalent in the logarithmic scale
to the corresponding Gromov-Witten invariants.
These results  not only  prove
the existence of interpolating real rational curves,
but also show their abundance.

In the present paper, we extend these results to almost all del Pezzo surfaces of degree $\geq 2$ (see Theorems \ref{t9}
and \ref{t11})
and,
in particular,
cover all the missing cases in degree $\geq 3$. The main novelty is the use of nodal del Pezzo surfaces
in a way which is similar to Vakil's approach to computation of Gromov-Witten invariants
of the plane blown up at six points \cite{Va}. We derive new real Caporaso-Harris type formulas
(see Theorems \ref{t1} and \ref{t3})
and real analogues of Abramovich-Bertram-Vakil formula
\cite{AB,Va} (see Theorems \ref{t2} and \ref{t4}).
These formulas combined together allow one to compute the
purely real Welschinger invariants of all real del Pezzo surfaces of degree $\ge2$ from  finitely many
explicitly determined initial values (see Propositions \ref{ini1} and \ref{ini2}).

As a technical tool,
we introduce certain numbers (called ordinary $w$-numbers and sided $w$-numbers) that count with signs some specifically constrained real rational curves
on real nodal del Pezzo surfaces, and exhibit a case when sided $w$-numbers are independent of the choice of point constraints
(see Corollary \ref{c2}).

A
new phenomenon for del Pezzo surfaces of degree $2$ is the absence of
real rational curves in some cases
(see
Section \ref{sec414}).
In this regard, note that in the case of multicomponent
real surfaces, the
original Welschinger invariants
often happen to
vanish
(see \cite[Proposition 3.3]{BP}).
However, by (\ref{ee20}) in Theorem \ref{t9},
in
many
situations
such a vanishing is not
related to the
non-existence of real rational curves, but only
states that the real rational curves under consideration cancel each
other when supplied with the original Welschinger signs.

Several
results related to that of the present paper should be mentioned here.
When we were working on Theorems \ref{t2} and \ref{t4}, E.~Brugall\'e
and
N.~Puignau communicated to us
similar real
versions of Abramovich-Bertram-Vakil formula in the case of del Pezzo surfaces of degree $\geq 3$;
afterwards, they extended
these formulas to the symplectic setting and arbitrary real rational
symplectic $4$-manifolds, see~\cite{BP}.\footnote{When our text was finished, there appeared a new preprint by E. Brugall\'e \cite{B-new},
where he uses a slightly different approach and besides introduces floor diagrams with respect to a conic;
this technique allows him to
advance in the calculation of Gromov-Witten and Welschinger invariants
of del Pezzo surfaces and to treat some del Pezzo surfaces of degrees $3$, $2$, and $1$;
in particular, he proves the non-negativity of Welschinger invariants in one of the two cases
not covered by our Theorem \ref{t9}.}

J.~Solomon \cite{Sol2, HS} suggested a completely different and very powerful recursive tool for computing Welschinger invariants of real blown ups of the projective plane.
His recursion is based on analogues of Kontsevich-Manin axioms and WDVV equation,
and involves the Gromov-Witten invariants and a finite number of initial values.
However, the presence of plenty of terms of opposite signs
(contrary to our formulas which, in most of cases,
contain only non-negative terms)
makes not evident
the use of these recursive formulas
for getting general statements on positivity and asymptotic behavior.

In an unpublished joint
work with R.~Rasdeaconu, J.~Solomon has
considered
a
kind of $w$-numbers
which count
curves subject to point constraints
and
odd
tangency conditions
to a fixed divisor,
and showed that some combinations of such numbers
are independent of point constraints.
Let us underline that
our sided $w$-numbers are defined via
even tangency conditions and,
in some cases, are individually invariant with respect to point constraints.

The
paper is organized as follows.
Section \ref{statements}
describes the main results of the paper;
these results concern positivity and asymptotics of Welschinger invariants.
In Section~\ref{Wel}, we define ordinary and sided $w$-numbers and prove Caporaso-Harris type recursive formulas
for these numbers in the case of real rational surfaces $Y$ with a given real smooth rational
curve $E$ such that the classes $-K_Y$ and $-K_Y-E$
are nef (we call $(Y,E)$ a monic log-del Pezzo pair).
In Section~\ref{sec-vakil}, we consider nodal degenerations of del Pezzo surfaces and derive Abramovich-Bertram-Vakil
type
formulas relating the ordinary and sided $w$-numbers of the degeneration with
Welschinger invariants.
Section~\ref{sec-pa} contains the proof of the main results of the paper.
The
further two sections are devoted to
other applications
of the results of Sections~\ref{Wel} and~\ref{sec-vakil}:
monotonicity behavior of Welschinger invariants is studied in Section~\ref{sec-mon}, and Mikhalkin
type
congruences in Section~\ref{secn22}.
The index of notations can be found at the end of the paper.

\medskip

{\bf Acknowledgments}. A considerable part of the work on this text
was done during
authors' {\it Research in Pairs} stay at the
{\it Mathematisches Forschungsinstitut Oberwolfach}.
The final version of the paper was prepared
during
author's stay at the {\it Max-Planck-Institut f\"ur Mathematik}
in Bonn. We thank these institutions
for hospitality and excellent working conditions.
The first two authors were partially funded by the ANR-09-BLAN-0039-01 grant
of {\it Agence Nationale de la Recherche}.
The third
author enjoyed a support from the Israeli Science Foundation grant
no. 448/09, from the German-Israeli Foundation research grant no.
1174-197.6/2011, and from the Hermann-Minkowski-Minerva Center for
Geometry
at the Tel Aviv University.

We
are deeply grateful to E. Brugall\'e
for the help in uncovering a mistake
in the first draft of this paper.
We are thankful to the referee for the helpful comments and suggestions on the manuscript that allowed us to improve the presentation.

\section{Positivity and asymptotics of Welschinger invariants:
statement of results}\label{statements}

\subsection{Purely real (modified) Welschinger invariants of real del Pezzo
surfaces}\label{intro1}
Let $X$ be a real del Pezzo surface ({\it i.e.},
a real smooth rational surface having an ample anticanonical class
$-K_X$) with a non-empty real point set $\R X$. Let $D\in\Pic(X)$ be
a real effective divisor class with $D^2\ge-1$, assumed to be
primitive in $\Pic(X)$ if $D^2=0$.

The set $R(X,D)$ of
reduced
irreducible rational curves in $|D|$ is a non-empty
quasi-projective variety of pure dimension $-DK_X-1$ with nodal
curves as generic elements (see, for instance, \cite[Lemma
4]{IKS8}). Denote by $\R R(X,D)$ the set of real rational curves in
$R(X,D)$.

We intend to count curves in $\R R(X,D)$
that match
a suitable number
of real point constraints. If $-DK_X>1$, we pick a generic
collection $\bw$ of $-DK_X-1$ points in $\R X$. Since a curve in $\R
R(X,D)$ passing through $\bw$ must contain all these points in its
(unique) real one-dimensional component, we have to suppose that
$\bw$ lies in one connected component of $\R X$. Notice also that if
$-DK_X=1$, each curve in the (finite) set $\R R(X,D)$ has a
one-dimensional real branch.
Indeed, a real curve with a finite real part must have an even
self-intersection, whereas
$D^2
=
-DK_X \mod2$ by the adjunction formula.

To introduce {\it {\rm (}modified{\rm )\/} purely real Welschinger numbers},
let us fix
a connected component $F$ of the real part
$\R X$ of $X$ and,
in addition,
a conjugation
invariant
class $\varphi\in H_2(X\setminus F,\Z/2)$.
If
$-DK_X=1$, we set $\R R(X,D,F)=\{C\in \R R(X,D)\ :\ |C\cap
F|=\infty\}$ and put
$$W(X,D,F,\varphi)=\sum_{C\in\R
R(X,D,F)}(-1)^{s(C)+C_{1/2}\circ\;\varphi}\ ,$$
where $C_{1/2}$ is
the
image of one of the halves of
$\PP^1 \setminus \R P^1$ by the normalization map
$\PP^1\to C$,
and
$s(C)$ is the number of real solitary
nodes of $C$.
If $-DK_X>1$, we pick a generic
collection $\bw$ of $-DK_X-1$ points of $F$, set $\R
R(X,D,\bw)=\{C\in\R R(X,D)\ :\ C\supset\bw\}$, and put
\begin{equation}W(X,D,F,\varphi,\bw)=\sum_{C\in\R
R(X,D,\bw)}(-1)^{s(C)+C_{1/2}\circ\;\varphi}\
.\label{enn1}\end{equation}

The following statement is a version of the
Welschinger theorem \cite{W1} ({\it cf}. also \cite[Theorem
7]{IKS8}).

\begin{theorem}\label{p1}
{\rm (1)} If $-DK_X>1$, the number $W(X,D,F,\varphi,\bw)$
does not depend on the choice of a generic
collection $\bw$ of $-DK_X-1$ points in $F$.

{\rm (2)} With the given data $X,D,F,\varphi$ as above, let $X_t$,
$t\in[0,1]$, $X_0=X$, be a
smooth family of smooth real rational surfaces with non-empty real
part such that for all but finitely many
$t\in[0,1]$,
$X_t$ is a real del Pezzo surface.
Let $\Theta_t:X_0\to X_t$, $t\in[0,1]$, $\Theta_0=\Id$, be a smooth
family of conjugation invariant $C^\infty$-diffeomorphisms
that trivializes our family of surfaces.
Then
$$
W(X,D,F,\varphi,\bw)=W(X_1,(\Theta_{1})_*(D),\Theta_1(F),(\Theta_{1})_*(\varphi),\Theta_1(\bw))\ .
$$
\end{theorem}

In the sequel we write $W(X,D,F,\varphi)$ omitting the notation of point constraints.

\subsection{Real del Pezzo surfaces of degree 2}\label{sec12}
The
classification of real del Pezzo surfaces up to equivariant
deformations
is well known and goes back to A. Comessatti;
the details can be found, for example,
in \cite[Section 17.3]{DIK}.

Here, we recall this classification in the case of del Pezzo surfaces of degree $2$.
According to this classification, a real del Pezzo surface of degree $2$ is determined up to equivariant deformation by the topology of
its real part.

The anticanonical linear system on a real del Pezzo surface $X$ of
degree $2$ defines a double covering $X\to\PP^2$
branched in a nonsingular real quartic curve $Q_X\subset\PP^2$, and
thus identifies $X$ with a hypersurface defined in the weighted
projective space $P^3(1,1,1,2)$ by
an equation $u^2=\eps f_X(x,y,z)$,
where
$f_X$ is a real defining polynomial of $Q_X$ and $\eps=\pm 1$.
Therefore, as a topological space, $\R X$ is
the result of gluing
two copies of $\R f_{X,\eps}=\{p\in\R\PP^2\ :\ \eps
f_X(p)\ge0\}$ along their common boundary, if this boundary is
non-empty,
and the disjoint
union of two copies otherwise.
Below we
always choose the sign for $f_X$ so that $\R f_{X,-}$ is
non-orientable.

As is known, the real part of a real non-singular quartic is
isotopic in $\R P^2$
either to the union of $0\le q\le 4$ null-homologous circles placed
outside each other (denote this isotopy type by $\langle q\rangle$),
or to a pair of null-homologous circles
placed one inside the other (denote this isotopy type by
$1\langle1\rangle$).
In accordance with this notation and the above
sign-convention, the
deformation types of real del
Pezzo surfaces $X$ of degree $2$ with $\R X\ne\emptyset$ are denoted
below by $\langle0\rangle^-$, $\langle q\rangle^\eps$, $1\le q\le4$,
and $1\langle1\rangle^\eps$. For example, the
deformation type
of
the plane blown up at $a$ real points and $b$ pairs of complex
conjugate points, $a+2b=7$, which we denote
by
$\PP^2_{a,b}$,
coincides with $\langle4-b\rangle^-$.

For surfaces $X$ of type
$\langle q \rangle^\eps$,
$1\le q\le4$,
$\langle0\rangle^-$, and $1\langle1\rangle^+$, the choice of a
connected component $F$ of $\R X$
does not affect the computation of Welschinger invariants;  indeed,
for
$X$ of type $\langle0\rangle^-$ the two connected components of $\R
X$
are interchanged by the deck transformation of the above double
covering, while for other types of $X$ with disconnected $\R X$ such
an independence follows from Theorem \ref{p1}(2).

As to surfaces $X$ of type $1\langle1\rangle^-$, they have two
connected components:
one, which we denote $F^{\;o}$, is orientable, and
the other one, $F^{\;no}$, is not.

Notice also that the $28$ bitangents of $Q_X$ lift into
the $56$ curves
in $X$ with self-intersection $-1$, and that the curves of the
linear system $|-K_X|$ are the pull-backs of the straight lines in
$\PP^2$.

\subsection{Main results}\label{sec414}
Let
$X$ be a real
del Pezzo surface.
Denote by $\Pic^\R(X)$ the
subgroup of $\Pic(X)$ formed by real divisor classes of $X$, and
denote
by
$\BH: \Pic^\R(X)\to H_1(\R X; \Z/2)$ the natural
homomorphism which sends
each
real effective divisor class  $D$
that is represented by a real reduced curve, say $C$,
to
$[\R C \cap \R X] \in H_1(\R X, \Z/2)$
(cf.~\cite{BH,Viro}). If $\cal F$ is a union of some connected components of $\R X$,
then denote by
$\BH_{\cal F}$ the composition of $\BH$
with the projection $H_1(\R X; \Z/2) \to H_1({\cal F}; \Z/2)$.

Let $F$ be a connected component of $\R X$. We say that a real effective divisor class
$D$ on $X$ is {\it $F$-compatible},
if $\BH_{\R X \setminus F}(D) = 0$. It is clear that if
a real effective divisor class $D$ is not $F$-compatible,
then
$W(X,D,F,\varphi)$ vanishes
for
any conjugation invariant
class
$\varphi\in H_2(X \setminus F, \Z/2)$.

\begin{theorem}\label{t9}
Let $X$ be a real del Pezzo surface of degree $\ge2$ with a non-empty
real point set,
and let $F$ be a connected component of $\R X$.
Assume
that $F\ne S^2$ if
$X$ is of degree $2$ and $\R X$ is $S^2$ or $S^2\sqcup \R P^2\#\R P^2$.
Then, for any $F$-compatible nef and big divisor class
$D$ on $X$, one has
\begin{equation}
W(X,D,F,[\R X\setminus F])>0\ .\label{ee20}\end{equation} In
particular, through any
collection of $-K_XD-1$ points of
$F$, one can trace a real rational curve $C\in|D|$.
Furthermore,
\begin{equation}
\log W(X,nD,F,[\R X\setminus F]) = -K_X D \cdot n \log n + O(n),
\;\;\; n \to +\infty.\label{ee21}
\end{equation}
\end{theorem}

The proof of Theorem \ref{t9} is presented in Sections
\ref{sec-pos1}
and \ref{sec-pos2}.

\begin{remark}\label{rr2}
{\rm
Theorem \ref{t9}
covers all the cases
studied in \cite{IKS2,IKS4,IKS6,IKS7,Sh2}
(notice that the proof
of
\cite[Theorem 2]{IKS7}
in the journal version
contains an inaccuracy in the case of
the plane blown up at two real points and two pairs of complex conjugate points;
this inaccuracy is corrected in the ArXiv version of \cite{IKS7}).}
\end{remark}

\begin{theorem}\label{t11}
Let $X$ be a real del Pezzo surface of degree $2$ with $\R X = S^2$.
Then,
\begin{enumerate}\item[{\rm (i)}] for any real effective divisor class
$D$ on $X$, we have $W(X,D,\R X,0)\ge0$;
\item[{\rm (ii)}] the big
and nef real effective divisor classes $D$ on $X$ such that $W(X,D,\R X,0)>0$
form a subsemigroup in $\Pic(X)$; this subsemigroup contains
$-mK_X$ with $m \geq 2$ and
all the divisor classes $D'$ and $D' - K_X$, where $D'$ is big, nef,
and disjoint from a pair of complex conjugate $(-1)$-curves;
\item[{\rm (iii)}] if a big and nef real effective divisor class $D$ on $X$
satisfies $W(X,D,\R X,0)>0$, then
\begin{equation}
\log W(X,nD,\R X,0) = -K_X D \cdot n \log n + O(n), \;\;\; n \to
+\infty ;\label{e2078}\end{equation}
\item[{\rm (iv)}] if a big and nef real effective divisor class
$D$ on $X$ satisfies $D^2\le2$, then $W(X,D,\R X,0)=0$ as long
as $-K_XD\ne4$.
\end{enumerate}
\end{theorem}

The proof is given in Section \ref{sec-s2}.

Notice that Theorem
\ref{t11} (iv) implies the following statement: for a real Del Pezzo
surface $X$ of degree $2$ with $\R X = S^2$ there are infinitely
many nef and big real divisors $D$  such that $W(X,D,\R X,0)=0$.
Indeed, represent $X$ as an ellipsoid (that is, a real quadric with spherical real part)
blown up at $3$ pairs of
complex conjugate points and choose
the
basis $L_1,L_2,E_1,...,E_6$ of
$\Pic(X)$
where $L_1,L_2$ are generators of the
quadric
and
$E_1,...,E_6$ are the exceptional divisors of the blow up; then,
each divisor $D=m(L_1+L_2)-n(E_1+...+E_6)$, where $m^2-3n^2=1$ and
$m\ne7$, is real and nef, and it satisfies $D^2=2$ and
$-K_XD=4m-6n\ne 4$.

Such a vanishing is sometimes "sharp":
if the only oval of a real plane
quartic of type $\langle1\rangle$ is convex, then there is no real
tangent through a point inside the oval, and hence there are no real rational curves $C\in|-K_X|$
at all.

The
case $\R X = S^2\sqcup \R P^2\#\R P^2$ and $F = S^2$, not covered by the above theorems,
is discussed in Section \ref{prooft12}.

The following table contains the values of Welschinger invariants $W(X,
D, F,\varphi)$ for $D=-K_X$ or $-2K_X$ and $\varphi=0$ or $\varphi=\varphi_F
=[\R X\setminus F]$
(for surfaces of types $\langle 2 \rangle^+$, $\langle 3 \rangle^+$
or $\langle 4 \rangle^+$, the invariants do not depend on the choice of $F$ among the components
of $\R X$).

{\scriptsize{
\begin{center}
\begin{tabular}{|l|c||c|c|c|c|c|c|c|c|c|c|c|c|c|}
\hline $D$ & $\varphi$ & $\langle 4 \rangle^-$ & $\langle 3
\rangle^-$ & $\langle 2 \rangle^-$ & $\langle 1 \rangle^-$ & $\langle 0\rangle^-$ &
$1\langle1\rangle^+$ & $1\langle1\rangle^-,F^{\;no}$ & $1\langle1\rangle^-,F^{\;o}$ &
$\langle1\rangle^+$ & $\langle2\rangle^+$ & $\langle3\rangle^+$ &
$\langle4\rangle^+$
  \\
\hline\hline $-K_X$ & $0$ & $8$ & $6$ & $4$ & $2$ & $0$ & $2$ & $0$ & $0$ &
$0$ & $-2$ & $-4$ & $-6$
  \\
\hline $-K_X$ & $\varphi_F$ & $8$ & $6$ & $4$ & $2$ & $0$ & $2$ & $4$ & $0$ & $0$
& $2$ & $4$ & $6$ \\ \hline
\hline $-2K_X$ & $0$ & $224$ & $128$ & $64$ & $24$ & $0$ & $32$ & $0$ & $0$ &
$8$ & $0$ & $0$ & $0$
  \\
\hline
$-2K_X$ & $\varphi_F$ & $224$ & $128$ & $64$ & $24$ & $32$ & $32$ & $48$ & $16$ &
$8$ & $16$ & $32$ & $64$
  \\
\hline
\end{tabular}
\end{center}
}}

Notice
that the original Welschinger invariants ($\varphi=0$) may take negative values or vanish
for the multi-component del Pezzo surfaces. This reflects the following general phenomenon.

\begin{proposition}\label{t10}
Let $X$ be a real del Pezzo
surface of degree $\ge2$
with disconnected real point set,
let $F$ and $F'$ be two distinct connected components of $\R X$,
and let $\varphi \in H_2(X \setminus (F \cup F'); \Z/2)$ be a conjugation invariant class.
Then,
$W(X, D, F, \varphi) = 0$
for any big and nef real effective divisor class
$D$ on $X$ such that $-K_XD\ge3$.
\end{proposition}

This proposition immediately follows from the formulas (\ref{ee3}) and (\ref{ee13})
proved, respectively, in Sections \ref{sec7} and \ref{sec9II}.
In a more general setting,
the
vanishing statement given by
Proposition \ref{t10}
is found in \cite[Proposition 3.3]{BP}.

\section{Recursive formulas for $w$-numbers of real monic log-del Pezzo pairs}\label{Wel}

\subsection{Surfaces under consideration}\label{log-nodal}
Let $Y$ be a smooth rational surface which is a blow-up of $\PP^2$,
and let $E\subset Y$ be a smooth rational curve.
Suppose that $-K_Y $ is positive on all curves different from $E$
and $K_YE\ge0$, and that the log-anticanonical class $-(K_Y+E)$ is
nef, effective, and satisfies
$(K_Y +E)^2=0$. We call such a pair $(Y,E)$ a {\it monic log-del Pezzo
pair}. Throughout Section~\ref{Wel}, we assume that $(Y,E)$ is a
monic log-del Pezzo pair.

Observe that  $-(K_Y+E)E=2$, $E^2\le-2$, and
$K_Y(K_Y+E)=2$,
so that the latter
implies, once more by adjunction, that $|-(K_Y+E)|$ is a
one-dimensional linear system, whose generic element is a smooth
rational curve.
This linear system contains precisely two smooth curves $L'$, $L''$
(quadratically) tangent to $E$, and $4-E^2$ reducible curves, all of
type $L_1+L_2$ where $L_1^2=L_2^2=-1$, $L_1L_2=1$, $L_1E=L_2E=1$.
In particular, it provides a conic bundle structure on $Y$
and shows that $Y$ can be regarded as the plane blown up at $\ge 6$
points on a smooth conic ($E$ is the strict transform of the conic)
and at one more point outside the conic. We will assume that the
blown up points are in general position subject to the above
allocation with respect to the conic.
The curves $L'$ and $L''$ are called {\it supporting curves}.

Introduce the sets
$$
{\mathcal E}(E)=\{E'\in\Pic(Y)\ :\ (E')^2=-1, \ E'K_Y = -1, \ E'E>0\}\ .
$$
$$
{\mathcal E}(E)^{\perp D} =\{E'\in{\mathcal E}(E)\ :\ E'D=0\},\quad
D\in\Pic(Y)\ .
$$

Suppose that $(Y,E)$ is equipped with a real structure such that $\R
Y\supset\R E\ne\emptyset$. Denote by $F$ the connected component of
$\R Y$ containing $\R E$. We also choose a
conjugation invariant
class $\varphi\in H_2(Y\setminus F,\Z/2)$.

Quadruples $(Y,E,F,\varphi)$ as above are called {\it basic
quadruples}.

\subsection{Some notations}\label{some_notations}
Let $\Z^\infty_+$ be the direct sum of countably many additive
semigroups \mbox{$\Z_+=\{
k \in\Z\ |\
k \ge 1\}$}, labeled by the
positive integer numbers, with the basis
formed by the summand generators $e_i$, $i=1,2,...$ For
$\alp=(\alp_1,\alp_2,...)\in\Z^\infty_+$, put
$$
\|\alp\|=\sum_{
i=1}^\infty
\alp_i, \quad I\alp=\sum_{
i=1}^\infty
i \alp_i, \quad I^\alp= \prod_{i=1}^\infty
i^{\alp_i},\quad\alp!=\prod_{i=1}^\infty\alp_i!\ .
$$
For $\alp, \beta \in \Z^\infty_+$, we write $\alp \geq \beta$ if
$\alpha_i \geq \beta_i$ for any positive integer number~$i$.
For $\alp^{(0)},...,\alp^{(m)},\alp\in\Z_+^\infty$ such that
$\alp^{(0)}+...+\alp^{(m)}\le\alp$, put
$$\left(\begin{matrix}\alp\\
\alp^{(0)},...,\alp^{(m)}\end{matrix}\right)=
\frac{\alp!}{\alp^{(0)}!...\alp^{(m)}!(\alp-\alp^{(0)}-...-\alp^{(m)})!}\
.$$ Introduce also the semigroups
\begin{eqnarray}\Z^{\infty,\;\odd}_+&=&\Span\{e_{2i+1}\ : \ i\ge0\}\
,\nonumber\\
\Z^{\infty,\;\even}_+&=&\Span\{e_{2i}\ :\ i\ge1\}\ ,\nonumber\\
\Z^{\infty,\;\odd\;\cdot\;\even}_+&=&\Span\{e_{4i+2}\ :\ i\ge0\}\
.\nonumber\end{eqnarray}

\subsection{Divisor classes}\label{div-cl}
Let $\Sig$ be a smooth real surface.
We denote by $\Pic^\R(\Sig)$ the
subgroup of $\Pic(\Sig)$ formed by real divisor classes of $\Sig$ and denote by $\PicPR(\Sig)$ the
subsemigroup of $\Pic^\R(\Sig)$ generated by
effective real divisor classes. Let $E\subset\Sig$ be a smooth real curve.
Put $\PicPP(\Sig,E)$ to be
the subsemigroup of $\Pic(\Sig)$ generated by complex irreducible
curves $C$ such that $CE\ge0$.
The involution of complex
conjugation~$\conj: \Sig \to \Sig$
naturally acts on $\Pic(\Sig)$
and preserves
$\PicPP(\Sig, E)$. Denote by
$\PicPPR(\Sig,E)$ the disjoint union of the sets
$$
\left\{D\in\PicPP(\Sig,E)\ :\ \conj \ D=D\right\}
$$
and
$$
\left\{\{D_1,D_2\}\in\Sym^2(\PicPP(\Sig,E))\ :\ \conj \ D_1
=D_2\right\}\ .
$$
For an element $\D\in\PicPPR(\Sig,E)$, define $[\D]\in\PicPP(\Sig,E)$ by
$$[\D]=\begin{cases}D,\quad & \D=D,\ \text{a divisor class},\\
D_1+D_2,\quad &\D=\{D_1,D_2\},\ \text{a pair of divisor
classes}.\end{cases}$$ For a element $\D\in\PicPPR(\Sig,E)$ and a
vector $\bet\in\Z^\infty_+$, put
$$R_\Sig(\D,\bet)=-[\D](K_\Sig+E)+\|\bet\|-\begin{cases}1,\quad & \D=D,\ \text{a divisor class},\\
2,\quad &\D=\{D_1,D_2\},\ \text{a pair of divisor
classes}.\end{cases}$$

\subsection{Families of real curves}\label{Welsch-inv}
Let $(Y,E)$ be a real monic log-del Pezzo pair.
An {\it admissible tuple} $(\D, \alp, \bet^{\re}, \bet^{\ima},
\bp^\flat)$ consists of an element $\D$ in $\PicPPR(Y,E)$,
vectors $\alp,\bet^{\re},\bet^{\ima}\in\Z^\infty_+$ satisfying
\mbox{$I(\alp+\bet^{\re}+2\bet^{\ima})=[\D]E$}, and a sequence
$\bp^\flat=\{p_{i,j}\ :\ i\ge1,\ 1\le j\le\alp_i\}$
of $\|\alp\|$ distinct real generic points on $E$.
Denote by $V^\R_Y (\D,\alp,\bet^{\re},\bet^{\ima},\bp^\flat)$
the closure in
the linear system $|[\D]|$ of the family of real
reduced curves~$C$ such that
\begin{enumerate}\item[(i)] if $\D=D$, a divisor class,
then $C\in|D|$ is an irreducible over $\C$ rational curve; if
$\D=\{D_1,D_2\}$, a pair of divisor classes, then $C=C_1\cup C_2$,
where $C_1\in|D_1|$, $C_2\in|D_2|$ are distinct, irreducible,
rational, complex conjugate curves; \item[(ii)] $C\cap E$ consists
of $\bp^\flat$ and of $\|\bet^{\re}+2\bet^{\ima}\|$ other points:
$\|\bet^{\re}\|$ of them real, and $2\|\bet^{\ima}\|$ form pairs of
complex conjugate points;
\item[(iii)] at each point of $C\cap E$,
the curve $C$ has one local branch, and the intersection
multiplicities of $C$ and $E$ are described as follows:
\begin{itemize}\item $(C\cdot E)(p_{i,j})=i$ for all $i\ge 1$, $1\le
j\le\alp_i$, \item for each $i\ge 1$, there are $\bet^{\re}_i$ real
points $q\in(C\cap E)\setminus\bp^\flat$ such that $(C\cdot
E)(q)=i$; \item for each $i\ge 1$ there are $\bet^{\ima}_i$ pairs
$q,q'$ of complex conjugate points of $C\cap E$ such that $(C\cdot
E)(q)=(C\cdot E)(q')=i$.
\end{itemize}
\end{enumerate}

If $D \in \PicPPR(Y, E)$ a divisor class,
introduce also the variety
$V_Y(D,\alp,\bet,\bp^\flat)$ which is the closure in $|D|$ of the
family of complex
reduced
irreducible rational curves~$C$ such that
$C \cap E$ consists of $\bp^\flat$ and of $\|\bet\|$ other points,
at each point of $C\cap E$, the curve $C$ has one local branch, and
the intersection multiplicities of $C$ and $E$ are as follows:
\begin{itemize}
\item $(C\cdot E)(p_{i,j})=i$ for all $i\ge 1$, $1\le j\le\alp_i$,
\item for each $i\ge 1$, there are
$\bet_i$ points $q\in(C\cap E)\setminus\bp^\flat$ such that
$(C \cdot E)(q) = i$.
\end{itemize}

\begin{lemma}\label{l1}
If $\D = D$ is a divisor class
and
$V^\R_Y(\D,\alp,\bet^{\re},\bet^{\ima},\bp^\flat)$ is nonempty,
then $R_Y(\D,\bet^{\re}+2\bet^{\ima})\ge0$, and each component of
$V^\R_Y(\D,
\alp,
\bet^{\re},\bet^{\ima},\bp^\flat)$
has dimension $\le R_Y(\D,\bet^{\re}+2\bet^{\ima})$.
Moreover, a generic element~$C$ of any
component of $V^\R_Y(\D,\alp,\bet^{\re},\bet^{\ima},\bp^\flat)$ of
dimension $R_Y(\D,\bet^{\re}+2\bet^{\ima})$ is an immersed curve,
nonsingular along $E$. If, in addition, $E^2 \geq -3$, then~$C$ is nodal.
\end{lemma}

{\bf Proof}. If $D$ is a multiple of a divisor class orthogonal to $K_Y + E$,
then $V^\R_Y(\D,\alp,\bet^{\re},\bet^{\ima},\bp^\flat)$ cannot be nonempty,
since such a linear system contains only non-reduced curves.
In the other case, the statement follows from~\cite[Proposition 2.1]{MS2}.
\proofend

\smallskip

Suppose that $R_Y(\D,\bet^{\re}+2\bet^{\ima})\ge 0$. Pick a set
$\bp^\sharp$ of $R_Y(\D,\bet^{\re}+2\bet^{\ima})$ generic points of
$F\setminus E$ and denote by $V^\R_Y
(\D,\alp,\bet^{\re},\bet^{\ima},\bp^\flat,\bp^\sharp)$
the set of
curves $C\in V^\R_Y(\D,\alp,\bet^{\re},\bet^{\ima},\bp^\flat)$
passing through $\bp^\sharp$.

\begin{lemma}\label{l2}
Assume that $V^{\R}_Y(\D, \alp, \beta^{\re}, \beta^{\ima},
\bp^\flat)$ is nonempty.

{\rm (1)}
If $\D=D$, a divisor class, then $V^\R_Y
(D,\alp,\bet^{\re},\bet^{\ima},\bp^\flat,\bp^\sharp)$ is a finite
set of real
immersed
irreducible
rational curves
which are nonsingular along $E$.

{\rm (2)}
If $\D=\{D_1,D_2\}$, a pair of divisor classes, then $V^\R_Y
(\D,\alp,\bet^{\re},\bet^{\ima},\bp^\flat,\bp^\sharp)$ is finite,
and it is nonempty only if $\alp=\bet^{\re}=0$, $R_Y
(\D,2\bet^{\ima})=0$, and $\bp^\flat=\bp^\sharp=\emptyset$.
\end{lemma}

{\bf Proof}. By Lemma \ref{l1} we have to show only that $R_Y
(\D,2\bet^{\ima})=0$ is necessary for the nonemptyness of $V^\R_Y
(\D,0,0,\bet^{\ima},\emptyset,\bp^\sharp)$ with $\D=\{D_1,D_2\}$,
and the proof of this fact literally coincides with the proof of
\cite[Lemma 3(2)]{IKS7}. \proofend

\begin{lemma}\label{l3}
The only nonempty sets
$V^\R_Y(\D,\alp,\bet^{\re},\bet^{\ima},\bp^\flat)$
for admissible tuples $(\D,\alp,\bet^\re,\bet^\ima,\bp^\flat)$ such
that
$ R_Y(\D,\bet^{re}+2\bet^{\ima})=0$ are the following ones:
\begin{enumerate}
\item[{\rm (1)}] if $\D=D$ is
a divisor class and
$I(\alp+\bet^{\re}+2\bet^{\ima})=DE>0$,
\begin{enumerate}
\item[{\rm (1i)}]
$V^\R_Y(E',0,e_1,0,\emptyset)$ consists of one element, where $E'$
is a real $(-1)$-curve crossing $E$;
\item[{\rm (1ii)}] $V^\R_Y
(-(K_Y+E),0,e_2,0,\emptyset)$ consists of
two elements $L',L''$, if
the supporting
curves
$L',L''$ are both real;
\item[{\rm (1iii)}] $V^\R_Y(-(K_Y
+E),e_1,e_1,0,\bp^\flat)$ consists of one element;
\item[{\rm (1iv)}]
$V^\R_Y(D,\alp,0,0,\bp)$ consists of one element, if $(K_Y+E)D=-1$,
$I\alp=DE$;
\end{enumerate}
\item[{\rm (2)}] if $\D=\{D_1,D_2\}$
is
a pair of divisor classes
and $I(\alp+\bet^{\re}+2\bet^{\ima})=[\D]E>0$,
\begin{enumerate}
\item[{\rm (2i)}] $V^\R_Y(\{E'_1,E'_2\},0,0,e_1,\emptyset)$ consists
of one element, where $E'_1,E'_2$ are complex conjugate
$(-1)$-curves crossing $E$;
\item[{\rm (2ii)}] $V^\R_Y(\{-(K_Y+E),-(K_Y
+E)\},0,0,e_2,\emptyset)$ consists of one element $\{L',L''\}$, if
$L',L''$ are complex conjugate;
\end{enumerate}
\item[{\rm (3)}] if $I(\alp+\bet^{\re}+2\bet^{\ima})=[\D]E=0$,
\begin{enumerate}
\item[{\rm (3i)}]
$V^\R_Y(E',0,0,0,\emptyset)$ consists of one element, where $E'$ of
a real $(-1)$-curve disjoint from $E$;
\item[{\rm (3ii)}]
$V^\R_Y(\{E'_1,E'_2\},0,0,0,\emptyset)$ consists of one element,
where $E'_1,E'_2$ are complex conjugate $(-1)$-curves disjoint from
$E$.
\end{enumerate}\end{enumerate}
\end{lemma}

{\bf Proof}. Straightforward.
\proofend

\subsection{Deformation diagrams and CH position}\label{sec5}

\subsubsection{Deformation
diagrams}\label{deform-diagrams}
Let $(Y, E)$ be a monic log-del
Pezzo pair such that $Y$ and $E$ are real and $\R E \ne
\emptyset$.
Denote by $F$ the connected component of $\R Y$
containing $\R E$ and pick a conjugation invariant
class $\varphi\in
H_2(\R Y\setminus F,\Z/2)$.
Let
$(D,\alp,\bet^{\re},\bet^{\ima},\bp^\flat)$ be an admissible tuple,
where $D \in \PicPPR(Y,E)$
is a divisor class and
$R_Y(D,\bet^{\re}+2\bet^{\ima})>0$. Pick a set
$\widetilde\bp^\sharp$ of $R_Y(D,\bet^{\re}+2\bet^{\ima})-1$ generic
real points of $F\setminus E$, a generic real
point $p\in
E\setminus\bp^\flat$, and a smooth real algebraic curve germ $\Lam$,
crossing $E$ transversally at $p$. Denote by $\Lam^+=\{p(t)\ :\
t\in(0,\eps)\}$ a parameterized connected component of
$\Lam\setminus\{p\}$ with $\lim_{t\to0}p(t)=p$. There exists
$\eps_0>0$ such that, for all $t\in(0,\eps_0]$, the sets $V_Y
(D,\alp,\bet^{\re},\bet^{\ima},\bp^\flat,\widetilde\bp^\sharp\cup\{p(t)\})$
are finite, their elements remain immersed, nonsingular along $E$ as
$t$ runs over the interval $(0,\eps_0]$, and the closure in $V_Y
(D,\alp,\bet^{\re},\bet^{\ima},\bp^\flat)$ of the family
\begin{equation}V=\bigcup_{t\in(0,\eps_0]}V_Y(D,\alp,\bet^{\re},
\bet^{\ima},\bp^\flat,
\widetilde\bp^\sharp\cup\{p(t)\})\label{eDD}\end{equation} is a
union of real algebraic arcs, disjoint for $t>0$. This closure is
called a \emph{deformation diagram} of
$(D,\alp,\bet^{\re},\bet^{\ima},\bp^\flat,\widetilde\bp^\sharp,p)$,
{\it cf}. \cite[Section 3.3]{IKS7},
and the real algebraic arcs under consideration
are called {\it branches} of the deformation diagram.
The elements of $V_Y
(D,\alp,\bet^\re,\bet^\ima,\bp^\flat,\widetilde\bp^\sharp\cup\{p(1)\})$
are called \emph{leaves} of the deformation diagram, and the
elements of $\overline V\setminus V$ are called \emph{roots} of the
deformation diagram.

\begin{lemma}\label{DC_new}
Each connected component of a deformation diagram
of $(D,\alp,\bet^{\re},\bet^{\ima},\bp^\flat,\widetilde\bp^\sharp,p)$
contains exactly one root.
Each root is
either a generic member of an
$(R_Y(D,\bet^{\re}+2\bet^{\ima})-1)$-dimensional component of one
of the families
$$V_Y(D,\alp+e_j,\bet^{\re}-e_j,\bet^{\ima},
\bp^\flat\cup\{p\}, \widetilde\bp^\sharp\})\ ,$$ where $j$ is a
natural number such that $\bet^{\re}_j>0$, or a reducible
curve having $E$ as a component.
\end{lemma}

{\bf Proof.}
The statement follows from \cite[Proposition 2.6]{MS2}.
\proofend

For
any root~$\rho$ of a deformation diagram,
the leaves belonging to the connected component of~$\rho$
is said to be {\it generated} by~$\rho$.

\subsubsection{CH position}\label{CH-position}
Pick a divisor class $D_0\in\PicPPR(Y,E)$ and put $N=\dim|D_0|$.
Note that the set
$$\Prec(D_0)=\{D \in \PicPPR(Y,E)
\ : \ D \ \text{\rm is a divisor class and} \ D_0 \ge D\}$$
is
finite, and we have $\dim|D|\le N$ for each $D\in\Prec(D_0)$.
Furthermore, for each nonempty variety
$V^\R_Y(D,\alp,\bet^{\re},\bet^{\ima},\bp^\flat)$ with
$D\in\Prec(D_0)$, we have
$$\|\alp\|+R_Y(D,\bet^{\re}+2\bet^{\ima})\le N\ .$$

\begin{lemma}\label{D0}{\rm ({\it cf}. \cite[Lemma 10]{IKS7})}
Let $D_0\in \PicPPR(Y,E)$ be a divisor
class with $N=\dim|D_0|>0$. Then, there exists a sequence
$\Lam(D_0)=(\Lam_i)_{i=1,...,N}$ of $N$ disjoint smooth real
algebraic arcs in $Y$, which are parameterized by $t \in
[-1, 1] \mapsto p_i(t)\in\Lam_i$, such that $p_i(0)\in E$, the arcs
$\Lam_i$ are transverse to $E$ at $p_i(0)$, $i=1,...,N$, and the
following condition holds:

for
any admissible tuple $(D,\alp,\bet^{\re},\bet^{\ima},\bp^\flat)$,
any disjoint subsets $J^\flat,J^\sharp\subset\{1,...,N\}$, any
positive integer $k\le N$, and any sequence $\sigma = (\sigma_i)_{i
= 1, \ldots, N}$ such that
\begin{enumerate}
\item[{\rm (i)}]
$D \in \Prec(D_0)$,
\item[{\rm (ii)}] $R_Y(D,\bet^{\re}+2\bet^{\ima})>0$,
\item[{\rm (iii)}] $i<k<j$ for
all $i\in J^\flat$, $j\in J^\sharp$,
\item[{\rm (iv)}] the number of
elements in $J^\sharp$ is equal to $R_Y
(D,\bet^{\re}+2\bet^{\ima})-1$,
\item[{\rm (v)}]
$\bp^\flat=\{p_i(0)\ :\ i\in J^\flat\}$,
\item[{\rm (vi)}] $\sigma_i = \pm1$ for any integer $1 \leq i \leq N$,
\end{enumerate}
the closure of the family
$$
\bigcup_{t\in(0,1]}V_Y(D,\alp,\bet^{\re},\bet^{\ima},\bp^\flat,
\widetilde\bp^\sharp\cup\{p_k(\sigma_k t)\})\ ,
$$
where $\widetilde\bp^\sharp = \{p_j(
\sigma_j)\}_{j\in J^\sharp}$, is a deformation diagram of
$(D,\alp,\bet^{\re},\bet^{\ima},\bp^\flat,\widetilde\bp^\sharp,p_k(0))$.
\end{lemma}

{\bf Proof}. Take a sequence $\widehat\Lam_i$, $i=1,...,N$, of
disjoint smooth real algebraic arcs in $Y$, which are parameterized
by $t \in [-1,1] \mapsto p_i(t)\in\Lam_i$, such that
$(p_i(0))_{i=1,...,N}$ is a generic sequence of points in $E$, and
the arcs $\widehat\Lam_i$ are transverse to $E$ at $p_i(0)$,
$i=1,...,N$. We will inductively shorten
these arcs in order to satisfy the
condition
required in Lemma.

Take an integer $1 \leq k \leq N$, and suppose that we
have
already constructed
arcs $\Lam_1,...,\Lam_{k-1}$
parameterized respectively by intervals $[0, \eps_i]$, $1 \leq i <
k$. There are finitely many admissible tuples
$(D,\alp,\bet^{\re},\bet^{\ima},\bp^\flat)$, subsets
$J^\flat,J^\sharp\subset\{1,...,N\}$, and sequences $\sigma$
satisfying the restrictions (i)-(vi)
above. Given such a
datum $D,\alp,\bet^{\re},\bet^{\ima},\bp^\flat,J^\flat,J^\sharp,
\sigma$, we take a small positive number $\eps_k$
such that the closure of the family
$$\bigcup_{t \in (0, \eps_k]} V_Y
(D,\alp,\bet^{\re},\bet^{\ima},\bp^\flat,\widetilde\bp^\sharp \cup\{
p_k(\sigma_k t)\})\ ,
$$
where $\widetilde\bp^\sharp=\{p_i(\eps_i)\}_{1\le i<k}$, is a
deformation diagram of
$(D,\alp,\bet^{\re},\bet^{\ima},\bp^\flat,\widetilde\bp^\sharp,p_k(0))$,
and put
$$
\Lam_k(D, \alp, \bet^{\re}, \bet^{\ima}, \bp^\flat, J^\flat,
J^\sharp, \sigma) = \bigcup_{t \in [-\eps_k, \eps_k]} p_k(t).
$$
Then, we define
$$
\Lam_k=\bigcap_{(D, \alp, \bet^{\re}, \bet^{\ima}, \bp^\flat,
J^\flat, J^\sharp, \sigma)}
\Lam_k(D,\alp,\bet^{\re},\bet^{\ima},\bp^\flat, J^\flat,J^\sharp,
\sigma).
$$

It remains now to reparameterize by the interval $[-1, 1]$ the arcs
$\Lam_1$, $\ldots$, $\Lam_N$ obtained. \proofend

Take a divisor class
$D_0\in\PicPPR(Y,E)$ such that $N = \dim |D_0| > 0$ and a sequence
of arcs $(\Lam_i)_{i=1,...,N}$ as in Lemma \ref{D0}. Given a
sequence $\sigma = (\sigma_i)_{i = 1, \ldots, N}$ of $\pm 1$ and two
subsets $J^\flat,J^\sharp\subset\{1,...,N\}$ such that $i<j$ for all
$i\in J^\flat$, $j\in J^\sharp$, we say that the pair of point
sequences
$$\bp^\flat=\{p_i(0)\ :\ i\in J^\flat\},\quad
\bp^\sharp=\{p_j(\sigma_j)\ :\ j\in J^\sharp\}$$
is \emph{in
$D_0$-CH position}.
A
pair of point sequences
$$(\bp^\flat)' = \{p_i(0)\ :\ i\in (J^\flat)'\},\quad
(\bp^\sharp)' = \{p_j(\sigma_j)\ :\ j\in (J^\sharp)'\}$$
in $D_0$-CH position
is said to be a {\it predecessor} of a pair of point sequences
$$\bp^\flat=\{p_i(0)\ :\ i\in J^\flat\},\quad
\bp^\sharp=\{p_j(\sigma_j)\ :\ j\in J^\sharp\}$$
in $D_0$-CH position
if $(J^\sharp)' = \{j \in J^\sharp : j > k\}$
for a certain integer~$k$.

Let $(\D,\alp,\bet^{\re},\bet^{\ima},\bp^\flat)$ be an admissible tuple
such that $\D = D$
is a divisor class.
Choose a sequence $\bp^\sharp$
of $R_Y(\D,\bet^{\re}+2\bet^{\ima})$ points in~$F$,
and assume that
the pair of point sequences $\bp^\flat,\bp^\sharp$ is in a $D_0$-CH
position. Then, $(D, \alp, \bet^{\re}, \bet^{\ima}, \bp^\flat, \bp^\sharp)$
is called a {\it $D_0$-proper tuple}.
The elements of
$V_Y^\R(D,\alpha,\bet^{\re},\bet^{\ima},\bp^\flat, \bp^\sharp)$
are called {\it interpolating curves} constrained
by the $D_0$-proper tuple
$(D, \alp, \bet^{\re}, \bet^{\ima}, \bp^\flat, \bp^\sharp)$.
We say that a $D_0$-proper tuple
$(D', \alp', (\bet^{\re})', (\bet^{\ima})', (\bp^\flat)', (\bp^\sharp)')$
{\it precedes} a $D_0$-proper tuple
$(D, \alp, \bet^{\re}, \bet^{\ima}, \bp^\flat, \bp^{\sharp})$
if $R_Y(D', (\bet^{\re})' + 2(\bet^{\ima})')
< R_Y(D,\bet^{\re}+2\bet^{\ima})$
and the pair $(\bp^\flat)', (\bp^\sharp)'$
is a predecessor of $\bp^\flat, \bp^\sharp$.

\begin{lemma}\label{imaginary}
Let $D_0\in\PicPPR(Y,E)$ be a divisor class,
and let
$(D, \alp, \bet^{\re}, \bet^{\ima}, \bp^\flat, \bp^\sharp)$
be a $D_0$-proper tuple such that
$R_Y(D, \bet^{\re} + 2\bet^{\ima}) > 0$ and $\bet^{\ima} \ne 0$.
Then,
$V_Y^\R(D,\alpha,\bet^{\re},\bet^{\ima},\bp^\flat,
\bp^\sharp) = \emptyset$.
\end{lemma}

{\bf Proof}.
Assume that $V_Y^\R(D,\alpha,\bet^{\re},\bet^{\ima}
\bp^\flat,\bp^\sharp)\ne\emptyset$,
and put
$k =
\min J^\sharp$, where $\bp^\sharp=\{p_j(\sigma_j)\ :\ j\in J^\sharp\}$
(see Lemma~\ref{D0}).
We obtain inductively a contradiction
showing that $V_Y^\R(D',\alpha',(\bet^{\re})',(\bet^{\ima})'
(\bp^\flat)',(\bp^\sharp)')\ne\emptyset$
for a certain $D_0$-proper tuple
$(D',\alp',(\bet^{\re})',(\bet^{\ima})',(\bp^\flat)',(\bp^\sharp)')$
that precedes
$(D, \alp, \bet^{\re}, \bet^{\ima}, \bp^\flat, \bp^\sharp)$
and such that $R_Y(D', (\bet^{\re})' + 2(\bet^{\ima})') > 0$ and
$(\bet^{\ima})' \ne 0$.

Consider the degeneration of
$C\in V_Y^\R(D,\alpha,\bet^{\re},\bet^{\ima},\bp^\flat,\bp^\sharp)$
when $p_k\in\bp^\sharp$ tends to $E$ along the arc $\Lam_k$.
By
\cite[Proposition 2.6]{MS2}, the degenerate curve
is
either
an
irreducible interpolating curve constrained by
a $D_0$-proper tuple
that precedes
$(D, \alp, \bet^{\re}, \bet^{\ima}, \bp^\flat, \bp^\sharp)$,
or
of the form
$E\cup C'$.
In the latter case,
the curve $C'$
has a real component
belonging to
$V_Y^\R(D',\alp',(\bet^{\re})',(\bet^{\ima})',(\bp^\flat)',(\bp^\sharp)')$,
where
$(D',\alp',(\bet^{\re})',(\bet^{\ima})',(\bp^\flat)',(\bp^\sharp)')$
is a $D_0$-proper tuple which precedes
$(D, \alp, \bet^{\re}, \bet^{\ima}, \bp^\flat, \bp^\sharp)$
and $R_Y(D', (\bet^{\re})' + 2(\bet^{\ima})') > 0$, $(\bet^{\ima})' \ne 0$..
This statement follows
from~\cite[Lemma 2.9]{MS2} and the fact that
any imaginary
component of $C'$
avoids $\bp^\sharp$, and thus
has a unique
intersection point with $E$ (see Lemma~\ref{l3}).
The former lemma states that
the intersection points of $C$ with $E\setminus\bp^\flat$ all
come from the intersection points of $C'$ with
$E\setminus\bp^\flat$, and
that, in the deformation of
$E\cup C'$ into $C$,
each component of $C'$ glues up with $E$ via
smoothing out one of its intersection points with
$E\setminus\bp^\flat$.
\proofend

\subsection{Ordinary $w$-numbers}\label{new-ordinary}

Let $C$ be a real curve on a real smooth surface $\Sig$,
and let $z$ be a real singular point of $C$
such that all local branches of $C$ at $z$ are smooth.
Denote by $s(C, z)$
the number of pairs
of imaginary complex conjugate local branches of $C$ at $z$,
each pair being counted with the weight equal to the intersection number of the branches.

\begin{lemma}\label{ordinary_invariance}
Let $C(t)$, $-\varepsilon < t < \varepsilon$, be a continuous family of real curves in~$\Sig$,
and let~$z_0$ be a real singular point of $C(0)$
such that all local branches of $C(0)$ at $z_0$ are smooth.
Assume that for a certain neighborhood $U(z_0) \subset \Sig$ of~$z_0$ and a sufficiently
small number $\varepsilon' > 0$, the curves
$C(t)$,
$-\varepsilon' < t < \varepsilon'$,
are transversal to the boundary of $U(z_0)$, and
the curves
$C(t) \cap U(z_0)$, $-\varepsilon' < t < \varepsilon'$,
admit simultaneous parametrizations by a continuous family of
immersions
$\Del_i(t) \to U(z_0)$, $i = 1$, $\ldots$, $b(z_0)$, where
$b(z_0)$ is the number of local branches of $C(0)$ at $z(0)$, and
$\Del_i(t)$, $-\varepsilon' < t < \varepsilon'$, is a continuous family of discs in~$\C$.
Then, $\sum_{z \in \Sing(C(t)) \cap U(z_0)}s(C(t), z)$
does not depend on~$t$.
\end{lemma}

{\bf Proof}. Straightforward.
\proofend

For an immersed real curve $C \subset \Sig$,
put
$s(C) = \sum_{z\in\Sing(C)}s(C,z)$.

Let $(\D, \alp, \bet^{\re}, \bet^{\ima}, \bp^\flat)$ be an admissible tuple such that
$\D = D\in\PicPPR(Y,E)$ is a divisor class, and
let $\bp^\sharp$ be a generic set of $R_Y(\D,\bet^{\re}+2\bet^{\ima})$ points
in $F \setminus E$.
The set
$V^\R_Y(\D,\alp,\bet^{\re},\bet^{\ima},\bp^\flat,\bp^\sharp)$
is finite
and consists of immersed curves (see Lemma~\ref{l2}).
We
put
\begin{equation}W_{Y,E,\varphi}(\D,\alp,\bet^{\re},\bet^{\ima},\bp^\flat,\bp^\sharp)
=\sum_{C \in V^\R_Y(\D,\alp,\bet^{\re},\bet^{\ima},\bp^\flat,\bp^\sharp)}\mu_\varphi(C)
\ ,\label{e2074}\end{equation} where
\begin{equation}
\mu_\varphi(C)=
(-1)^{s(C)+C_{1/2}\circ\;\varphi}\ . \label{e2203}
\end{equation}
and $C_{1/2}$ is
the
image of one of the halves of
$\PP^1 \setminus \R P^1$ by the normalization map
$\PP^1\to C$
if $C$ is an irreducible real curve,
and
one of the
irreducible components of $C$ if $C$
is a pair of complex conjugate irreducible curves.
The number $\mu_\varphi(C)$ is called
{\it {\rm (}modified{\rm )} Welschinger sign}.

The proof of the following proposition
literally
coincides with the proof of \cite[Proposition 11]{IKS7}.

\begin{proposition}\label{l12} Let $(Y,E,F,\varphi)$
be a basic quadruple.
Fix a tuple $(D, \alp, \bet^\re, \bet^\ima)$,
where~$D \in \PicPPR(Y , E)$ is a divisor class,
$\alp,\bet^\re\in\Z_+^{\infty,\;\odd}$, and $\bet^\ima \in
\Z_+^{\infty}$ are such that $R_Y(D, \bet^\re + 2\bet^\ima) > 0$.
Choose two point sequences $\bp^\flat$ and $\bp^\sharp$ satisfying
the following restrictions:
\begin{enumerate}
\item[{\rm (r1)}] the tuple $(D, \alp, \bet^\re, \bet^\ima, \bp^\flat)$
is admissible,
\item[{\rm (r2)}] the number of points in $\bp^\sharp$
is equal to $R_Y(D, \bet^\re + 2\bet^\ima)$,
\item[\rm {(r3)}] the pair $(\bp^\flat, \bp^\sharp)$ is in $D_0$-CH position
for some divisor class $D_0 \in \PicPPR(Y, E)$, $D_0 \geq D$.
\end{enumerate}
Then, the number $W_{Y,E,\varphi}(D, \alp, \bet^\re, \bet^\ima,
\bp^\flat, \bp^\sharp)$ does not depend on the choice of sequences
$\bp^\flat$ and $\bp^\sharp$ subject to
{\rm (r1)-(r3)}. \proofend
\end{proposition}

\begin{proposition}\label{ini1} The only non-zero
numbers
$W_{Y,E,\varphi}(\D,\alp,\bet^{\re},\bet^{\ima},\bp^\flat,\emptyset)$
for admissible tuples $(\D,\alp,\bet^\re,\bet^\ima,\bp^\flat)$ such
that $\D\in\PicPPR(Y,E)$, $\alp,\bet^{\re} \in \Z^{\infty,\; \odd}$,
$\bet^{\ima}\in\Z^\infty_+$, and
$$I(\alp+\bet^{\re}+2\bet^{\ima})=[\D]E>0,
\quad R_Y(\D,\bet^{re}+2\bet^{\ima})=0\ ,$$ are the following ones:
\begin{enumerate}
\item[{\rm (1)}] if $\D=D$ is a divisor class,
\begin{enumerate}
\item[{\rm (1i)}]
$W_{Y,E,\varphi}(E',0,e_1,0,\emptyset,\emptyset) =
(-1)^{E'_{1/2}\circ\varphi}$, where $E'\in{\mathcal E}(E)$ is real;
\item[{\rm (1ii)}] $W_{Y,E,\varphi}(-(K_Y
+E),e_1,e_1,0,\bp^\flat,\emptyset) = (-1)^{L_{1/2}\circ\varphi}$, where
$L\in|-(K_Y+E)|$ is real, $\R L\subset F$;
\item[{\rm (1iii)}]
$W_{Y,E,\varphi}(D,\alp,0,0,\bp^\flat,\emptyset) =
(-1)^{C_{1/2}\circ\varphi}$, where $-(K_Y+E)D=1$, $I\alp=DE$, $C\in
V_Y^\R(D,\alp,0,0,\bp^\flat,\emptyset)$;
\end{enumerate}
\item[{\rm (2)}] if $\D$ is a pair of divisor classes,
\begin{enumerate}
\item[{\rm (2i)}] if $E'_1,E'_2\in{\mathcal E}(E)$ are complex conjugate,
$E'_1E'_2=1$, then
$$W_{Y,E,\varphi}(\{E'_1,E'_2\},0,0,e_1,\emptyset,\emptyset) =-(-1)^{E'_1\circ\varphi}\ ,$$
\item[{\rm (2ii)}] if $E'_1,E'_2\in{\mathcal E}(E)$ are disjoint complex conjugate, then
$$W_{Y,E,\varphi}(\{E'_1,E'_2\},0,0,e_1,\emptyset,\emptyset) =(-1)^{E'_1\circ\varphi}\ ,$$
\item[{\rm (2iii)}]
$W_{Y,E,\varphi}(\{-(K_Y+E),-(K_Y+E)\},0,0,e_2,\emptyset,\emptyset)=1$,
if $L',L''$ are complex conjugate.
\end{enumerate}
\end{enumerate}
\end{proposition}

{\bf Proof}. Proposition \ref{ini1} can easily be derived from Lemma
\ref{l3}. Notice only that, in case (2iii),
$L'\circ\varphi
=
0\mod2$ since the linear system
$|-(K_Y+E)|$, which contains $L'$, contains also a real rational
curve
whose complex locus is divided into two halves by its real locus located in $F$.
\proofend

The numbers
$W_{Y,E,\varphi}(\D,\alp,\bet^{\re},\bet^{\ima},\bp^\flat,\emptyset)$
in Proposition \ref{ini1} do not depend on
the choice of $\bp^\flat$.

We skip $\bp^\flat$ and $\bp^\sharp$ in the notation of the
numbers appearing in Propositions~\ref{l12} and~\ref{ini1}, and
write $W_{Y,E,\varphi}(\D,\alp,\bet^{\re},\bet^{\ima})$ for these
numbers calling them
{\it ordinary $w$-numbers}.

\subsection{Formula for ordinary $w$-numbers}\label{secRF}

\begin{theorem}\label{t1} Let $(Y,E,F,\varphi)$
be a basic quadruple.

{\rm (1)} For any divisor class $D\in\PicPPR(Y,E)$ and vectors
$\alp,\bet^{\re}\in\Z_+^{\infty,\;\odd}$,
$\bet^{\ima}\in\Z_+^\infty$ such that
$I(\alp+\bet^{\re}+2\bet^{\ima})=DE$, $R_Y
(D,\bet^{\re}+2\bet^{\ima})\ge0$, and $\bet^{\ima}\ne 0$, one has
\begin{equation}W_{Y,E,\varphi}(D,\alp,\bet^{\re},\bet^{\ima})=0\
.\label{e43}\end{equation}

{\rm (2)} For any divisor class $D\in\PicPPR(Y,E)$ and vectors
$\alp,\bet\in\Z_+^{\infty,\;\odd}$ such that
\mbox{$I(\alp+\bet)=DE$} and $R_Y(D,\bet)>0$, one has
$$W_{Y,E,\varphi}(D,\alp,\bet,0)=\sum_{j\ge 1,\
\bet_j>0}W_{Y,E,\varphi}(D,\alp+e_j,\bet-e_j,0)$$
$$+(-1)^{E_{1/2}\circ\varphi}\sum(-1)^{(I\bet^{(0)}+I\alp^{(0)})(L_{1/2}\circ\varphi)}\cdot
\frac{2^{\|\bet^{(0)}\|}}{\bet^{(0)}!}(l+1)\left(\begin{matrix}\alp\\
\alp^{(0)}\alp^{(1)}...\alp^{(m)}\end{matrix}\right)
\frac{(n-1)!}{n_1!...n_m!}$$
\begin{equation}\times\prod_{i=1}^m\left(\left(
\begin{matrix}(\bet^{\re})^{(i)}\\
\gam^{(i)}\end{matrix}\right)W_{Y,E,\varphi}(\D^{(i)},\alp^{(i)}
,(\bet^{\re})^{(i)},(\bet^{\ima})^{(i)})\right)
,\label{e44}\end{equation} where $L$ is any real curve in
$|-(K_Y+E)|$ with $\R L\subset F$, $$n = R_Y(D, \bet),\quad n_i=R_Y
(\D^{(i)}, (\bet^{\re})^{(i)}+2(\bet^{\ima})^{(i)}),\ i=1,...,m\ ,
$$
and the second sum in~(\ref{e44}) is taken \begin{itemize}\item over
all integers $l\ge0$ and vectors $\alp^{(0)}\le\alp$,
$\bet^{(0)}\le\bet^{\re}$; \item over all sequences
\begin{equation}(\D^{(i)},\alp^{(i)},(\bet^{\re})^{(i)},(\bet^{\ima})^{(i)}),\
1\le i\le m \ ,\label{e45}
\end{equation} such that, for all $i=1,....,m$,
\begin{enumerate}
\item[{\rm (1a)}] $\D^{(i)}\in\PicPPR(Y,E)$, and $\D^{(i)}$ is neither the
divisor class $-(K_Y+E)$, nor the pair $\{-(K_Y+E),-(K_Y+E)\}$,
\item[{\rm (1b)}]
$I(\alp^{(i)}+(\bet^{\re})^{(i)}+2(\bet^{\ima})^{(i)})=[\D^{(i)}]E$,
and $R_Y(\D^{(i)}, (\bet^{\re})^{(i)}+2(\bet^{\ima})^{(i)})\ge0$,
\item[{\rm (1c)}] $\D^{(i)}$ is a pair of divisor classes
if and only if $(\bet^{\ima})^{(i)}\ne0$,
\item[{\rm (1d)}] if $\D^{(i)}$ is a pair of divisor classes, then
$n_i=0$ and $\alp^{(i)}=(\bet^{\re})^{(i)}=0$,
\end{enumerate}
and
\begin{enumerate}
\item[{\rm (1e)}]
$D-E=\sum_{i=1}^m[\D^{(i)}]- (2l+I\alp^{(0)}+I\bet^{(0)})(K_Y+E)$,
\item[{\rm (1f)}] $\sum_{i=0}^m\alp^{(i)}\le\alp$, $\sum_{i=0}^m(\bet^{\re})^{(i)}\ge\bet$,
\item[{\rm (1g)}] each
tuple $(\D^{(i)},0,(\bet^{\re})^{(i)},(\bet^{\ima})^{(i)})$ with
$n_i=0$ appears in (\ref{e45}) at most once,
\end{enumerate}
\item over all sequences
\begin{equation}\gam^{(i)}\in\Z_+^{\infty,\;
\odd},\quad
\|\gam^{(i)}\|=\begin{cases}1,\ &\D^{(i)}\ \text{is a divisor class}, \\
0,\ &\D^{(i)}\ \text{is a pair of divisor classes},\end{cases}\quad
i=1,...,m\ ,\label{e41n}
\end{equation} satisfying
\begin{enumerate}
\item[{\rm (2a)}] $(\bet^{\re})^{(i)}\ge\gam^{(i)}$, $i=1,...,m$, and $\sum_{i=1}^m\left((\bet^{\re})^{(i)}-
\gam^{(i)}\right)=\bet^{\re}-\bet^{(0)}$,
\end{enumerate}
\end{itemize}
and the second sum in (\ref{e44}) is factorized by simultaneous
permutations in the sequences (\ref{e45}) and (\ref{e41n}).

{\rm (3)} All ordinary $w$-numbers $W_{Y,E,\varphi}(D,\alp,\bet,0)$, where
$D\in\PicPPR(Y,E)$ is a divisor class and
$R_Y(D, \bet) > 0$, are recursively determined by the formula (\ref{e44}) and
the initial conditions given by Proposition \ref{ini1}.
\end{theorem}

\begin{remark}
\label{r2} It is easy to verify that $n-1=\sum_in_i+\|\bet^{(0)}\|$
{\rm (}in the notation of Theorem \ref{t1}{\rm )}.
\end{remark}

The proof of Theorem \ref{t1}
literally coincides with the proof of \cite[Theorem 1 and Corollary
14]{IKS7} (notice that in the journal version of \cite[Section 3.3]{IKS7}
one deformation label of type (DL1) is missing;
the complete list of the deformation labels can be found in the ArXiv version of \cite{IKS7}).

We present here an immediate consequence that will be used below.

\begin{corollary}\label{c1}
Under the hypotheses of Theorem \ref{t1}(2), assume in addition that
$F\setminus\R E$ is disconnected, $DE=0$, and $R_Y(D,0)\ge2$. Then
$$W_{Y,E,\varphi}(D,0,0,0)=0\ .$$
\end{corollary}

{\bf Proof}. This follows from Proposition \ref{l12}: indeed, we may
choose two of the points of $\bp^\sharp$ in different components of
$F\setminus\R E$ making the set $V^\R_Y
(\D,0,0,0,\emptyset,\bp^\sharp)$ empty, since a real rational curve
cannot have two one-dimensional real components. \proofend

\subsection{Sided $w$-numbers}\label{sec6}
Let $(Y,E,F,\varphi)$ be a basic quadruple. Suppose in addition that
$F\setminus\R E$ splits into two connected components $F_+$ and
$F_-$. In this case, $(Y,E,F,\varphi)$ is called {\it dividing basic
quadruple}.

Let $(\D,\alpha,\beta^{\re},\beta^{\ima},\bp^\flat)$
be an admissible tuple.
Choose a sequence $\bp^\sharp$
of $R_Y(\D,\bet^{\re}+2\bet^{\ima})$ points in~$F_+$. Suppose that
the pair of point sequences $\bp^\flat,\bp^\sharp$ is in a $D_0$-CH
position with respect to some divisor class $D_0\in\PicPP(Y,E)$,
$D_0\ge[\D]$.
Put
$$V^\R_{Y,F_+}(\D,\alpha,\beta^{\re},\beta^{\ima},\bp^\flat,\bp^\sharp)
=\{C\in
V^\R_Y(\D,\alpha,\beta^{\re},\beta^{\ima},\bp^\flat,\bp^\sharp)\ :\
\card(C \cap F_-) < \infty\}\ .$$
Clearly,
if $\D$ is a pair of divisor classes, then
$$V^\R_{Y,F_+}(\D,\alpha,\beta^{\re},\beta^{\ima},\bp^\flat,\bp^\sharp)
=V^\R_Y(\D,\alpha,\beta^{\re},\beta^{\ima},\bp^\flat,\bp^\sharp)\
.$$

Set
\begin{equation}W_{Y,F_+,\varphi}^\eps(\D,\alpha,\beta^{\re},\beta^{\ima},
\bp^\flat,\bp^\sharp)=\sum_{C\in
V^\R_{Y,F_+}(\D,\alpha,\beta^{\re},\beta^{\ima},\bp^\flat,\bp^\sharp)}\mu^\eps_\varphi(C),\quad\eps=\pm
\ ,\label{e2200}\end{equation}
where $\mu^+_\varphi(C)=\mu_\varphi(C)$ is defined by (\ref{e2203})
and
\begin{equation}\mu^-_\varphi(C)=
(-1)^{s(C) +
C_{1/2}
\circ\varphi+
\card(C_{1/2} \cap F_-)}. \label{e2202}\end{equation}

\begin{remark}\label{r1}
By Lemma \ref{l2}(2),
if $\D\in\PicPPR(Y,E)$ is a pair of divisor classes, then
$$W_{Y,F_+,\varphi}^+(\D,\alpha,\beta^{\re},\beta^{\ima},\bp^\flat,
\bp^\sharp)=0$$ as long as $\alpha+\beta^{\re}>0$ or
$R_Y(\D,\beta^{\re}+2\beta^{\ima})>0$.
\end{remark}

\begin{proposition}\label{ll12} Let $(Y,E,F,\varphi)$
be a dividing basic quadruple.
Fix an admissible tuple $(D, \alp, \bet^\re, \bet^\ima,\bp^\flat)$,
where~$D \in \PicPPR(Y , E)$ is a divisor class,
$\alp,\bet^\re\in\Z_+^{\infty,\;\even}$,  $\bet^\ima \in
\Z_+^{\infty}$, and
$R_Y(D, \bet^\re + 2\bet^\ima) > 0$.
Choose two point sequences $\bp^\flat\subset\R E$ and
$\bp^\sharp\subset F_+$ satisfying the restrictions {\rm (r1)-(r3)} of
Proposition \ref{l12}. Then, the numbers $W_{Y,F_+,\varphi}^\pm(D,
\alp, \bet^\re, \bet^\ima, \bp^\flat, \bp^\sharp)$ do not depend on
the choice of sequences $\bp^\flat$ and $\bp^\sharp$ subject to
{\rm (r1)-(r3).}
\end{proposition}

The
proof of Proposition~\ref{ll12}
is given in Section~\ref{sec7}.

\begin{proposition}\label{ini2} Let $(Y,E,F,\varphi)$ be
a dividing basic quadruple,
and let $(\D,\alp,\bet^\re,\bet^\ima,\bp^\flat)$
be an admissible tuple
such
that
$\alp,\bet^{\re} \in \Z^{\infty,\; \even}_+$
and
$$I(\alp+\bet^{\re}+2\bet^{\ima})=[\D]E>0,
\quad R_Y(\D,\bet^{re}+2\bet^{\ima})=0\ .$$

\begin{enumerate}
\item[{\rm (1)}]
Assume that
$W_{Y,F_+,\varphi}^\pm(\D,\alp,\bet^{\re},\bet^{\ima},\bp^\flat,\emptyset)
\ne 0$
and $\D=D$
is a divisor class.
Then,
\begin{enumerate}
\item[{\rm (1i)}]
either $D = -(K_Y + E)$, $\alp = 0$, $\bet^{\re} = e_2$,
$\bet^{\ima} = 0$, $\bp^\flat = \emptyset$,
and the supporting curves $L'$, $L''$ are both real,
\item[{\rm (1ii)}]
or
$-(K_Y+E)D=1$, $I\alp=DE$, and $D$ is represented by
a curve $C \in V^\R_{Y,F_+}(\D,\alpha,0,0,\bp^\flat,\emptyset)$
with $\R C\subset \overline F_+$.
\end{enumerate}
In the first case,
$W_{Y,F_+,\varphi}^\pm(-(K_Y +E),0,e_2,0,\emptyset,\emptyset) =
\lambda(-1)^{L'_{1/2}\circ\varphi}$,
where $\lambda$ is the number of supporting curves $L'$, $L''$
whose real part is contained in $\overline F_+$.
In the second case,
$W_{Y,F_+,\varphi}^\pm(D,\alp,0,0,\bp^\flat,\emptyset) =
(-1)^{C_{1/2}\circ\varphi}$.
\item[{\rm (2)}]
Assume that
$W_{Y,F_+,\varphi}^\pm(\D,\alp,\bet^{\re},\bet^{\ima},\bp^\flat,\emptyset)
\ne 0$
and
$\D$
is a pair of divisor classes.
Then,
\begin{enumerate}
\item[{\rm (2i)}]
either
$\D = \{E'_1, E'_2\}$, where
$E'_1,E'_2\in{\mathcal E}(E)$
are complex conjugate,
\item[{\rm (2ii)}]
or
$\D = \{L', L''\}$ and
the supporting curves
$L',L''$ are complex conjugate.
\end{enumerate}
In the first case,
\begin{eqnarray}
W_{Y,F_+,\varphi}^+(\{E'_1,E'_2\},0,0,e_1,\emptyset,\emptyset)
&=&(-1)^{E'_1\circ\varphi + E'_1 \circ E'_2}
,\nonumber
\\
W_{Y,F_+,\varphi}^-(\{E'_1,E'_2\},0,0,e_1,\emptyset,\emptyset)
&=&(-1)^{E'_1\circ\varphi + \card(E'_1\cap F_-) + E'_1 \circ E'_2}\ .
\nonumber
\end{eqnarray}
In the second case,
$$W_{Y,F_+,\varphi}^\pm(\{-(K_Y+E),-(K_Y+E)\},0,0,e_2,\emptyset,\emptyset)
=1\ .$$
\end{enumerate}
\end{proposition}

{\bf Proof}.
The statement
can be easily derived from Lemma~\ref{l3},
taking into account that $L'\circ\varphi
=
0\mod2$ in
(2ii) ({\it cf.}, the proof of Proposition \ref{ini1}).
\proofend

We skip $\bp^\flat$ and $\bp^\sharp$ in the notation of the
numbers
appearing in Propositions~\ref{ll12} and~\ref{ini2}, and write
$W_{Y,F_+,\varphi}^\pm(\D,\alp,\bet^{\re},\bet^{\ima})$ for these
numbers calling them
{\it sided $w$-numbers}.

\subsection{Sided $w$-numbers
in deformation diagrams}\label{new-sided}

Let $(Y,E,F,\varphi)$ be a dividing basic
quadruple,
and let $(\D,\alpha,\beta,
0, \bp^\flat)$
be an admissible tuple
such that $\D = D$
is a divisor class and $R_Y(\D,\bet)
> 0$.
Choose a sequence $\bp^\sharp$
of $R_Y(\D,\bet)$
points in~$F_+$,
and assume that
the pair of point sequences $\bp^\flat,\bp^\sharp$ is in a $D_0$-CH
position with respect to some divisor class $D_0\in\PicPPR(Y,E)$.

Put
$k =
\min J^\sharp$, where $\bp^\sharp=\{p_j(\sigma_j)\ :\ j\in J^\sharp\}$

(see Lemma~\ref{D0}),
and denote by $\overline{V}$ a deformation diagram
of $(D, \alpha, \beta,
0, \bp^\flat, \bp^\sharp
\setminus \{p_k\}, p_k)$.

\begin{lemma}\label{l4}
Let $\bet_j>0$,
and $C\in V_{Y,F_+}^\R(D,\alp+e_j,\bet-e_j,0,\bp^\flat\cup\{p_k(0)\},
\bp^\sharp\setminus\{p_k\})$
intersects~$E$ at $p_k(0)$ with multiplicity $j$. Then,
the
real leaves of $\overline{V}$
that are generated by the root~$C$ consist of two curves
$C_1,C_2 \in V^\R_{Y, F_+}(D, \alp, \bet, 0, \bp^\flat, \bp^\sharp)$,
and $\mu^\pm_\varphi(C_i)=\mu^\pm_\varphi(C)$, $i = 1, 2$.
\end{lemma}

{\bf Proof.} Choose local coordinates $x,y$ in a neighborhood of
$p_k(0)$ so that $E=\{y=0\}$, $F_+=\{y>0\}$, $\Lam_k=\{(0,t)\ :\
t\ge0\}$, and
$$C=\left\{ay+bx^j+\sum_{m+jn > j}c_{mn}
x^my^n=0\right\},\quad a,b\in\R^*\ .$$
Since $
\card(\R C \cap F_-)
<\infty$,
the multiplicity
$j$ is even,
and $ab<0$.
Hence, the root~$C$
has two real branches
given in a neighborhood of $p_k(0)$
by ({\it cf.} \cite[Formulas (22)
and (23)]{MS2})
$$
C_i(t)=\left\{ay+b(x+\tau)^j+\sum_{m+jn > j}O(\tau)\cdot
(x+\tau)^my^n=0\right\},
$$
where
$\tau=\left(-\frac{a}{b}\right)^{1/j}t^{1/j}+
o(t^{1/j})$ for $i = 1$,
and
$\tau=-\left(-\frac{a}{b}\right)^{1/j}t^{1/j}+
o(t^{1/j})$
for $i = 2$.
For each curve $C_i(t)$, $i=1,2$, one has $\mu^\pm_\varphi(C_i)=\mu^\pm_\varphi(C)$.
Indeed, the above local formula insures
that the topology of the curves is preserved in a neighborhood of $p_k(0)$; outside of a
neighborhood of $p_k(0)$, the equality required follows from Lemma~\ref{ordinary_invariance}.
\proofend

\begin{lemma}\label{l5}
Let $C = E \cup \check{C}$ be a root of $\overline{V}$
such that $C$ generates at least one leaf belonging to
$V_{Y,F_+}^\R(D,\alp,\bet,0,\bp^\flat,\bp^\sharp)$.
Then,
$\check{C}$ splits in {\rm primary components}
from the following list:

\begin{enumerate}
\item[{\rm (i)}] pairs of reduced complex conjugate components as described in Lemma \ref{l3}(2);
\item[{\rm (ii)}] real reduced components, whose
all intersection points with $E$ are real and have even
multiplicity; each of these components belongs to $V^\R_{Y,F_+}(D', \alp', \bet', 0, (\bp^\flat)', (\bp^\sharp)')$
for a certain $D_0$-proper tuple $(D', \alp', \bet', 0, (\bp^\flat)',(\bp^\sharp)')$;
\item[{\rm (iii)}] non-reduced components $s'L'$, $s''L''$,
where $L',L''\in|-(K_Y+E)|$ are
the supporting curves,
and, in
addition,
$s'=s''$ if $L',L''$ are complex
conjugate, and $s'$ (respectively, $s''$) is even if $L'$ (respectively, $L''$) is real with
$\R L'\subset\overline F_-$ (respectively, $\R L''\subset\overline F_-$);
\item[{\rm (iv)}] non-reduced components $sL(z)$, where $s$ is even, $z\in\bp^\sharp\setminus\{p_k\}$,
and $L$ is the curve belonging to $|-(K_Y+E)|$
and passing through $z$;
\item[{\rm (v)}] non-reduced components $iL(p_{ij})$, where $i$ is even, $p_{ij}\in\bp^\flat$, and
$L(p_{ij})$ is the curve belonging to $|-(K_Y+E)|$
and passing through $p_{ij}$.\end{enumerate}
Moreover,\begin{enumerate}
\item[{\rm (1)}] pairwise intersections of distinct
primary components are either empty, or transversal contained
in $Y\setminus E$,
and all reduced primary components are immersed and are nonsingular
along $E$,
\item[{\rm (2)}]
the intersection multiplicities of $\check C$ and $E$ at the points $\check C\cap\bp^\flat$ are encoded by a vector $\check\alp\le\alp$,
\item[{\rm (3)}] there exists a conjugation invariant set $\bz\subset\check C\cap E\setminus\bp^\flat$ containing
exactly one point of each irreducible component of $\check C$ and such that the intersection multiplicities
of $\check C$ and $E$ at the points of $\check C\cap E\setminus(\bp^\flat\cup\bz)$ are encoded by the vector $\beta$.
\end{enumerate}
\end{lemma}

{\bf Proof.} The statement
follows from
\cite[Proposition 2.6]{MS2}.
The fact that each real reduced primary component of $\check{C}$
intersects $E$ only at real points is guaranteed by Lemmas~\ref{l3} and~\ref{imaginary};
these intersection points have even
multiplicity due to the assumption that the $C$ generates a leaf belonging to
$V_{Y,F_+}^\R(D,\alp,\bet,0,\bp^\flat,\bp^\sharp)$.
\proofend

To describe the leaves of deformation diagrams with reducible roots,
we
use
certain
deformation labels
from the list (DL1)-(DL9)
introduced in \cite[Section 3.3]{IKS7}.
Each label
is seen as a curve on the toric surface
that is determined by the Newton polygon of the polynomial defining the label.
We recall these defining polynomials (some of them are slightly modified by a conjugation-invariant coordinate change),
where $\cheb_k(t)=\cos\arccos kt$ is the $k$-th Chebyshev polynomial and $y_k$ is the only simple positive root of $\cheb_k(t)-1$:
\begin{enumerate}\item[(DL2)$_j$] for an even $j>0$, two deformation labels defined by the equations
$$\psi_1(x,y)=y^2+1-y\cdot\cheb_j(x),\quad\psi_2(x,y)=\psi_1(x\sqrt{-1},y)\ ,$$
\item[(DL3)$_i$] for an even $i>0$, a deformation label $$(x-1)(y^i-x)=0\ ,$$ \item[(DL5)$_s$] for an even $s>0$, two deformation labels
$$(x-1)(x((y\pm1)^s+1)-1)=0\ ,$$
\item[(DL6)$_{s,\eps_0}$] for an integer $s>0$ and $\eps_0=\pm1$,
a deformation label
$$\frac{y+x^2}{2y}\left(\cheb_{s+1}\left(y_{s+1}-\frac{y\eps_0}{2^{(s-1)/(s+1)}}\right)-1\right)+1=0\ ,$$
\item[(DL7)$_s$] for $s>0$, $s+1$ pairs of complex conjugate deformation labels $$\frac{y+\sqrt{-1}x^2}{2y}\left(\cheb_{s+1}\left(\frac{y\eps}{2^{(s-1)/(s+1)}}
+y_{s+1}\right)-1\right)+1=0\ ,$$
$$\frac{y-\sqrt{-1}x^2}{2y}\left(\cheb_{s+1}\left(\frac{y\overline\eps}{2^{(s-1)/(s+1)}}
+y_{s+1}\right)-1\right)+1=0\ ,$$ where $\eps^{s+1}=1$.
\end{enumerate}

\begin{lemma}\label{l6}
Let a curve $C = E \cup \check{C}\in|D|$ be such that the primary components of $\check C$ belong to the list
(i)-(v) and satisfy conditions (1)-(3) in the statement of Lemma \ref{l5}.
Then,
the curve $C$ is a root of $\overline{V}$. It generates at least one leaf that belongs to
$V_{Y,F_+}^\R(D,\alp,\bet,0,\bp^\flat,\bp^\sharp)$ if and only if
the following condition holds:
$\check C$
has no any real primary component
$s'L'$ or $s''L''$
whose multiplicity ($s'$ or $s''$) is odd
and the real point set is contained in $\overline F_-$.
If the set
of leaves generated by~$C$ and belonging to $V_{Y,F_+}^\R(D,\alp,\bet,0,\bp^\flat,\bp^\sharp)$ is non-empty, then it is in one-to-one correspondence with the following data: \begin{itemize}\item a set $\bz\subset\check C\cap E\setminus\bp^\flat$
satisfying the condition (3) of Lemma \ref{l5}, \item a collection ${\mathcal{DL}}$ of deformation labels chosen as follows:

- a deformation label of type (DL2)$_j$ for each point $q\in\bz'$ satisfying $(\check C\cdot E)_q=j$, where $\bz'\subset\bz$
consists of the points which do not lie on the primary components $s'L'$, $s''L''$,
$iL(p_{ij})$, $sL(z)$ of $\check C$,

- a deformation label of type (DL3)$_i$ for each primary component $iL(p_{ij})$ of $\check C$,

- a deformation label of type (DL5)$_s$ for each primary component $sL(z)$ of $\check C$,

- a deformation label of type (DL6)$_{s',
+1}$ (respectively, (DL6)$_{s'',
+1}$) for a real primary component $s'L'$ (respectively, $s''L''$)
of $\check C$ with $\R L'\subset\overline
F_+$
(respectively, $\R L''\subset\overline
F_+$), if $s'$ (respectively, $s''$) is even,

- a deformation label of type (DL6)$_{s', -1}$ (respectively, (DL6)$_{s'',-1}$)
for a real primary component $s'L'$ (respectively, $s''L''$)
of $\check C$ with $\R L'\subset\overline F_-$
(respectively, $\R L''\subset\overline F_-$), if $s'$ (respectively, $s''$) is even,

- a deformation label of type (DL6)$_{s',
+1}$ or (DL6)$_{s',
-1}$ (respectively, (DL6)$_{s'',
+1}$
or (DL6)$_{s'',
-1}$) for a real primary component $s'L'$ (respectively, $s''L''$)
of $\check C$ with $\R L'\subset\overline F_+$
(respectively, $\R L''\subset\overline F_+$), if $s'$ (respectively, $s''$) is odd,

- a pair of complex conjugate deformation labels of type (DL7)$_s$
for a pair of complex conjugate primary components $sL'$, $sL''$ of $\check C$.
\end{itemize}
\end{lemma}

{\bf Proof.} Since the primary components of $\check C$ belong to the list
(i)-(v) and satisfy conditions (1)-(3) in the statement of Lemma \ref{l5},
we obtain
(see \cite[Lemma 2.20]{MS2})
that $C$ is a root of the deformation diagram $\overline V$, and
the leaves of $\overline V$ are in one-to-one correspondence with the set of pairs $(\bz,{\mathcal{DP}})$,
where $\bz\subset\check C\cap E\setminus\bp^\flat$
is a set of points
satisfying the condition (3) of Lemma \ref{l5}, and ${\mathcal{DP}}$ is a conjugation-invariant collection  of deformation
patterns
$$\{\Psi_q\ :\ q\in\bz'\},\quad\{\Psi_{iL(p_{ij})}\ :\ iL(p_{ij})\subset\check C\},\quad\{\Psi_{sL(z)}\ :\ sL(z)\subset\check C\},
$$ $$\{\Psi_{s'L'}\ :\ s'L'\subset\check C\},\quad\{\Psi_{s''L''}\ :\ s''L''\subset\check C\},$$ where $iL(p_{ij})$, $sL(z)$, $s'L'$,
or $s''L''$ run over the corresponding primary components of $\check C$, and $\Psi_*$
denote specific curves on toric surfaces
introduced in \cite[Section 2.5.3]{MS2}.
We consider these pairs $(\bz,{\mathcal{DP}})$ in detail
and prove that they induce leaves
belonging to $V_{Y,F_+}^\R(D,\alp,\bet,0,\bp^\flat,\bp^\sharp)$ if and only if $\check C$ does not have
a primary component
$s'L'$ or $s''L''$
such that this component is real,
its real point set is contained in $\overline F_-$,
and its multiplicity ($s'$ or $s''$) is odd.

A deformation of $C$ into any leaf-curve $\widetilde C\in
V_{Y,F_+}(D,\alp,\bet,0,\bp^\flat,\bp^\sharp(t))$, where
$\bp^\sharp(t)=(\bp^\sharp\setminus\{p_k\})\cup\{p_k(t)\}$,
$t>0$, can be described by a family of sections of
$H^0(Y,{\mathcal O}_Y(D))$: \begin{equation}S\check T_\tau+\tau^\kappa
T_\tau,\ \tau\in(\C,0),\quad \tau^\kappa=t\ ,\label{e2206}\end{equation}
where $S,\check T_0,T_0$ are real, $S^{-1}(0)=E$, $\check
T_0^{-1}(0)=\check C$, and $\kappa$ is the least common multiple of all the parameters
$j$, $i$, $s$, $s'+1$, $s''+1$ in the assertion of Lemma \ref{l6} ({\it cf.} \cite[Formula (57)]{MS2}).
In addition, $S$ and $\check T_0$ are negative
in $F_-$ (except possibly for a finite set
of points and lines),
and
\begin{equation}\begin{cases}&S(p_k(t))>0,\ \check T_0(p_k(t))<0,\
T_0(p_k(t))>0\quad\text{as}\ t>0,\\
&T_0(z)>0\quad\text{for all}\quad z\in\R E\setminus\bp^\sharp\ .\end{cases}\label{e2206b}\end{equation}
Indeed, formula (\ref{e2206}) follows from \cite[Lemma 2.10]{MS2}. Observe that $S$ does not divide $T_0$ in view of $t=\tau^k$ (see (\ref{e2206})). Hence
$T_0^{-1}(0)$ intersects $E$ at finitely many points and with even multiplicities,
since the same holds for each curve $C_\tau=\{S\check T_\tau+\tau^\kappa
T_\tau=0\}$, $\tau\ne0$. This claim combined with the facts that
$
\card(\R\check C\cap
F_-)
<\infty$,
$
\card(\R\widetilde C\cap F_-)
<\infty$,
$(S\check T_\tau+\tau^\kappa T_\tau)(p_k(t))=0$ for all $t>0$, and with the assumption $\alp+\bet\in\Z_+^{\infty,\;\even}$ yields
all sign conditions (\ref{e2206b}).

Given a point $q\in\bz'$, the intersection multiplicity $j=(\check C\cdot E)_q$ is even by
Lemma \ref{l5}. Choose local real coordinates $x,y$ in a neighborhood of $q$ in $Y$ so that
$E=\{y=0\}$, $q=(0,0)$, $F_-=\{y<0\}$. Then in formula (\ref{e2206}), we get
$\widetilde
T_0=y-2x^j+\text{h.o.t.}$, $T_0=a+\text{h.o.t.}$, where $a>0$ due to (\ref{e2206b}). Thus, by
\cite[Lemma 2.11]{MS2},
there are two real
deformation patterns $\Psi_q$ given by $y^2-2yP(x)+a=0$ with $P(x)=x^j+\text{l.o.t.}$, and they can be brought to the form (DL2)$_j$ by a transformation \begin{equation}\psi(x,y)\mapsto\lambda\psi(\mu x,\nu y),\quad \lambda,\mu,\nu>0\ .\label{e2206c}\end{equation}

Given a primary component $iL(p_{ij})$ of $\check C$ (with $i$ even by Lemma
\ref{l5}), it has a unique deformation pattern $\Psi_{iL(p_{ij})}$ (see
\cite[Lemma 2.15]{MS2}) which is real and can be brought to the form (DL3)$_i$ by transformation (\ref{e2206c}).

Given a primary component $sL(z)$ of $\check C$ (with $s$ even by
Lemma \ref{l5}) and an intersection point $q\in L(z)\cap E$ belonging to $\bz$, we are interested in deformation patterns that
describe a deformation of $C$ in a neighborhood of $L(z)$ such that the intersection point $q$ of $sL(z)$ and $E$
smoothes out, and the other intersection point of $sL(z)$ and $E$ turns into an intersection point of multiplicity $s$ of the deformed curve with $E$. Choose real coordinates $x,y$ in a neighborhood of $L(z)$ so that $L(z)=\{y=0\}$, $q=(0,0)$,
$z=(x_0,0)$,
$E=\{x^2-x+xy=0
\}$,
where $L(z)\cap F_+=\{(x,0)\ :\
0<x<1\}$. In particular, $0<x_0<1$, since $z\in\bp^\sharp\subset F_+$,
and in
(\ref{e2206}) we have $S(x,y)=x-x^2-xy$ and $\check T_0(x,y)=y^s(\varphi(x)-y\psi(x,y))$ with
$\varphi(0)=-1$, and hence $S(x,y)\check T_0(x,y)=y^s(x^2-x+xy)(y\psi(x,y)-\varphi(x))$.
By \cite[Section 2.5.3(8) and Lemma 2.16(1)]{MS2},
deformation of the curve $E\cup\check C$ in a neighborhood of $L(z)$
can be viewed as a patchworking of the curve $E\cup\check C$, given by $S(x,y)\check T_0(x,y)=0$,
and a curve, given by a polynomial $-(x-1)h(x,y)\varphi(x)$, where the factor
$\Psi_{sL(z)}=(x-1)h(x,y)$
(the deformation pattern for the pair $(sL(z),(0,0))$)
satisfies the relations $$h(x,y)=xf(y)+a,\quad f(y)+a=(y+\xi)^s,\ \xi\in\C, \quad
x_0f(0)+a=0\ .$$ Notice that the coefficients of the common monomials $x^iy^s$, $i\ge0$, for
$S(x,y)\check T_0(x,y)$ and $-(x-1)h(x,y)\varphi(x)$ respectively coincide. Furthermore, in this presentation,
$T_0=-(x-1)h(x,0)\varphi(x)$,
and, in particular,
$-a=T_0(q) > 0$
by (\ref{e2206b}).
From this we easily derive that $\xi^s=-a\frac{1-x_0}{x_0}>0$, obtaining two real deformation patterns that can be brought to the form (DL5)$_s$ via transformation (\ref{e2206c}).

If $\check C\supset s(L'\cup L'')$, where $L',L''$ are complex conjugate, then
there are $s$ pairs of complex conjugate deformation patterns for
these primary components (see \cite[Lemma 2.13]{MS2}), which can be brought to the form
(DL7)$_s$.

Let $L'$ be real, $\R L'\subset\overline F_+$, and $\check C\supset sL'$. By (\ref{e2206}) and (\ref{e2206b}), we can choose
real coordinates $x,y$ in a
neighborhood of $L'$ so that
\begin{equation}\begin{cases}&L'=\{y=0\},\quad S=y+x^2+xy,\quad q=(0,0)=E\cap L',\\
&F_+=\{S>0\},\quad \check
T_0=y^s((-1)^{s+1}+\text{h.o.t.}),\quad c=T_0(q)>0\ .\end{cases}\label{e2206d}\end{equation} Substituting this data to
the formulas of
\cite[Lemma 2.13]{MS2}, we obtain that the (complex)
deformation patterns of $sL'$ are given by the formula
\begin{equation}\Psi_{sL'}=\{(y+x^2)f(y)+c\;(-1)^{s+1}=0\}\ ,\label{e2209X}\end{equation} where
\begin{equation}yf(y)+c\;(-1)^{s+1}=\frac{c\;(-1)^{s+1}}{2}\left(\cheb_{s+1}\left(\frac{\xi y}{(2^{s-1}c\;(-1)^{s+1})^{1/(s+1)}}+y_{s+1}\right)+1\right)
\ ,\label{e2210X}\end{equation} $\xi^{s+1}=1$.
If $s$ is even, then there exists a unique real deformation pattern, and via transformation (\ref{e2206c}) preserving the terms $y+x^2$ in the above equation $S(x,y)=0$ for $E$ we can bring it to the form (DL6)$_{s,+1}$.
If $s$ is odd, then there exist two real deformation patterns.
Indeed, the equation for $\Psi_{sL'}$ can be rewritten as
\begin{equation}x^2=-y\cdot\frac{\cheb_{s+1}\left(\frac{\xi_0
y}{(2^{s-1}c)^{1/(s+1)}}+y_{s+1}\right)+1} {\cheb_{s+1}\left(\frac{\xi_0
y}{(2^{s-1}c)^{1/(s+1)}}+y_{s+1}\right)-1}\label{e2211X}\end{equation}
with $(2^{s-1}c)^{1/(s+1)}>0$, $\xi_0=\pm1$. It is easy to bring them to the form (DL6)$_{s,-\xi_0}$
via transformation of type (\ref{e2206c}).

Suppose that $L'$ is real, $\R L'\subset \overline F_-$, and $\check C$ contains a primary component $sL'$, $s>0$.
Then $s$ must be even, since otherwise the leaf-curves would contain a
one-dimensional branch in the domain $F_-$, which is not possible by the definition of
$V_{Y,F_+}^\R(D,\alp,\bet,0,\bp^\flat,\bp^\sharp)$. In view of sign conditions (\ref{e2206b}),
formulas (\ref{e2206d}) of the preceding paragraph turn into
$$\begin{cases}&L'=\{y=0\},\quad S=-(y+x^2+xy),\quad q=(0,0)=E\cap L',\\
&F_+=\{S>0\},\quad \check
T_0=y^s(-1+\text{h.o.t.}),\quad c=T_0(q)>0\ .\end{cases}$$
Substituting these data to
the formulas of
\cite[Lemma 2.13]{MS2}, we obtain the (complex)
deformation patterns of $sL'$ in the form
$$\Psi_{sL'}=\{(y+x^2)f(y)+c=0\}\ ,$$ where
$$yf(y)+c=\frac{c}{2}\left(\cheb_{s+1}\left(\frac{\xi y}{(2^{s-1}c)^{1/(s+1)}}+y_{s+1}\right)+1\right),
\quad\xi^{s+1}=1
\ ;$$ hence there exists a unique real deformation pattern, corresponding to
$\xi=1$, and it can easily be transformed to the type (DL6)$_{s,-1}$.
\proofend

Introduce the following numbers:
\begin{itemize}
\item for a deformation label $\Psi$ of type (DL2)$_j$, put
$\mu^+(\Psi)=(-1)^{s(\Psi)}$ and $\mu^-(\Psi)=(-1)^{s^-(\Psi)}$,
where
$s(\Psi)$ is the number of solitary nodes of $\Psi$, and $s^-(\Psi)$
is the number of solitary nodes lying in the domain $y>0$;
\item for a deformation label $\Psi$ of type (DL3)$_i$ or (DL5)$_s$, put
$\mu^+(\Psi)=\mu^-(\Psi)=1$;
\item for a deformation label $\Psi$ of type (DL6)$_{i, \eps_0}$ with even $i$, put
$\mu^+(\Psi)=(-1)^{s(\Psi)}$ and $\mu^-(\Psi)=(-1)^{s^-(\Psi)}$, where
$s(\Psi)$ is the number of solitary nodes of $\Psi$,
and $s^-(\Psi)$ is the number of solitary nodes lying in the domain $\eps_0(y+x^2)>0$;
\item for a deformation label $\Psi$ of type (DL6)$_{i, \eps_0}$ with odd $i$, put
$\mu^+(\Psi)=(-1)^{s(\Psi)}$ and $\mu^-(\Psi)=(-1)^{s^-(\Psi)}$, where
$s(\Psi)$ is the number of solitary nodes of $\Psi$,
and $s^-(\Psi)$ is the number of solitary nodes lying in the domain $y+x^2>0$.
\end{itemize}
Let
\begin{itemize}\item $\mu^\pm_{2,j}$, $\mu^\pm_{3,i}$, $\mu^\pm_{5,s}$ be the sums of the numbers $\mu^\pm(\Psi)$ over
all deformation labels of type (DL2)$_j$, (DL3)$_i$, (DL5)$_s$, respectively,
\item $\mu^\pm_{6,s,\eps_0}$, where $s$ is even, be the value of $\mu^\pm(\Psi)$ for the deformation label $\Psi$
of type (DL6)$_{s,\eps_0}$, $\eps_0=\pm1$,
\item
$\mu^\pm_{6,s}$, where $s$ is odd, be the sum of the numbers $\mu^\pm(\Psi)$ over the
two deformation labels of type (DL6)$_{s,+1}$ and (DL6)$_{s,-1}$.
\end{itemize}

\begin{lemma}\label{l7}
We have
 \begin{equation}\mu^+_{2,j}=0,\quad
\mu^-_{2,j}=\begin{cases}0,\quad j
=
0\mod4,\\ 2,\quad
j
=
2\mod4,\end{cases}\label{e2204}\end{equation}
\begin{equation}\mu^\pm_{3,i}=1,\quad\mu^\pm_{5,s}=2\ ,\label{e2207}\end{equation}
\begin{equation}\mu^\sig_{6,s,\eps_0}=\begin{cases}1,\quad & s
=
0\mod2,\ \eps_0=1,\ \sig=\pm\\
1,\quad & s
=
0\mod2,\ \eps_0=-1,\ \sig=+,\\ (-1)^{s/2},\quad & s
=
0\mod2,\ \eps_0=-1,\ \sig=-\end{cases}
\label{e2212}\end{equation}\begin{equation}
\mu^\sig_{6,s}=\begin{cases}0,\quad & s
=
1\mod2,\ \sig=+,\\
2,\quad & s
=
1\mod4,\ \sig=-,\\
0,\quad & s
=
3\mod4,\ \sig=-. \end{cases}\label{e2212a}\end{equation}
\end{lemma}

{\bf Proof.} All the relations follow from a direct computation.
\proofend

Let $C = E \cup \check{C}$ be a root of $\overline{V}$
such that $C$ generates a leaf $\widetilde C\in V_{Y,F_+}^\R(D,\alp,\bet,0,\bp^\flat,\bp^\sharp)$ corresponding to a pair $(\bz,{\mathcal{DL}})$
(see Lemma \ref{l6}). Introduce the following numbers: \begin{itemize}\item for each point $q\in\bz'$ (where $\bz'\subset\bz$
consists of the points which do not lie on the primary components $s'L'$, $s''L''$,
$iL(p_{ij})$, $sL(z)$ of $\check C$), put $\mu^\pm(C,q)=\mu^\pm_{2j}$, where $j=(\check C\cdot E)_q$;
\item for each primary component $sL(z)$ of $\check C$, put $\mu^\pm(C,sL(z))=\mu^\pm_{5,s}$;
\item for each real primary component $sL'$ (respectively, $sL''$) of $\check C$
such that $\R L' \subset \overline F_+$ (respectively, $\R L'' \subset \overline F_+$),
put $\mu^\pm(C,sL')$ (respectively, $\mu^\pm(C,sL'')$) equal to $\mu^\pm_{6,s,+1}$ if $s$ is even,
and
equal to $\mu^\pm_{6,s}$ if $s$ is odd,
\item for each real primary component $sL'$ (respectively, $sL''$) of $\check C$
such that $\R L' \subset \overline F_-$ (respectively, $\R L'' \subset \overline F_-$) and $s$ is even,
put $\mu^\pm(C,sL')$ (respectively, $\mu^\pm(C,sL'')$) equal to $\mu^\pm_{6, s, -1}$.
\end{itemize}

\begin{lemma}\label{sum-leaves}
Let $C = E \cup \check{C}$ be a root of $\overline{V}$
such that $C$ generates at least one leaf belonging to
$V_{Y,F_+}^\R(D,\alp,\bet,0,\bp^\flat,\bp^\sharp)$, and let $\Leaf(C)$ be the set of all such leaves. Then,
\begin{equation}\sum_{\widetilde C\in\Leaf(C)}\mu^\pm_\varphi(\widetilde
C)=(-1)^{E_{1/2}\circ\varphi}\cdot\mu^\pm_\varphi(\check C^{red})\cdot2^m\cdot M^\pm(C)\cdot\sum_{\bz}\prod_{q\in\bz'}\mu^\pm(C,q)
\ ,\label{e2213}\end{equation} where $\check C^{red}$ is the union of all reduced primary components
of $\check C$ different from $L',L''$, the exponent $m$ is the number of primary components $sL(z)$
of $\check C$, the factor $M^\pm(C)$ equals $(-1)^{(s'+s'')(L'_{1/2}\circ\varphi)}\mu^\pm(C,s'L')\mu^\pm(s''L'')$ if $s'L'$ and $s''L''$ are real primary components of
$\check C$, and equals $s+1$ if $\check C$ contains a pair of complex conjugate primary components $sL',sL''$,
and finally,
$\bz$ runs over all subsets
of $\check C\cap E\setminus\bp^\flat$ satisfying condition (3) of Lemma \ref{l5},
and $\bz'\subset\bz$
consists of the points which do not lie on the primary components $s'L'$, $s''L''$,
$iL(p_{ij})$, $sL(z)$ of $\check C$.
\end{lemma}

{\bf Proof.} Let $\widetilde C\in\Leaf(C)$. By \cite[Lemma
2.9]{MS2}, singular points of $\widetilde C$ (regarded as a small deformation of $C$) appear in a neighborhood
of $\Sing(C)$. Furthermore, the local branches do not glue up in local deformation of singular points in
$\Sing(\check C^{red})$, of intersection points $q\in\bp^\flat\cup(\check C\cap E)\setminus\bz$, and of the intersection points
of $\check C^{red}$ with the other primary components of $\check C$.
In particular, first, local deformations of the intersection points of $\check C^{red}$
with other primary components of $\check C$ and of the points $q\in\bp^\flat\cup(\check C\cap E)\setminus\bz$
do not bear solitary nodes,
and, second, due to Lemma \ref{ordinary_invariance}, the multiplicative contribution of $\Sing(\check C^{red})$ to
$\mu^\pm_\varphi(\widetilde C)$ is $\mu^\pm_\varphi(\check C^{red})$. Local deformations
of the primary components $iL(p_{ij})$, $sL(z)$, $s'L'$, $s''L''$ of $\check C$ and
of the points $q\in\bz'$ are determined by the corresponding deformation labels so that the solitary nodes of $\widetilde C$, which appear in these
deformations are in one-to-one correspondence with the solitary nodes of all deformation labels involved. Deformation labels of
type (DL3)$_i$, (DL5)$_s$, and (DL7)$_s$ do not have solitary nodes. The solitary nodes of the other deformation labels, which belong to the domains indicated in
the definition of numbers
$\mu^-(\Psi)$, correspond precisely to the nodes of $\widetilde C$ belonging to $F_+$. It
follows from the fact that, in the coordinates $x,y$ in the equations of a deformation label $\Psi$, this domain defines an intersection of
$F_+$ with a neighborhood of a point $q\in\bz'$ or
with a neighborhood of real primary components $s'L''$, $s''L''$ of $\check C$.

Then, formula (\ref{e2213}) immediately follows from
Lemmas \ref{l5} and \ref{l6}.
\proofend

\subsection{Formula for sided $w$-numbers}\label{sec7}

\begin{theorem}\label{t3} Let $(Y,E,F,\varphi)$ be a dividing basic
quadruple.

{\rm (1)} For a divisor class $D\in\PicPPR(Y,E)$ and vectors
$\alp,\bet^{\re},\bet^{\ima}\in\Z_+^\infty$ such that
$I(\alp+\bet^{\re}+2\bet^{\ima})=DE$
and $R_Y (D,\bet^{\re}+2\bet^{\ima})
> 0$,
one has
\begin{equation}W_{Y,F_+,\varphi}^\pm(D,\alpha,\beta^{\re},\beta^{\ima})=0\
,\label{e43o}\end{equation}
provided that either
$\alpha\not\in\Z_+^{\infty,\;\even}$, or
$\beta^{\re}\not\in\Z_+^{\infty,\;\even}$, or $\beta^{\ima}\ne0$.

{\rm (2)} For a divisor class $D\in\PicPPR(Y,E)$ and vectors
$\alp,\bet\in\Z_+^{\infty,\;\even}$ such that
\mbox{$I(\alp+\bet)=DE$} and $R_Y(D,\bet)>0$, one has:
\begin{enumerate}
\item[{\rm (2i)}] If $(K_Y+E)D=0$ or $(K_Y+E)D<-2$, then
\begin{equation}W_{Y,F_+,\varphi}^+(D,\alp,\bet,0)=0\
.\label{ee3}\end{equation}
\item[{\rm (2ii)}] If $(K_Y+E)D=-1$, then
\begin{equation}W_{Y,F_+,\varphi}^+(D,\alp,\bet,0)=2^{\|\bet\|}
W_{Y,F_+,\varphi}^+(D,\alp+\bet,0,0)\
.\label{ee4}\end{equation}
\item[{\rm (2iii)}] If $(K_Y+E)D=-2$, then
$$W_{Y,F_+,\varphi}^+(D,\alp,\bet,0)=2\sum_{j\ge 2,\
\bet_j>0}W_{Y,F_+,\varphi}^+(D,\alp+e_j,\bet-e_j,0)$$
\begin{equation}+4^{n-1}
(-1)^{E_{1/2}\circ\varphi}\cdot\sum(-1)^{I\alp^{(0)} \cdot (L_{1/2}\circ\varphi)} (l/2+1)
\left(\begin{matrix}\alp\\
\alp^{(0)}\end{matrix}\right)\prod_{i=1}^m
W_{Y,F_+,\varphi}^+(\D^{(i)},0 ,0,e_1)\ ,\label{e44o}\end{equation}
where $L\in|-(K_Y+E)|$ is real with $\R L\subset F$,
$n = R_Y(D, \bet)$,
and the second
sum in~(\ref{e44o}) is taken over all even integers $l \geq 0$, vectors
$\alp^{(0)}\le\alp$, and sequences of distinct tuples
\begin{equation}
(\D^{(i)},0,0,e_1),\
1\le i\le m \ ,\label{e45o}
\end{equation}
such that
\begin{itemize}
\item
each $\D^{(i)}\in\PicPPR(Y,E)$ is a pair of divisor classes that is different
from \mbox{$(-(K_Y+E),-(K_Y+E))$} and
satisfies $[\D^{(i)}]E=2$, $R_Y(\D^{(i)}, 2e_1)=0$,
\item $D-E=\sum_{i=1}^m[\D^{(i)}]-
(l+I\alp^{(0)}+I\bet)(K_Y+E)$;
\end{itemize}
the second sum in (\ref{e44o}) is factorized by permutations of
sequences (\ref{e45o}).
\end{enumerate}

{\rm (3)} For any divisor class $D\in\PicPPR(Y,E)$ and vectors
$\alp,\bet\in\Z_+^{\infty,\;\even}$ such that
\mbox{$I(\alp+\bet)=DE$} and $R_Y(D,\bet)>0$, one has
$$W_{Y,F_+,\varphi}^-(D,\alp,\bet,0)=2\sum_{j\ge 2,\
\bet_j>0}W_{Y,F_+,\varphi}^-(D,\alp+e_j,\bet-e_j,0)$$
$$+(-1)^{E_{1/2}\circ\varphi}\cdot\sum(-1)^{(I\alp^{(0)}+\bet^{(0)})(L_{1/2}\circ\varphi)}
\cdot\frac{4^{\|\bet^{(0)}\|}}{\bet^{(0)}!}\eta
(l)\cdot2^m\left(\begin{matrix}\alp\\
\alp^{(0)}\alp^{(1)}...\alp^{(m)}\end{matrix}\right)
\frac{(n-1)!}{n_1!...n_m!}$$
\begin{equation}\times\prod_{i=1}^m\left(\left(
\begin{matrix}
(\bet^{\re})^{(i)}\\
\gam^{(i)}
\end{matrix}
\right)W_{Y,F_+,\varphi}^-(\D^{(i)},\alp^{(i)}
,(\bet^{\re})^{(i)},(\bet^{\ima})^{(i)})\right)
,\label{e44om}
\end{equation}
where $L\in|-(K_Y+E)|$ is real with $\R
L\subset F$, $$n = R_Y(D, \bet),\quad n_i=R_Y (\D^{(i)},
(\bet^{\re})^{(i)}+2(\bet^{\ima})^{(i)}),\ i=1,...,m\ ,
$$ \begin{equation}\eta
(l)=\begin{cases}1,\quad & \text{if}\ l=0,\\ l/2+1,\quad
&\text{if}\ l\ \text{is even},\ L',L''\ \text{are imaginary},\\
0,\quad & \text{if}\ l\ \text{is odd},\ L',L''\ \text{are
imaginary}\\ &\quad\text{or are real with}\
\R L',\R L''\subset\overline F_-,\\
(l/2+1)(2-(-1)^{l/2}),\quad &\text{if}\ l\ \text{is even},\ L',L''\
\text{are real},\\ &\quad \R L',\R L''\subset\overline F_+,\\
4([l/4]+1)
(-1)^{L'_{1/2}\circ\varphi},\quad &\text{if}\ l\ \text{is odd},\ L',L''\ \text{are real},\\
&\quad \R L',\R L''\subset\overline F_+,\\ (1+(-1)^{l/2})/2, \quad &\text{if}\ l\
\text{is even},\ L',L''\ \text{are real},\\ &\quad \R
L'\subset\overline F_\pm,\ \R L''\subset\overline F_\mp,\\
2([l/4]+1)(-1)^{(l-1)/2+L'_{1/2}\circ\varphi},\quad &\text{if}\ l\ \text{is odd},\ L',L''\ \text{are real},\\
&\quad \R L'\subset\overline F_\pm,\ \R L''\subset\overline
F_\mp,\\ (-1)^{l/2}(l/2+1),\quad &\text{if}\ l\ \text{is even},\ L',L''\ \text{are real},\\
&\quad \R L',\R L''\subset\overline F_-, \end{cases} \label{ee6}
\end{equation}
and the second sum
in~(\ref{e44om}) is taken
\begin{itemize}\item over all integers $l\ge0$ and vectors
$\alp^{(0)}\le\alp$, $\bet^{(0)}\le\bet^{\re}$;
\item over all
sequences
\begin{equation}(\D^{(i)},\alp^{(i)},(\bet^{\re})^{(i)},(\bet^{\ima})^{(i)}),\
1\le i\le m \ ,\label{e45om}
\end{equation} such that, for all $i=1,...,m$,
\begin{enumerate}
\item[{\rm (3a)}] $\D^{(i)}\in\PicPPR(Y,E)$, and $\D^{(i)}$ is neither the
divisor class $-(K_Y+E)$, nor the pair $\{-(K_Y+E),-(K_Y+E)\}$,
\item[{\rm (3b)}]
$I(\alp^{(i)}+(\bet^{\re})^{(i)}+2(\bet^{\ima})^{(i)})=[\D^{(i)}]E$,
and $R_Y(\D^{(i)}, (\bet^{\re})^{(i)}+2(\bet^{\ima})^{(i)})\ge0$,
\item[{\rm (3c)}] $\D^{(i)}$ is a pair of divisor classes
if and only if
$(\bet^{\ima})^{(i)} \ne 0$,
\item[{\rm (3d)}] if $\D^{(i)}$ is a pair of divisor classes, then
$n_i=0$, $\alp^{(i)}=(\bet^{\re})^{(i)}=0$,
and $(\bet^{\ima})^{(i)}=e_1$,
\end{enumerate}

and
\begin{enumerate}
\item[{\rm (3e)}]
$D-E=\sum_{i=1}^m[\D^{(i)}]- (l+I\alp^{(0)}+I\bet^{(0)})(K_Y+E)$,
\item[{\rm (3f)}] $\sum_{i=0}^m\alp^{(i)}\le\alp$, $\sum_{i=0}^m(\bet^{\re})^{(i)}\ge\bet$,
\item[{\rm (3g)}] each
tuple $(\D^{(i)},0,(\bet^{\re})^{(i)},(\bet^{\ima})^{(i)})$ with
$n_i=0$ appears in (\ref{e45o}) at most once,
\end{enumerate}
\item over all sequences
\begin{equation}\gam^{(i)}\in\Z_+^{\infty,\;
\odd\cdot\;\even},\quad i=1,...,m\ ,\label{e41nom}\end{equation}
satisfying
\begin{enumerate}
\item[{\rm (3h)}] $\|\gam^{(i)}\|=\begin{cases}1,\ &\D^{(i)}\ \text{is a divisor class}, \\
0,\ &\D^{(i)}\ \text{is a pair of divisor classes},\end{cases}\quad
i=1,...,m$,
\item[{\rm (3i)}] $(\bet^{\re})^{(i)}\ge\gam^{(i)}$,
$i=1,...,m$, and $\sum_{i=1}^m\left((\bet^{\re})^{(i)}-
\gam^{(i)}\right)=\bet^{\re}-\bet^{(0)}$;
\end{enumerate}
\end{itemize}
the second sum in (\ref{e44om}) is factorized by simultaneous
permutations in the sequences (\ref{e45om}) and (\ref{e41nom}).

{\rm (4)} All sided $w$-numbers $W_{Y,F_+,\varphi}^\pm(D,\alp,\bet,0)$,
where $D\in\PicPPR(Y,E)$ is a divisor class
and $R_Y(D, \bet) > 0$, are recursively determined
by the formulas (\ref{ee3}), (\ref{ee4}), (\ref{e44o}), (\ref{e44om}) and
the initial conditions given by Proposition~\ref{ini2}.
\end{theorem}

{\bf Proof}.
We follow the main lines of the proof of the recursive
formula in \cite[Section 3]{IKS7}.

{\it Proof of {\rm (1)}}.
The
statement
is clear for
$\alp\not\in\Z_+^{\infty,\;\even}$ or
$\bet^{\re}\not\in\Z_+^{\infty,\;\even}$, since the curves in count
must have a non-empty one-dimensional part in $F_-$ contrary to the
definition of
$V_{Y,F_+}^\R(D,\alp,\bet^{\re},\bet^{\ima},\bp^\flat,\bp^\sharp)$.
In the case $\bet^{\ima} \ne 0$, the statement follows
from Lemma~\ref{imaginary}.

\smallskip

{\it Proof of {\rm (2i)}.} If $(K_Y+E)D=0$, then either
$D^2=-1$, $DE=1$, or $D=-(K_Y+E)$, $DE=2$, and in both
situations, $V_{Y,F_+}^\R(D,\alp,\bet,0,\bp^\flat)=\emptyset$.
Indeed, in the former case, we have $DE=-DK_Y=1$; in the latter
case, the condition $R_Y(-(K_Y+E),\bet)>0$ yields $\bet=2e_1$, and
both conclusions contradict the assumption
$\bet\in\Z_+^{\infty,\;\even}$.

Let $(K_Y+E)D<-2$. The leaf-curves from
$V_{Y,F_+}^\R(D,\alp,\bet,0,\bp^\flat,\bp^\sharp)$ generated by any
reducible root-curve $C=E\cup \check C$ do not contribute
to $W_{Y,F_+,\varphi}^+(D,\alp,\bet,0,\bp^\flat,\bp^\sharp)$.
Indeed, $\check C$ must contain a reduced real primary component, since
$(K_Y+E)(D-E)<0$ and $K_Y+E$ vanishes on all imaginary or
non-reduced primary components of $\check C$ ({\it cf.} Lemma \ref{l5}). Hence, the
total contribution of the leafs $\widetilde C\in\Leaf(C)$ is zero in view of the factor $\mu^+(C,q)=0$ (see formula
(\ref{e2204})) in (\ref{e2213}). Thus, by Lemma \ref{l4},
\begin{equation}W_{Y,F_+,\varphi}^+(D,\alp,\bet,0,\bp^\flat,\bp^\sharp)
=2^{\|\bet\|}W_{Y,F_+,\varphi}^+(D,\alp+\bet,0,0,(\bp^\flat)',(\bp^\sharp)')\
.\label{e2214}\end{equation} However, $R_Y(D,0)=-(K_Y+E)D-1>0$, that
is $(\bp^\sharp)'\ne\emptyset$, and, as explained above, the
right-hand side of (\ref{e2214}) must vanish.

\smallskip

{\it Proof of {\rm (2ii).}} Notice that $D-E$ is not
effective, since $-(K_Y+E)(D-E)=-1$ and $-(K_Y+E)$ is nef. Hence,
there are only irreducible root-curves, and the formula follows from Lemma
\ref{l4}.

\smallskip

{\it Proof of Proposition \ref{ll12} and assertions {\rm (2iii)} and {\rm (3).}} All these statements follow by induction on
$n=R_Y(D,\bet)$ from Lemmas
\ref{l4} and \ref{sum-leaves}. Proposition \ref{ini2} serves as the base
of induction.

In the right-hand
side of formulas (\ref{e44o}) and (\ref{e44om}), the first sum runs over irreducible root-curves and the second
sum runs over root-curves containing $E$. We only explain notations and coefficients occurring in the second sum:
\begin{itemize}\item the vector $\alp^{(0)}$ encodes the multiplicities of the primary components of
type $iL(p_{ij})$, \item the vector $\bet^{(0)}$ encodes the
multiplicities of the primary components of type $sL(z)$, and their multiplicative contribution amounts to $4^{\|\bet^{(0)}\|}$,
\item the factors
$(l/2 + 1)$
and $\eta
(l)$ are the sums of the contributions of
primary components
$sL'$, $s''L''$
(computed
using (\ref{e2212}) and (\ref{e2212a})) over the range
$s'+s''=l$,
\item the vectors $\gam^{(i)}$ encode
the intersection multiplicities $j=(\check C\cdot E)_q$ of the points
$q\in\bz'$ (see Lemma \ref{l6}) for the reducible root-curves $C=E\cup \check C$.
\end{itemize}

The statement {\rm (4)} is straightforward.
\proofend

\section{ABV formula over
the reals}\label{sec-vakil}
Let $Y$ be a smooth rational surface, $E\subset Y$ a smooth rational
curve. If the anti-canonical class $-K_Y $ is effective, positive on
all curves different from $E$, and $K_YE=0$, we call the pair
$(Y,E)$ a {\it nodal del Pezzo pair}. (It follows from the
adjunction formula that $E^2=-2$.) Notice that a nodal del Pezzo
pair may be not monic log-del Pezzo, and vice versa. Throughout this
section we assume that $(Y,E)$ is a nodal del Pezzo pair.

An example of a nodal del Pezzo pair is provided by
the plane
blown up at a generic collection of $\leq 8$ points subject to the
condition that
six of them
belong to a
conic.

A
nodal del Pezzo pair $(Y, E)$ is an almost Fano surface in the sense
of \cite[Section 4.1]{Va}, and thus by \cite[Theorem 4.2]{Va} we
have the following
Abramovich-Bertram-Vakil formula (briefly {\it ABV formula}):
\begin{equation}GW_0(Y,D)=\sum_{m\ge0}\binom{DE+2m}{m}N_Y(D-mE,0,(DE+2m)e_1)\ ,
\label{ee1}\end{equation} where $D\in\Pic(Y)$ is a divisor class, $GW_0(Y,D)$ states for the
genus zero Gromov-Witten invariant,
and $N_Y(D',0,(D'E)e_1)$ is the number of
rational curves $C\in|D'|$ passing through a generic
collection of $-D'K_Y-1$
points in $Y\setminus E$.

\subsection{Deformation representation
of ABV formula}\label{sec8}

ABV formula (\ref{ee1}) can be represented geometrically. Let
$\pi:{\mathfrak X}\to
(\C,0)$ be a proper
holomorphic submersion of a smooth three-dimensional variety ${\mathfrak X}$  (with $(\C,0)$ being understood as a disc germ),
where
each fiber ${\mathfrak X}_t, t \ne 0$, is a del Pezzo
surface
and the central fiber $Y
={\mathfrak X}_0$ contains a smooth rational curve $E$ such that
$(Y,E)$ is a nodal del Pezzo pair.

\begin{remark}\label{r3} There is a natural isomorphism
\begin{equation}\Pic({\mathfrak
X}_t)\overset{\simeq}{\longrightarrow}\Pic(Y), \quad t\ne0\
,\label{e2072}\end{equation} preserving the intersection form; for
the sake of
brevity we use the same symbol for corresponding
divisor classes in ${\mathfrak X}_t$, $t\in
(\C,0)$.
To distinguish
linear systems themselves
we use the notation $|D|_{{\mathfrak X}_t}$.
\end{remark}

Let $D\in\Pic({\mathfrak X}_t)$ be effective for $t\ne0$ and
satisfy
$-K_{{\mathfrak X}_t}D-1\ge0$. Pick $r=-K_{{\mathfrak X}_t}D-1$ disjoint
sections $z_i:(\C,0)\to{\mathfrak X}$, $1\le i\le r$, so that
$\bp^\sharp(0)=\{z_i(0),\ 1\le i\le r\}$ is a generic
point collection
in $Y\setminus E$, and $\bp^\sharp(t)=\{z_i(t),\ 1\le i\le r\}$ is a
generic
point collection
in ${\mathfrak X}_t$, $t\ne0$. For each
$t\in(\C,0)$, denote by $V_t(D,\bp^\sharp(t))$ the set of reduced
irreducible rational curves $C\in|D|_{{\mathfrak X}_t}$
that
pass
through $\bp^\sharp(t)$. It is well-known (see, for instance,
\cite{GP}) that $V_t(D,\bp^\sharp(t))$, $t\ne0$, is finite, contains
$GW_0(Y, D)$ elements, and each element is a nodal curve. By \cite[Proposition 2.1]{MS2}, for each
$m\ge0$, the set $V_0(D-mE,\bp^\sharp(0))$ is finite, and its
elements are immersed curves crossing $E$ transversally at $DE+2m$
distinct points. Thus, we have a diagram
\begin{equation}\begin{matrix}\widetilde{\mathcal C}' & \overset{\nu'}{\longrightarrow} &
{\mathcal C}' & \hookrightarrow & {\mathfrak X} \\
\;\;\downarrow\widetilde\pi' & & \downarrow & & \;\;\;\downarrow\pi \\
(\C,0)\setminus \{0\}
& = &
(\C,0)\setminus \{0\} & \hookrightarrow &
(\C,0)
\end{matrix}\label{ee10}\end{equation} where ${\mathcal C}'$ is the union of
$GW_0(Y,D)$ families of curves $C\in
V_t(D,\bp^\sharp(t))$, $t\in(\C,0)\setminus\{0\}$, and $\widetilde
{\mathcal C}'$ is its normalization. The following statement
follows from \cite[Theorem 4.2]{Va}.

\begin{proposition}\label{t6}
There exists a diagram
\begin{equation}\begin{matrix}\widetilde{\mathcal C} &
\overset{\nu}{\longrightarrow} & {\mathcal C} & \hookrightarrow & {\mathfrak X}\\
\;\;\downarrow\widetilde\pi & & \downarrow & & \;\;\;\downarrow\pi\\
(\C,0) & = & (\C,0) & = &
(\C,0)\end{matrix}\label{e2071}\end{equation} which extends the
diagram (\ref{ee10}) so that
\begin{enumerate}\item[{\rm (1)}]
\begin{itemize}\item ${\mathcal C}$ is the closure of ${\mathcal C}'$ in
${\mathfrak X}$; \item the fiber over $0$ of each component of
${\mathcal C}$ is $C_0\cup mE$ for some $m\ge 0$, where $C_0\in
V_0(D,\bp^\sharp(0))$; \item each curve $C_0\cup mE$ with $m\ge0$,
$C_0\in V_0(D-mE,\bp^\sharp(0))$ appears as the fiber over $0$ for
exactly $\binom{DE+2m}{m}$ components of ${\mathcal
C}$;\end{itemize}
\item[{\rm (2)}] \begin{itemize}\item $\widetilde{\mathcal C}$ is the union of $GW_0(Y,D)$ disjoint nonsingular
surfaces; \item the fiber over $0$ of each component of
$\widetilde{\mathcal C}$ is either isomorphic to $\PP^1$ with
$\nu:\PP^1\to C_0\in V_0(D,\bp^\sharp(0))$ birational, or is a nodal
reducible rational curve $\bigcup_{i=0}^m\PP^1_{(i)}$ with some
$m\ge1$, $\PP^1_{(i)}\simeq\PP^1$ for all $i=0,...,m$,
$\PP^1_{(1)},...,\PP^1_{(m)}$ disjoint from each other,
$\PP^1_{(0)}$ intersecting each $\PP^1_{(1)},...,\PP^1_{(m)}$ at one
point, and such that $\nu:\PP^1_{(0)}\to C_0\in
V_0(D-mE,\bp^\sharp(0))$ is birational, $\nu:\PP^1_{(i)}\to E$ is an
isomorphism for all $i=1,...,m$; \item for each $C_0\in
V_0(D-mE,\bp^\sharp(0))$, $m\ge0$, there are exactly
$\binom{DE+2m}{m}$ components of $\widetilde{\mathcal C}$ whose
fiber $\bigcup_{i=0}^m\PP^1_{(i)}$ over $0$ covers $C_0$, and they
differ from each other by the image of the $m$-tuple
$\PP^1_{(0)}\cap(\bigcup_{i=1}^m\PP^1_{(i)})$ in the
\mbox{$(DE+2m)$}-tuple $C_0\cap E$.\end{itemize}
If the families ${\mathfrak X}$, ${\mathcal C}'$, and $\widetilde{\mathcal C}'$
are defined over the reals, then so are the families
${\mathcal C}$ and $\widetilde{\mathcal C}$.\proofend
\end{enumerate}
\end{proposition}

\subsection{Nodal degenerations}\label{sec-nod}
Let $\pi':{\mathfrak X}'\to(\C,0)$ be a holomorphic map of a smooth three-dimensional variety
${\mathfrak X}'$ such that
the
fibers ${\mathfrak X}'_\tau$,
$\tau\in(\C,0)\setminus\{0\}$, are del Pezzo surfaces,
the central fiber
${\mathfrak X}'_0$ is a
surface with one singular point $z$ of type
$A_1$ (node),
at the point $z$ the map $\pi'$ is given, with respect to appropriate local
coordinates $x_1,x_2,x_3$,
by
\begin{equation}
\pi'(x_1,x_2,x_3)=a_1x_1^2+a_2x_2^2+a_3x_3^2,
\quad
a_1a_2a_3\ne0
\ ,\label{ee11}\end{equation}
and $\pi' $ is a submersion
at each point of ${\mathfrak X}'\setminus \{z\}$. Such a family is called {\it nodal degeneration}.

Make the base change $\tau=t^2$, perform the blow up $\widetilde{\mathfrak X}'\to{\mathfrak X''}$
at the node of the new family,
$\mathfrak X''=
(\C, 0) \times_{t^2=\pi'}\mathfrak X'$,
and obtain a
family $\widetilde\pi':\widetilde{\mathfrak X}'\to(\C,0)$, whose
fibers $\widetilde{\mathfrak X}'_t$, $t\ne0$, are del Pezzo
surfaces, and $\widetilde{\mathfrak X}'_0=Y\cup Z$, where
$Z\simeq(\PP^1)^2$, $E=Y\cap Z$ is a smooth rational $(-2)$-curve in
$Y$, and $(Y,E)$ is a nodal del Pezzo pair. Here, $E$ represents the
class $C_1+C_2$ in $\Pic(Z)$, $C_1,C_2$ being the generators of the
two rulings of $Z$. We call the family $\widetilde\pi':\widetilde{\mathfrak X}'\to(\C,0)$
the {\it unscrew} of the nodal degeneration $\pi':{\mathfrak X}'\to(\C,0)$.

Contracting $Z$ to $E$ along the lines
of one of the rulings (see \cite{A}), say, generated by $C_1$, we get a family
\begin{equation}\pi:{\mathfrak X}\to(\C,0)\label{ee12}\end{equation}
of smooth surfaces as in Section \ref{sec8}.
Let $D\in\Pic({\mathfrak X}_t)$ be effective for $t\ne0$ and
satisfy
$-K_{{\mathfrak X}_t}D-1\ge0$. Pick $r=-K_{{\mathfrak X}_t}D-1$ disjoint
sections $z_i:(\C,0)\to{\mathfrak X}$, $1\le i\le r$, so that
$\bp^\sharp(0)=\{z_i(0),\ 1\le i\le r\}$ is a generic
point collection
in $Y\setminus E$, and $\bp^\sharp(t)=\{z_i(t),\ 1\le i\le r\}$ is a
generic
point collection
in ${\mathfrak X}_t$, $t\ne0$.
The following lemma is straightforward.

\begin{lemma}\label{limits}
Let ${\mathcal C}'$ be
a family of rational curves ${\mathcal C}'_t\in |D|_{{\mathfrak X}_t}$
which
pass
through $\bp^\sharp(t)$,
$t\in(\C,0)\setminus\{0\}$.
The family
${\mathcal C}'$
lifts
to a family of curves in
$\widetilde\pi':\widetilde{\mathfrak X}'\to(\C,0)$ as follows: if
${\mathcal C}'$
closes up in ${\mathfrak X}$ with the central fiber $C'_m\cup mE$,
$C'_m\in V_0(D-mE,\bp^\sharp(0))$,
then
${\mathcal C}'$ closes up
in $\widetilde{\mathfrak X}'$ with a central fiber
$C'_m\cup C^{(DE+m)}_1\cup C_2^{(m)}$, where
$C^{(m)}_2$ is the union of $m$ lines in $|C_2|_Z$
attached to $m$ intersection points of $C'_m$ and $E$, and $C_1^{(DE+m)}$ is the union of
$DE+m$ lines from $|C_1|_Z$ attached to the remaining points of
$C'_m\cap E$.\proofend
\end{lemma}

\begin{remark}\label{r4}
A family of plane quartics with the central fiber $Q$ having one node
$z$ induces a family of del Pezzo surfaces of degree $2$ degenerating into a nodal
del Pezzo pair. In this setting, $E$ is the exceptional divisor of
the blow up of the node of the double cover of the plane ramified
along the nodal quartic, the six pairs of intersecting $(-1)$-curves
crossing $E$ respectively cover the six lines in the plane passing
though $z$ and tangent to $Q$ outside $z$, and, finally, the supporting curves
$L',L''$ doubly cover the tangent lines to $Q$ at the node $z$.
\end{remark}

Assume that a nodal degeneration $\pi':{\mathfrak X}'\to(\C,0)$
possesses a real structure $\Conj$ which
lifts the standard complex conjugation $\cconj:(\C,0)\to(\C,0)$.
Then, the point $z$ is real and with respect to appropriate real local coordinates at $z$ the map $\pi'$ is given by
\begin{equation}
\pi'(x_1,x_2,x_3)=a_1x_1^2+a_2x_2^2+a_3x_3^2,
\quad a_1, a_2, a_3 \in \R, \quad
a_1a_2a_3\ne0.\label{ee11_new}\end{equation}
The real structure $\Conj$
gives rise
to two real structures
on
the unscrew $\widetilde\pi':\widetilde{\mathfrak X}'\to(\C,0)$ of $\pi':{\mathfrak X}'\to(\C,0)$.
One real structure covers the complex conjugation
$t \mapsto \overline t$,
and we call the resulting real unscrew a $\theta${\it -unscrew}, where $\theta$ is the signature of
the quadratic form (\ref{ee11_new}).
The quadric $Z\subset\widetilde{\mathfrak X}'_0$ is real, and $$\R Z\simeq\begin{cases}S^2,\quad
 &\text{if}\ \theta = 3 \ \text{or}\ -1,\\
 (S^1)^2,\quad &\text{if}\ \theta = 1,\\
 \emptyset,\quad &\text{if}\ \theta = -3.\end{cases}$$
The other real structure on $\widetilde\pi':\widetilde{\mathfrak X}'\to(\C,0)$
covers the conjugation $t \mapsto -\overline t$
and defines a {\it mirror $(-\theta)$-unscrew}.

\begin{proposition}\label{new-proposition}
Let $\widetilde\pi':\widetilde{\mathfrak X}'\to(\C,0)$ be a $\theta$-unscrew.
Then, the isomorphism $\Pic({\mathfrak X}_t) \overset{\sim}{\to} \Pic(Y)$
is conjugation invariant. If $\theta$ is equal to $1$ or $-3$, this isomorphism induces an isomorphism
$\Pic^\R({\mathfrak X}_t)\overset{\sim}{\to}\Pic^\R(Y)$,
$t\in(\R,0)$.
If $\theta$ is equal to $-1$ or $3$, the isomorphism $\Pic(\widetilde{\mathfrak X}'_t) \overset{\sim}{\to} \Pic(Y)$, $t\ne0$,
induces a monomorphism $\Pic^\R(\widetilde{\mathfrak X}'_t)\to\Pic^\R(Y)$, $t\in(\R,0)\setminus\{0\}$,
and the image of the latter is orthogonal to $[E]\in \Pic^\R(Y)$.
\end{proposition}

{\bf Proof}.
If $\theta$ is equal to $1$ or $-3$,
the action of  $\Conj$ in  $\Pic(\widetilde{\mathfrak X}'_t)=H^2(\widetilde{\mathfrak X}'_t; \Z)$
and $\Pic(Y)=H^2(Y;\Z)$ commutes
with the natural (as in Remark \ref{r3}) isomorphism $H^2(\widetilde{\mathfrak X}'_t; \Z)\to H^2(Y;\Z)$.
If $\theta$ is equal to $-1$ or $3$,
the action of  $\Conj$ in  $\Pic(\widetilde{\mathfrak X}'_t)$
and $\Pic(Y)$ does not commute
with the isomorphism $H^2(\widetilde{\mathfrak X}'_t; \Z)\to H^2(Y;\Z)$,
the defect being the composition with the reflection
in $[E]\in H^2(Y;\Z)$.
\proofend

\subsection{Real versions of ABV formula}\label{sec9}

\subsubsection{Ordinary and sided $u$-numbers}\label{unum}
Let $(Y, E)$ be a nodal del Pezzo pair such that
$Y$ and $E$ are real, and $\R E\ne\emptyset$. Denote by $F$ the
connected component of $\R Y$ containing $\R E$ and pick a
conjugation invariant
class $\varphi\in H_2(Y\setminus F,\Z/2)$.
Let
$D\in\PicPPR(Y,E)$.
Choose a generic collection
$\bp^\sharp$ of $-K_YD-1$
points
in $F\setminus E$.

By \cite[Proposition 4.1(b)]{Va}, the set $V_Y(D,
\bp^\sharp)$ of rational curves in the linear system $|D|$ which
pass through the points of $\bp^\sharp$
is finite and consists of immersed curves crossing $E$ transversally
at $DE$ distinct points.

For any nonnegative integers $k$ and $l$ such that $k + 2l = DE$,
define an {\it ordinary $u$-number} $U_{Y,F,\varphi}(D, ke_1, le_1,
\bp^\sharp)$ putting ({\it cf.} the definition of ordinary $w$-numbers in Section \ref{new-ordinary})
\begin{equation}
U_{Y,E,\varphi}(D, ke_1, le_1, \bp^\sharp) =\sum_{C\in V^\R_Y(D,
ke_1, le_1, \bp^\sharp)}\mu_\varphi(C)
\ ,\label{eu1}
\end{equation}
where $V^\R_Y(D, ke_1, le_1, \bp^\sharp) \subset V_Y(D,
\bp^\sharp)$ is formed by the curves intersecting $\R E$ in $k$
points (and intersecting $E \setminus \R E$ in $l$ pairs of complex
conjugate points) and $\mu_\varphi(C)$ is defined by
(\ref{e2203}).
If
$F \setminus E$ splits into two components $F_+$ and $F_-$,
the configuration
$\bp^\sharp$ lies in $F_+$, and $DE$ is even, then
define a {\it sided $u$-number} $U^{\pm}_{Y,F_+,\varphi}(D, 0,
(DE/2)e_1, \bp^\sharp)$,
putting
\begin{equation}
U_{Y,F_+,\varphi}^+(D, 0, (DE/2)e_1, \bp^\sharp) =\sum_{C\in
V^\R_{Y,F_+}(D, 0, (DE/2)e_1, \bp^\sharp)}\mu_\varphi(C)
\ ,\label{eu2}
\end{equation}
\begin{equation}
U_{Y,F_+,\varphi}^-(D, 0, (DE/2)e_1, \bp^\sharp) =\sum_{C\in
V^\R_{Y,F_+}(D, 0, (DE/2)e_1, \bp^\sharp)}\mu^-_\varphi(C)
\ ,\label{eu3}
\end{equation}
where
$$
V^\R_{Y, F_+}(D, 0, (DE/2)e_1, \bp^\sharp) = \{C \in V^\R_Y(D, 0,
(DE/2)e_1, \bp^\sharp) \ : \ \card(C \cap F_-) < \infty\} \
$$ and $\mu^-_\varphi(C)$ is defined by (\ref{e2202}).
\smallskip

We say that a quadruple $(Y, E, F, \varphi)$ has {\it
property}~($R$) if for any
divisor class $D\in\PicPPR(Y,E)$
and for any connected component $F_+$ of $F \setminus E$, there
exists a generic collection
$\bp^\sharp$ of $-K_YD-1$ points in $F_+$ (referred to as {\it
$R_{D, F_+}$-collection} or {\it $R_D$-collection}) such that, for
any $m\ge0$ with $D-mE\in\PicPPR(Y,E)$,
the following holds:
\begin{enumerate}
\item[(R1)] $U_{Y,F,\varphi}(D-mE, ke_1, le_1, \bp^\sharp) = 0$
whenever $l>0$,
\item[(R2)] if
$F\setminus E$ splits into two components
and the intersection $DE$ is even,
then $U^+_{Y,F_+,\varphi}(D-mE,0,(DE/2)e_1, \bp^\sharp)=0$.
\end{enumerate}

\begin{proposition}\label{u=w} Let $(Y, E)$ be a nodal del Pezzo pair such that
$Y$ and $E$ are real, and $\R E\ne\emptyset$. Denote by $F$ the
connected component of $\R Y$ containing $\R E$ and pick a
conjugation invariant
class $\varphi\in H_2(Y\setminus F,\Z/2)$.
Assume in addition
that $(Y, E)$ is monic log-del Pezzo.

{\rm (1)} Pick a divisor class $D_0 \in \PicPPR(Y, E)$ such that $\dim
|D_0| > 0$. Let $D \in \Prec(D_0)$ be a divisor class, and let
$\bp^\sharp$ be a collection of $-K_Y D - 1$ points in $D_0$-CH
position in $F \setminus E$. Then, for any nonnegative integers~$k$
and~$l$ such that $k + 2l = DE$, one has
$$
U_{Y, E, \varphi}(D, ke_1, le_1, \bp^\sharp) = W_{Y, E, \varphi}(D,
0, ke_1, le_1, \emptyset, \bp^\sharp).
$$
If $F \setminus E$ splits into two components $F_+$ and $F_-$, the
collection $\bp^\sharp$ is contained in $F_+$, and $DE$ is even,
then
$$
U^\pm_{Y, F_+, \varphi}(D, 0, (DE/2)e_1, \bp^\sharp) = W^\pm_{Y,
F_+, \varphi}(D, 0, 0, (DE/2)e_1, \emptyset, \bp^\sharp).
$$
The numbers $U_{Y, E, \varphi}(D, ke_1, le_1, \bp^\sharp)$ and
$U^\pm_{Y, F_+, \varphi}(D, 0, (DE/2)e_1, \bp^\sharp)$ do not depend
on the choice of a collection $\bp^\sharp$ in $D_0$-CH position.

{\rm (2)} The quadruple $(Y, E, F, \varphi)$ has property (R).
\end{proposition}

{\bf Proof}. The equality of ordinary (respectively, sided) $u$- and
$w$-numbers is tautological. The invariance of the $u$-numbers
considered follows from the invariance of $w$-numbers; for the
latter invariance see Propositions~\ref{l12} and~\ref{ll12}.
\proofend

\subsubsection{External $u$-numbers}\label{eun}

Let $(Y, E)$ be a nodal del Pezzo pair such that
$Y$ and $E$ are real, and $\R Y\ne\emptyset$. Let $F\subset\R Y$ be a connected component
such that $F\cap\R E=\emptyset$. Let
$D\in\PicPPR(Y,E)$.
Choose a generic collection
$\bp^\sharp$ of $-K_YD-1$
points
in $F$. Notice that $DE\ge0$ is even, and each real curve in $V_Y(D,
\bp^\sharp)$ intersects $E$ at $DE/2$ distinct pairs of complex conjugate points.

If $\R E\ne\emptyset$, we denote by $F'$ the connected component of $\R Y$ containing
$\R E$. If $\R E=\emptyset$, we put $F'=\emptyset$.
Choose a conjugation invariant class $\varphi\in H_2(Y\setminus(F\cup F'),\Z/2)$ and
define {\it external $u$-numbers}
$$U_{Y,E,\varphi'}(D,\bp^\sharp)=\sum_{C\in V^\R_Y(D, 0, (DE/2)e_1, \bp^\sharp)}\mu_{\varphi'}(C),\quad \text{for}\ \varphi'
=\varphi\ \text{or}\ \varphi+
[F']\ .$$

\subsubsection{ABV formulas for Welschinger invariants, I}\label{sec9I}
As in Section \ref{sec8}, let
$\pi:{\mathfrak X}\to(\C,0)$ be
holomorphic submersion of a smooth three-dimensional variety ${\mathfrak X}$,
where
each fiber ${\mathfrak X}_t, t \ne 0$, is a del Pezzo
surface
and the central fiber $Y
={\mathfrak X}_0$ contains a smooth rational curve $E$ such that
$(Y,E)$ is a nodal del Pezzo pair.

Suppose that ${\mathfrak
X}$ possesses a real structure $\conj: {\mathfrak X} \to {\mathfrak
X}$
such that
\begin{equation}\pi \circ \conj = \cconj \circ \pi\label{ers}\end{equation}
(where $\cconj$ is the standard real structure on $(\C, 0)$). We get
a family
$\pi:{}
\R{\mathfrak X}=\pi^{-1}(\R,0)\to(\R,0)\hookrightarrow(\C,0)$
of real surfaces, the fibers ${\mathfrak X}_t$,
$t\in(\R,0)\setminus\{0\}$, being real del Pezzo surfaces, and
$(Y,E)$ being a real nodal del Pezzo pair.
Assume that $\R E\ne\emptyset$. Denote by $F_0$ the
connected component of $\R Y$ containing $\R E$ and pick a
conjugation invariant
class $\varphi_0\in H_2(Y\setminus F_0,\Z/2)$.
The family $\R{\mathfrak
X}\to(\R,0)$ is topologically trivial. We extend $F_0$
to a continuous family $F_t$ of
connected components of $\R{\mathfrak X}_t$, and extend $\varphi_0$ to a continuous family of
conjugation invariant classes
$\varphi_t\in H_2(\R{\mathfrak
X}_t \setminus F_t, \Z/2)$, $t\in(\R,0)$.

\begin{theorem}\label{t2}
{\rm (1)} For any
real effective
divisor class $D$ on ${\mathfrak X}_t$, $t\ne0$, one has
$$W({\mathfrak X}_t,D,F_t,\varphi_t)=W({\mathfrak
X}_t,D+(DE)E,F_t,\varphi_t)\ .$$

{\rm (2)} Assume that the quadruple $(Y, E, F_0, \varphi_0)$ has property (R).
Then, for any $t\in(\R,0)$, $t\ne0$, any
divisor class $D\in\Pic^\R({\mathfrak X}_t)$, and
any $R_D$-collection
$\bp^\sharp \subset F_0 \setminus E$,
the following equality holds:
\begin{equation}
W({\mathfrak X}_t,D,F_t,\varphi_t)
=\sum_{m\ge0}(-1)^{m(E_{1/2}\circ\varphi_0)}\binom{DE+2m}{m}U_{Y,E,\varphi_0} (D-mE,(DE+2m)e_1,0,
\bp^\sharp) \ . \label{ee2}
\end{equation}
\end{theorem}

{\bf Proof}.
For the
first statement, without loss of generality, assume that
$DE=-d<0$ and choose a continuous family of
collections $\bp^\sharp_t \subset F_t$ of $-K_XD-1$ distinct points
so that
$\bp^\sharp_0$ is generic in $F_0\setminus E$. We establish a
one-to-one correspondence between the sets $M_1$ and $M_2$ of real
rational curves in $|D|_{{\mathfrak X}_t}$
and $|D-dE|_{{\mathfrak X}_t}$, respectively, passing through
$\bp^\sharp_t$, such that the correspondence preserves the Welschinger signs.
Indeed, by Proposition \ref{t6}, degenerations of curves from $M_1$
are of type $C\cup(d+s)E$, $s\ge0$, where $C\in|D-(d+s)E|_Y$ is a real
rational curve passing through
$\bp^\sharp_0$; furthermore, to each such a curve $C$
and a conjugation invariant subset $w \subset C \cap E$ of $d + s$
points there corresponds a unique curve in $M_1$, and its
Welschinger sign coincides with that of $C$. Similarly, degenerations of curves from
$M_2$ are of type $C\cup sE$, $C\in|D-(d+s)E|_Y$ as above, and to each
subset $(C \cap E) \setminus w\subset C\cap E$ of $s$ points there corresponds a unique curve in $M_2$, and
its
Welschinger sign coincides with that of $C$.

To prove the
second
statement, notice that, since the quadruple $(Y, E, F_0, \varphi_0)$ has
property (R), the degenerate real curves to consider are of type $C\cup mE$, $m\ge0$, where
$C \in V^\R_Y(D - mE, (DE+2m)e_1, 0, \bp^\sharp)$, and the
Welschinger sign of each real rational curve in ${\mathfrak
X}_t$, $t\ne0$, appearing in a deformation of
$C \cup mE$
is $(-1)^{m(E_{1/2}\circ\varphi)}\mu_\varphi(C)$. Thus,
formula (\ref{ee2}) follows from Proposition~\ref{t6}.
\proofend

\subsubsection{ABV formulas for Welschinger invariants, II}\label{sec9II}
Assume that
$\widetilde\pi':\widetilde{\mathfrak X}'\to(\C,0)$ is a real unscrew
of a real nodal degeneration $\pi':{\mathfrak X}'\to(\C,0)$, and
$\R E\ne\emptyset$ (see Section \ref{sec-nod}).
In this case, the signature $\theta$ of the unscrew is equal either to $1$, or to $-1$.
If $\theta = 1$ (respectively, $\theta = -1$), one has $\R Z\simeq(S^1)^2$
(respectively, $\R Z\simeq S^2$).
Let $F\subset\R Y$ be the connected component containing $\R E$.
Pick a conjugation invariant
class $\varphi\in H_2(Y\setminus F,\Z/2)$.

In the case of a
$1$-unscrew,
both rulings of $Z$ are real.
The contraction of $Z$ along one of the rulings leads to
a family of smooth real surfaces $\pi:{\mathfrak X}\to(\C,0)$ ({\it cf.} Section
\ref{sec-nod}).
Thus, Theorem \ref{t2}
applies and gives Welschinger invariants of the real del Pezzo
surfaces ${\mathfrak X}_t$, $t\in(\R,0)\setminus\{0\}$, via
$w$-numbers of the real nodal del Pezzo pair $(Y,E)$ by formula
(\ref{ee2}).

\begin{theorem}\label{t4} Assume that the unscrew $\widetilde\pi':\widetilde{\mathfrak X}'\to(\C,0)$
of a real
nodal degeneration $\pi':{\mathfrak X}'\to(\C,0)$ is of signature $-1$. Then, the following holds.

{\rm (1)} Suppose that $F\setminus E$ is connected.
Let $F_t$ be a component of
$\R\widetilde{\mathfrak X}'_t$ merging to $F$ as $t\to0$.
Choose a class $\varphi\in H_2(\R Y\setminus F,\Z/2)$,
and denote by $\varphi_t$ the
class in $H_2(\R\widetilde{\mathfrak X}'_t\setminus F_t,\Z/2)$,
$t\in(\R,0)\setminus\{0\}$, which converges to $\varphi$ as $t \to
0$. If the quadruple $(Y, E, F, \varphi)$ has property (R),
then for any $t\in(\R,0)$, $t\ne0$, any divisor class
$D\in\Pic_+^\R(\widetilde{\mathfrak X}'_t)$
and any
$R_D$-collection $
\bp^\sharp \subset F\setminus E$,
one has
\begin{equation}W(\widetilde{\mathfrak
X}'_t,D,F_t,\varphi_t)=U_{Y,E,\varphi}(D,0,0,
\bp^\sharp) \ .\label{e2077}
\end{equation}

{\rm (2)} Suppose that $F\setminus E$ splits into two components
$F_+,F_-$, and $\R\widetilde{\mathfrak X}'_t$ contains two connected
components $F_{+,t}$ and $F_{-,t}$ which merge to $F_+$ and $F_-$,
respectively.
Choose a class $\varphi \in H_2(\R Y\setminus F,\Z/2)$, and denote
by $\varphi_t$ be the
class in $H_2(\R\widetilde{\mathfrak X}'_t\setminus F_t,\Z/2)$,
$t\in(\R,0)\setminus\{0\}$, which converges to $\varphi$ as $t\to0$.
\begin{enumerate}
\item[{\rm (i)}] If
the quadruple $(Y, E, F, \varphi)$ has property (R), then, for
any $t\in(\R,0)$, $t\ne0$, any divisor class $D\in\Pic_+^\R(\widetilde{\mathfrak X}'_t)$,
and any
$R_{D, F_+}$-collection $
\bp^\sharp \subset F_+$,
one has
\begin{equation}
W(\widetilde{\mathfrak X}'_t,D,F_{+,t},\varphi_t) =
U^+_{Y,F_+,\varphi}(D,0,0,\bp^\sharp)\
. \label{ee13}
\end{equation}
\item[{\rm (ii)}]
For any $t\in(\R,0)$, $t\ne0$, any divisor class
$D\in\Pic_+^\R(\widetilde{\mathfrak X}'_t)$,
and any generic configuration $
\bp^\sharp \subset F_+$ of $-K_YD-1$ points, one has
\begin{equation}
W(\widetilde{\mathfrak X}'_t,D,F_{+,t},F_{-,t} \cup \varphi_t) =
U^-_{Y,F_+,\varphi}(D,0,0,
\bp^\sharp) \ . \label{ee14}\end{equation}\end{enumerate}
\end{theorem}

{\bf Proof}.
Pick a divisor class
$D\in\Pic_+^\R(\widetilde{\mathfrak X}'_t)$, $t\ne0$,
take disjoint sections $z_i:(\R,0) \to \R\widetilde{\mathfrak X}'$, $1\le i\le-K_YD-1$,
and consider the
limits
that the
real rational curves in $|D|_{\widetilde{\mathfrak X}'_t}$ passing
through $\{z_i(t)\}_{1\le i\le-K_YD-1}$, $t\ne0$, have in the central fiber
$Y\cup Z$. These limits are of type $C'_m\cup C_1^{(DE+m)}\cup C_2^{(m)}$
(see Lemma \ref{limits}).

To prove statements (1) and (2i), assume that the quadruple $(Y, E, F, \varphi)$ has property (R), and
that $\bp^\sharp=\{z_i(0)\}_{1\le i\le-K_YD-1}\subset F\setminus E$ is an
$R_D$-collection. Due to property (R), the components $C'_m\in|D-mE|_Y$, $m>0$, of the limits
of real rational curves in $|D|_{\widetilde{\mathfrak X}'_t}$ must have real intersection points with $E$.
However, such curves $C'_m$ cannot be completed to real curves in $|D|_Y$. Taking this into account,
we prove statements (2) and (3i)
in the same manner as Theorem \ref{t2}(2).

To prove statement (2ii), put $\bp^\sharp=\{z_i(0)\}_{1\le i\le-K_YD-1}$, and
notice that
if the limit curve contains as a component a real rational
curve $C\in|D-mE|_Y$, $m>0$, then $C\cap\R E=\emptyset$. In the
family $\{\widetilde{\mathfrak X}'_t\}_{t\in[0,1]}$, the surface $\widetilde{\mathfrak
X}'_0=Y\cup Z$
deforms so that $F_+$ glues up with a component $Z_+$ of $\R
Z\setminus E$, whereas $F_-$ glues up with the other component $Z_-$
of $\R Z\setminus E$. In turn, each real rational curve
$C\in|D-mE|_Y$ with $C\cap \R E=\emptyset$ can be completed up to a
real curve on $Y\cup Z$ in $2^m$
ways,
when attaching to each pair $z',z''\in C\cap E$ of complex conjugate
points either the pair $C_1'\supset\{z'\}$, $C''_2\supset\{z''\}$,
or the pair $C''_1\supset\{z''\}$, $C'_2\supset\{z'\}$, where $C'_1,C''_1$ belong to
one ruling of $Z$, and $C'_2,C''_2$ to the other. Observe that one of the pairs $(C'_1,C''_2)$, $(C''_1,C'_2)$
has a solitary node in $Z_+$, which contributes the factor $(-1)$ to
the Welschinger sign $\mu^-_\varphi$ of the corresponding deformed
curve in $|D|_{{\mathfrak X}'_t}$, $t>0$,
whereas the other pair has a solitary node in $Z_-$
that does not affect $\mu^-_\varphi$. Hence, the total contribution
$W(\widetilde{\mathfrak X}'_t,D,F_{+,t},F_{-,t}\cup \varphi_t)$ of the curves
coming from $C$ is zero, which proves formula
(\ref{ee14}).
\proofend

\begin{corollary}\label{c2}
Let $(Y,E)$ be a real nodal del Pezzo pair such that $\R E$ divides
a connected component $F$ of $\R Y$ into two parts, $F_+$ and $F_-$.
Let $D\in\Pic^\R(Y)$ be an effective divisor class such that $DE=0$.
Let $\varphi\in H_2(Y\setminus F,\Z/2)$ be a conjugation invariant class.
Then, the
sided $u$-number $U^-_{Y,F_+,\varphi}(D,0,0,
\bp^\sharp)$ does not depend on the choice of a generic
configuration $\bp^\sharp \subset F_+$ of $-K_YD-1$ points.\proofend
\end{corollary}

Let $\pi':{\mathfrak X}'\to(\C,0)$ be a real nodal degeneration. Assume that
its unscrew is
of signature $1$,
denote
this unscrew
by $\pi^h:{\mathfrak X}^h\to(\C,0)$. Its mirror unscrew
$\pi^e:{\mathfrak X}^e\to(\C,0)$ is
of signature $-1$.
Let $(Y,E)$ be the real nodal del Pezzo pair that appears as
a component of both ${\mathfrak X}^h_0$ and ${\mathfrak X}^e_0$. Let $F$ be the connected component of $\R Y$
containing $\R E$, and let $F^h_t$, respectively, $F^e_t$, be the component
(or the union of components) of $\R{\mathfrak X}^h_t$, respectively,
$\R{\mathfrak X}^e_t$, $t\ne0$, which merges to $F$ as $t\to0$. Let
$\varphi^h_t\in H_2({\mathfrak X}^h_t\setminus F^h_t,\Z/2)$, respectively,
$\varphi^e_t\in H_2({\mathfrak X}^e_t\setminus F^e_t,\Z/2)$,
$t\ne0$, be families of conjugation invariant
classes converging to
the same class $\varphi\in H_2(Y\setminus F,\Z/2)$ as $t\to0$.

\begin{corollary}\label{c3}
Assume that the quadruple $(Y,
E, F, \varphi)$ has property (R). Let $D$ be a
real effective
divisor class on $Y$
such that $DE=0$.
Then,
$$W({\mathfrak
X}^e_t,D,F^e_{+,t},\varphi^e_t)=\sum_{m\in\Z}(-1)^mW({\mathfrak
X}^h_t,D-mE,F^h_t,\varphi^h_t),\quad t\ne0\ ,$$ where $F^e_{+,t}$ is
any component of $F^e_t$.
\end{corollary}

{\bf Proof}. This is an immediate consequence of
Theorems \ref{t2}(1,2) and \ref{t4}(1,2i), and of the relation
$\sum_{m\in\Z}(-1)^m\binom{2k}{m-k}=0$ which holds for any integer $k\ne0$. \proofend

\subsubsection{ABV formulas for Welschinger invariants, III}\label{sec9III}
Let
$\widetilde\pi':\widetilde{\mathfrak X}'\to(\C,0)$ be
a real unscrew
of a real nodal degeneration $\pi':{\mathfrak X}'\to(\C,0)$
such that
$\R Y\ne\emptyset$ (see Section \ref{sec-nod}).
Assume that there is a connected component $F\subset\R Y$
disjoint from $\R E$. If $\R E\ne \emptyset$, denote by $F'$ the connected component of $\R Y$ containing $\R E$; if $\R E=
\emptyset$, put $F'=\emptyset$. Choose a conjugation invariant class $\varphi\in H_2(Y\setminus(F\cup F'),\Z/2)$. Extend
$F$ to a family $F_t$ of connected components of $\R\widetilde{\mathfrak X}'_t$, $t\in (\R,0)$, and denote by
$F'_t$ the part of $\R\widetilde{\mathfrak X}'_t$, $t\in(\R,0)\setminus\{0\}$, which converges to
$F'\cup\R Z$ as $t\to0$. Extend
$\varphi$ to a family of conjugation invariant classes $\varphi_t\in H_2(\R\widetilde{\mathfrak X}'_t\setminus(F_t\cup F'_t),\Z/2)$.

\begin{theorem}\label{t5}
{\rm (1)} Let the signature of the unscrew $\widetilde\pi':\widetilde{\mathfrak X}'\to(\C,0)$ be $1$ or $-3$. Then,
the following holds.
\begin{enumerate}
\item[{\rm (1i)}] For any
real effective
divisor class $D$ on $\widetilde{\mathfrak X}'_t$, $t\ne0$, one has
$$W(\widetilde{\mathfrak X}'_t,D,F_t,\varphi_t)=W(\widetilde{\mathfrak
X}'_t,D+(DE)E,F_t,\varphi_t)\ ,$$
$$W(\widetilde{\mathfrak X}'_t,D,F_t,\varphi_t+[F'_t])=W(\widetilde{\mathfrak
X}'_t,D+(DE)E,F_t,\varphi_t+[F'_t])\ .$$
\item[{\rm (1ii)}] For any
real effective
divisor class $D$ on $\widetilde{\mathfrak X}'_t$, $t\ne0$, and any generic collection
$\bp^\sharp$ of $-DK_Y-1$ distinct points of $F$, one has
$$W(\widetilde{\mathfrak X}'_t,D,F_t,\varphi_t)=\sum_{m\ge0}
\binom{DE/2+2m}{m}U_{Y,E,\varphi}(D-2mE,\bp^\sharp)\ ,$$
$$W(\widetilde{\mathfrak X}'_t,D,F_t,\varphi_t+[F'_t])=\sum_{m\ge0}
\binom{DE/2+2m}{m}U_{Y,E,\varphi+[F']}(D-2mE,\bp^\sharp)\ .$$
\end{enumerate}

{\rm (2)} Let the signature of the unscrew $\widetilde\pi':\widetilde{\mathfrak X}'\to(\C,0)$ be $3$ or $-1$. Then,
for any
real effective
divisor class $D$ on $\widetilde{\mathfrak X}'_t$, $t\ne0$, one has
$$W(\widetilde{\mathfrak X}'_t,D,F_t,\varphi_t)=\sum_{m\ge0}(-2)^m
U_{Y,E,\varphi}(D-mE,\bp^\sharp)\ ,$$
$$W(\widetilde{\mathfrak X}'_t,D,F_t,\varphi_t+[F'_t])=\sum_{m\ge0}2^m
U_{Y,E,\varphi+[F']}(D-mE,\bp^\sharp)\ .$$
\end{theorem}

{\bf Proof.}
Statement (1i)
can be proved as
Theorem \ref{t2}(1), if we notice that both the rulings of the quadric
$Z$ are conjugation-invariant in the considered situation. Formulas of (1ii) follow from
Proposition \ref{t6} and Lemma \ref{limits}: a curve $C\in|D-2mE|_Y$ intersects $E$ at $DE/2+2m$ pairs of
complex conjugate points, to each pair we attach two lines belonging to the same ruling of
$Z$, and, finally, one has to choose $m$ pairs of conjugate intersection points to attach the lines of
a marked ruling (denoted $|C_2|_Y$ in Lemma \ref{limits}).

In the situation of assertion (2), the complex conjugation
interchanges the rulings of $Z$. Thus,
the formulas required
follow from
Proposition \ref{t6} and Lemma \ref{limits}:
to each of the $m$ pairs of complex conjugate intersection points we attach a line from $|C_1|_Y$
and a line from $|C_2|_Y$, and we have two ways to do so. For the sign relations, notice that $E\circ\varphi=0\in\Z/2$,
and that a pair of complex conjugate lines on
$Z$ has a solitary intersection point in $\R Z$.
\proofend

\begin{proposition}\label{ex-u}
External $u$-numbers do not depend on the choice of
$\bp^\sharp$.
\end{proposition}

{\bf Proof.}
From formulas of Theorem \ref{t5}(1ii), we can express
$U_{Y,E,\varphi}(D,\bp^\sharp)$ (resp. $U_{Y,E,\varphi+[F']}(D,\bp^\sharp)$) as
a linear combination of $W(\widetilde{\mathfrak X}'_t,D-2m,F_t,\varphi_t)$, $m\ge0$,
(resp. $W(\widetilde{\mathfrak X}'_t,D-2m,F_t,\varphi_t+[F'_t])$, $m\ge0$).
\proofend

\begin{corollary}\label{c4}
Let $\pi':{\mathfrak X}'\to(\C,0)$ be a real nodal degeneration,
$\pi^{\theta}:{\mathfrak X}^{\theta}\to(\C,0)$ and $\pi^{-\theta}:{\mathfrak X}^{-\theta}\to(\C,0)$ the
mirror unscrews obtained from $\pi'$ (of signature $\theta$ and $-\theta$, respectively).
Let the (common for ${\mathfrak X}^{\theta}_0$ and ${\mathfrak X}^{-\theta}_0$)
real nodal del Pezzo pair $(Y,E)$ have a connected component $F$ of $\R Y\ne
\emptyset$ disjoint from $E$. Choose a conjugation invariant class $\varphi\in H_2(Y\setminus(F\cup F'),\Z/2)$.

If $\theta = 3$ or $-1$, then, for any real effective divisor class $D\in\Pic(Y)$ such that
$DE=0$ and and $t\in(\R,0)\setminus\{0\}$, one has
$$W({\mathfrak X}^{\theta}_t,D,F_t,\varphi_t)=
\sum_{m\in\Z}(-1)^mW({\mathfrak X}^{-\theta}_t,D+mE,F_t,\varphi_t)\ ,$$
$$W({\mathfrak X}^{\theta}_t,D,F_t,\varphi_t+[F'_t])=
\sum_{m\in\Z}W({\mathfrak X}^{-\theta}_t,D+mE,F_t,\varphi_t+[F'_t])\ .$$
\end{corollary}

{\bf Proof.}
This is an immediate consequence of Theorem \ref{t5}.
\proofend


\section{Proof
of positivity and asymptotics statements}\label{sec-pa}

\subsection{ABV families} \label{sec-deg}
Any two real del Pezzo surfaces
of degree $2$ that
have homeomorphic real parts
are deformation equivalent in the class of such surfaces
(see, for example, \cite{DIK}, Theorem 17.3), and, as a result,
they
have the same system of Welschinger invariants. Therefore, in the proof of
Theorems \ref{t9} and \ref{t11}, for each topological type,
it is sufficient to
pick a
particular
real del Pezzo surface, $X$,
that we include
into an appropriate family with a special fiber containing
a real nodal del Pezzo pair
$(Y,E)$;
the family depends on the choice of a connected component $F$ of $\R X$.

If $X$ is of type $\PP^2_{a,b}$, $a+2b=7$, (or, equivalently, of
type $\langle q\rangle^-$, $1\le q \le4$), we include $X$ into a family
${\mathfrak X}\to(\C,0)$, which is a holomorphic submersion possessing
a real structure
subject to (\ref{ers})
and whose central fiber is a real nodal del Pezzo pair
$(Y,E)$ ({\it cf.}
Section \ref{sec9I}). Namely, we specialize a conjugation-invariant set of
$6$ blown up points on a real conic $C_2$, $E$ being the strict
transform of $C_2$.
We call ${\mathfrak X}\to(\C,0)$ a {\it regular ABV family}
of
$X$.

\begin{figure}
\begin{center}
\epsfxsize 145mm \epsfbox{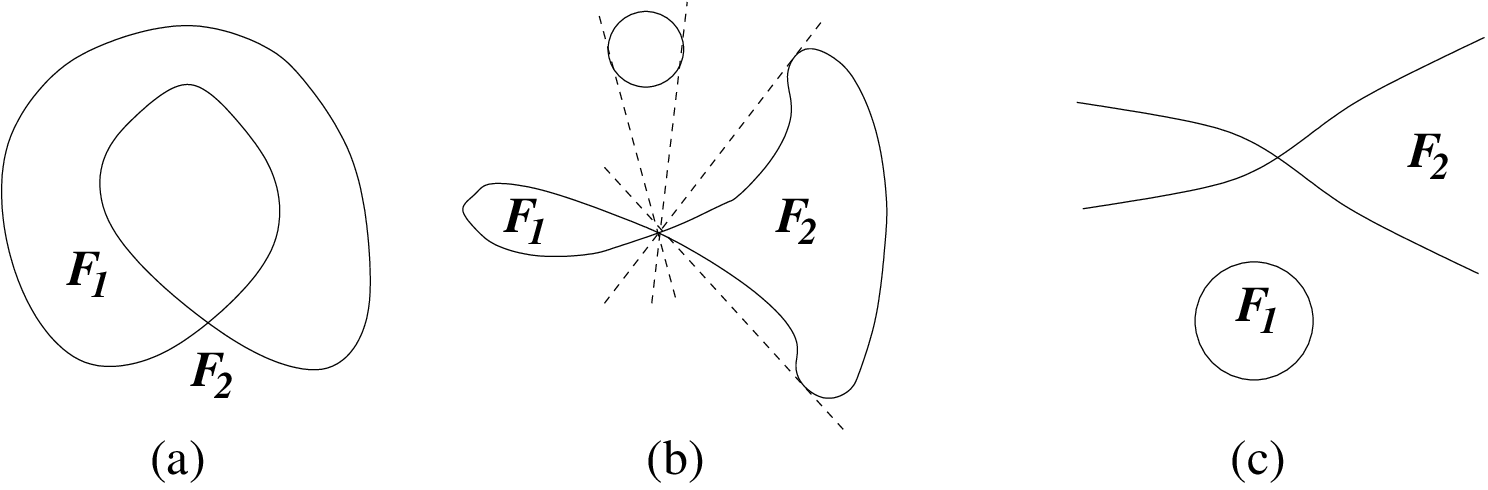}
\end{center}
\caption{Nodal quartics}\label{f2}
\end{figure}

The real del Pezzo surfaces $X$ of degree $2$
of other
deformation
types
can be included into the following unscrews of real nodal degenerations corresponding to
nodal degenerations of quartics $Q_X$:
\begin{enumerate}\item[(a)]
{\it hyperbolic ABV families}: \begin{itemize}\item if $X$ is of type $1\langle1\rangle^+$,
we degenerate $Q_X$ into a
nodal quartic shown in Figure \ref{f2}(a) and choose a
$1$-unscrew
assuming that the component $F\subset\R X$ doubly covers
the
annulus merging to
the domain
$F_1$; \item
if $X$ is of type $\langle 2\rangle^+$,
we degenerate $Q_X$ into
a real nodal quartic as shown in Figure \ref{f2}(b) and choose a $1$-unscrew,
assuming that $F$ doubly covers
the disc merging to the domain $F_1\cup F_2$;
\end{itemize} \item[(b)] {\it elliptic ABV families}:
\begin{itemize}\item if $X$ is of type
$1\langle1\rangle^-$, $F=F^{\;no}$,
we degenerate $Q_X$ into a
nodal quartic shown in Figure \ref{f2}(a) and choose a $(-1)$-unscrew,
assuming that $F$ doubly covers
the M\"obius band merging to the domain
$F_2$; \item if $X$ is of type
$\langle1\rangle^+$,
we degenerate $Q_X$ into a
nodal quartic shown in Figure \ref{f2}(a) and choose a $(-1)$-unscrew,
assuming that $F$ doubly covers
the disc merging to the domain
$F_1$; \item
if $X$ is of type $\langle0\rangle^-$, we degenerate
the quartic $Q_X$ (having an empty real part) into a nodal quartic with a one-point real part,
and choose a $(-1)$-unscrew;\end{itemize} \item[(c)] {\it
external ABV families}: \begin{itemize} \item
if $X$ is of type
$1\langle1\rangle^-$, $F=F^{\;o}$,
we degenerate $Q_X$ into a
nodal quartic shown in Figure \ref{f2}(c) and choose a $(2,1)$-unscrew,
assuming that $F$ doubly covers
the disc merging to
$F_1$, and $F'$ doubly covers a M\"obius band merging to the domain $F_2$; \item
if $X$ is of type $\langle q\rangle^+$, $q=3,4$, we degenerate $Q_X$ into
a nodal quartic so that one of the ovals collapses to a point and
choose a $(3,0)$-unscrew, assuming that $F'$ doubly covers
the disc merging to a point.
\end{itemize}\end{enumerate}

\subsection{$F$-compatible divisor classes}\label{F-compatible}
The
$F$-compatibility condition (see Section \ref{sec414})
holds for all real divisor classes if $H_1(\R X\setminus F,\Z/2)=0$.
Hence, the only cases with a non-trivial condition
for del Pezzo surfaces of degree $2$
are as follows:
either $X$ is of type $\langle0\rangle^-$,
or $X$ is of type $1\langle1\rangle^-$ and $F=F^{\;o}$. For example, $-K_X$ is not $F$-compatible in either
of these two cases.

We
use below the following
characterization of the $F$-compatibility condition
for the two
cases
mentioned above.

The elliptic ABV family for a surface of type $\langle0\rangle^-$ (see Section \ref{F-compatible})
has the central fiber $Y
\cup Z$, where $Z$ is a quadric and $Y$
is the plane blown up at three
pairs of complex conjugate points on a real conic $C_2$ such that
$\R C_2\ne\emptyset$, and at one more real point belonging to the orientable component of
$\R\PP^2\setminus\R C_2$. Consider the basis of $\Pic(Y)$ consisting of the
pull-back $L$ of a generic line
and the exceptional divisors
$E_1,...,E_7$, where $E_{2i},E_{2i+1}$ are complex conjugate,
$i=1,2,3$, and $E_1$ is real. Since the components
of $\R X$ for $X$ of type $\langle0\rangle^-$ are interchanged by an automorphism,
we can choose any of them, and we assume that $F$ is disjoint from $\R E_1$.
Following
Proposition \ref{new-proposition},
let us identify $\Pic^\R (X)$  with a subgroup of  $\Pic^\R(Y)$.

\begin{proposition}\label{remark-prop}
For a surface $X$ of type $\langle0\rangle^-$, a divisor class
$D\in\Pic(X)$, represented as $D=dL-d_1E_1-...-d_7E_7$ in $\Pic(Y)$,
\begin{itemize}
\item is real if and only if
$d_{2i}=d_{2i+1}$, $i=1,2,3$, and $2d=d_2+...+d_7$,
\item is $F$-compatible if and only if the number
$d_1=DE_1$ is even.
\end{itemize}
\end{proposition}

{\bf Proof}. Straightforward.
\proofend

We say that a divisor
class $D\in\Pic^\R(Y)$ is {\it $F_+$-compatible} if it satisfies conditions
of
Proposition \ref{remark-prop}.

\smallskip

Let $X$ be a surface of type $1\langle1\rangle^-$. Consider a degeneration of $Q_X$ into $C_1\cup C_3$, where
$C_3$ is a real two-component cubic, and $C_1$ a line crossing the one-sided component of $C_3$ in three real
points. This degeneration induces a conjugation invariant family of surfaces,
in which the component $F=F^{\;o}\subset\R X$ merges to the sphere doubly covering the disc bounded by the oval of
$C_3$. Making
the base change $\tau=t^2$ and blowing up the three nodes of the central fiber, we can realize a triple
$(2,1)$-unscrew with the central fiber $Y\cup Z^{(1)}\cup Z^{(2)}\cup Z^{(3)}$, where $Y$ is a smooth real surface
with $\R Y$ diffeomorphic to $\R X$, $Z^{(i)}$,
$i=1,2,3$, are disjoint real quadrics with $\R Z^{(i)}\simeq(S^1)^2$, intersecting $Y$ along real
$(-2)$-curves $E^{(i)}$, $i=1,2,3$, respectively. Contracting each quadric along one of its rulings, we obtain a
real family $\pi:{\mathfrak X}\to(\C,0)$, where ${\mathfrak X}$ is smooth, $\pi$ is a submersion, and ${\mathfrak X}_0=Y$.

\begin{proposition}\label{p2}
The subgroup $\Pic^\R(X)\simeq\Pic^\R(Y)\subset\Pic(Y)$ is generated by $K_X=K_Y$ and $E^{(i)}$, $i=1,2,3$.
A divisor class $D=-dK_X+d_1E^{(1)}+d_2E^{(2)}+d_3E^{(3)}$ is $F$-compatible if and only if $d+d_1+d_2+d_3$ is even.
\end{proposition}

{\bf Proof.} Straightforward.\proofend

\subsection{Auxiliary statements}\label{sec-as}
Let
$\pi:{\mathfrak X}\to
(\C,0)$ be a proper
holomorphic submersion of a smooth three-dimensional variety ${\mathfrak X}$  (with $(\C,0)$ being understood as a disc germ),
where
each fiber ${\mathfrak X}_t, t \ne 0$, is a del Pezzo
surface of degree $2$
and the central fiber $Y
={\mathfrak X}_0$ contains a smooth rational curve $E$ such that
$(Y,E)$ is a monic log-del Pezzo pair. In what follows
we
identify the Picard groups of the fibers
as in Remark \ref{r3}.

\begin{lemma}\label{ll1}
Let $X={\mathfrak X}_{t}$ for some
$t\ne 0$, and let $D\in \Pic(X)$.
\begin{enumerate}
\item[{\rm (i)}] If $D$ is big and $X$-nef, then $-K_XD>1$, and the linear system $|D|_X$ contains an irreducible rational curve.
\item[{\rm (ii)}] The divisor class $D$ is $X$-nef if and only if its
intersection with any $(-1)$-curve on $X$ is non-negative. In this
case $D^2\ge0$.
The divisor class
$D$ is $Y$-nef if and only if its intersection with $E$ and any
$(-1)$-curve on $Y$ is non-negative. If $D$ is $Y$-nef then it is
$X$-nef.
\item[{\rm (iii)}]
If $D$ is nonzero and $X$-nef, and satisfies $D^2=0$, then $D=kD'$,
where $D'$ is primitive {\rm (}i.e., not multiple of another
divisor{\rm )}. Furthermore, $-K_XD'=2$, $\dim|D'|_X=1$, and a
generic element of $|D'|_X$ is a smooth connected rational curve. If
$D'E\ge0$, then $|D'|_Y$ is one-dimensional with a smooth,
connected, rational curve as a generic element. Furthermore, if
$D'E>1$, then $D'=-(K_Y+E)$.
\item[{\rm (iv)}] If
$D\in\Pic(Y,E)$ a $Y$-nef and big divisor, satisfying $R_Y(D,0)>0$.
Then the divisor class $D'=D-E-\sum_{E'\in\;{\mathcal E}(E)^{\perp D}}E'$ is
$Y$-nef and satisfies $D'E'=0$ for all $E'\in{\mathcal E}(E)^{\perp D}$;
furthermore, if $D'\ne0$, then $D'$ is presented by the union of
curves different from $E$ and crossing $E$ positively.
\item[{\rm (v)}] If $D$ is big and $X$-nef such that
$
{\mathcal E}(E)^{\perp D}
\ne \emptyset$,
then $D - mE$ with $m > 0$
cannot be represented by an irreducible curve in $Y$.
\end{enumerate}
\end{lemma}

{\bf Proof}. It is known that big and nef divisors on del Pezzo
surfaces are effective and can be represented by irreducible
rational curves (see, for instance, \cite[Theorems 3, 4, and Remark
3.1.4]{GLS}). Hence $-K_XD>0$. In the case when $D=-K_X$ or when
$-K_X-D$ is effective, the inequality $-K_XD>1$ can easily be
verified. If $-K_X-D$ is not effective, then $-K_XD>1$ due to
$\dim|-K_X|=2$.

Statement (ii) on the $X$-nefness (respectively, $Y$-nefness)
follows from the fact that the effective cone in $\Pic(X)$
(respectively, $\Pic(Y)$) is generated by $(-1)$-curves
(respectively, by $(-1)$-curves and $E$).

If $D$ is $Y$-nef, then $DE'\ge0$ for all $(-1)$-curves $E'$ on $Y$,
and $DE\ge0$. Any $(-1)$-curve $E''$ in $X$ degenerates either into
a $(-1)$-curve of $Y$, or into a curve $E+E'$ with a $(-1)$-curve
$E'$ on $Y$, and hence in both the cases $DE''\ge0$.

The nonnegativity of $D^2$ in statement (ii) and the part of
statement (iii), concerning divisors and linear systems on $X$,
follow, for instance, from \cite[Theorems 3, 4, and Remark
3.1.4]{GLS}. In particular, if $D'$ is primitive and satisfies
$(D')^2=0$, then a general curve in $|D'|_X$ is non-singular,
rational.

Let $D'E\ge0$ in statement (iii). Suppose that $D'E'=0$ for some
$(-1)$-curve $E'\in{\mathcal E}(E)$. Then we can blow down $E'$ and
reduce the degenerating family to a family of del Pezzo surfaces,
which immediately yields that $\dim|D'|_Y=1$ as well as the fact
that a generic element of $|D'|_Y$ can be chosen to be a
smooth rational curve. Suppose that $D'E'>0$
for all $E'\in{\mathcal E}(E)$.
By
Proposition \ref{t6}, a general
curve $C_t\in|D'|_{{\mathfrak X}_t}$, $t\ne0$, degenerates into a
curve $C_0+mE\in|D'|_Y$ with some $m\ge0$ and $C_0\not\supset E$. If
$m$ were positive, we would have $(C_0)^2=(D')^2-2D'E-2m^2\le-2$,
and, in view of $-K_YD'=2$ (comes from genus formula) and
$-K_YC\ge-1$ for all irreducible curves $C\ne E$, we would get $C_0$
consisting of components with negative self-intersection, a
contradiction to $\dim|D'|_Y\ge1$.

Let $D'E>1$ in statement (iii), then
$-(K_Y+E)D'=-K_YD'-D'E=2-D'E\le0$, which in view of the nefness of
$-(K_Y+E)$, yields $D'=-(K_Y+E)$.

In view statement (iii), to prove (iv) it is enough to check that
$D-E$ non-negatively crosses each $(-1)$-curve of $Y$. If ${\mathcal
E}(E)^{\perp D}=\emptyset$, then this immediately follows from the fact that
$EE'=1$ for all $E'\in{\mathcal E}(E)$.

In the case of ${\mathcal E}(E)^{\perp D}\ne\emptyset$, we have $D'E'<0$ for
all $E'\in{\mathcal E}(E)^{\perp D}$, and hence $D-mE$ with $m>0$ cannot be
represented by an irreducible curve in $Y$: any curve in $|D-mE|_Y$
must contain all $E'\in{\mathcal E}(E)^{\perp D}$ as components, and such a
component cannot be unique, otherwise $D$ would not be big.
This proves (v). Furthermore, formula (\ref{ee1}) and statement (i)
yield
\begin{equation}
N_Y(D,0,(DE)e_1)=GW_0(X,D)>0\ .\label{e2056}\end{equation} Since
$R_Y(D,0)>0$, computing $N_Y(D,0,(DE)e_1)$ via a sequence of
formulas (66) from \cite{MS2} written in the form
\begin{equation}N_Y(D,0,je_1,(DE-j)e_1)=N_Y(D,0,(j+1)e_1,(DE-j-1)e_1)+S^\C_j\
,\label{e2058}\end{equation}
$$j=0,...,DE\ ,$$ where $S^\C_j$ stands for the second sum in the
right-hand side of the cited formula, and $N_Y(D,0,(DE+1)e_1,-e_1)$
is zero by definition, we get $S^\C_0+...+S^\C_{D'E}=GW_0(X,D)>0$.
That means the divisor $D-E$ is effective, and the divisor class
$D'=D-E-\sum_{E'\in\;{\mathcal E}(E)^{\perp D}}E'$ is represented by a curve
$C'$, whose all components are disjoint from the $(-1)$-curves
$E'\in{\mathcal E}(E)^{\perp D}$ and intersect with $E$. Notice, first, that
$C'$ does not contain $(-1)$-curves disjoint from $E$, and hence
$D'E''\ge0$ for all $(-1)$-curves $E''$ with $E''E=0$, and, second,
$(D')^2\ge0$, since otherwise, $C'$ would contain a $(-1)$-curve
crossing $E$ and disjoint from the other components of $C'$ and from
$E'\in{\mathcal E}(E)^{\perp D}$, contrary to the definition of ${\mathcal
E}(E)^{\perp D}$. Altogether this yields the required statement. \proofend

\begin{remark}\label{lemma22}
Notice that Lemma \ref{ll1}
can be applied
to all ABV families introduced in Section \ref{sec-deg}.
Namely, over $\C$ we can
contract the quadric surface in the central fiber of the family along one of the rulings and thus obtain
a family exactly as in Lemma \ref{ll1}.
\end{remark}

The following two claims will be used in the proof of the asymptotic statements
in Theorems \ref{t9} and \ref{t11}.

\begin{lemma}\label{ll1a}
Let $\{a_n\}_{n\ge0}$ be a sequence of positive numbers, $a_0=1$,
and let $0\le f(n)\le n$ an integral-valued function. If
\begin{itemize} \item either
\begin{equation}a_{n+1}\ge\lam a_{f(n)}a_{n-f(n)},\quad \text{for all}\ n\ge n_0\ \text{and some}\ \lam>0\
,\label{e2062}\end{equation}
\item or
\begin{equation}a_n\ge\lam a_{f(n)}a_{n-f(n)},\quad \text{for all}\ n\ge n_0\ \text{and some}\
\lam>0\ , \label{e2063}\end{equation}\end{itemize} then there exist
$\xi_1,\xi_2>0$ such that $a_n\ge\xi_1\xi_2^n$ for all $n\ge n_0$.
\end{lemma}

{\bf Proof}. Straightforward induction on $n$ with $\xi,\eta$ found
from the equations $$\lam\xi_1=\xi_2,\quad \xi_1\xi_2^{n_0}=a_0$$ in the
first case, and the equations $$\lam\xi_1=1,\quad
\xi_1\xi_2^{n_0}=a_{n_0}$$ in the second case.
\proofend

\begin{lemma}\label{las}
Asymptotic relations (\ref{ee21}) and (\ref{e2078}) follow from
$$
\log W(X,nD,\R X,[\R X\setminus F]) \ge-K_XD \cdot n \log n + O(n), \;\;\; n \to + \infty\ .$$
\end{lemma}

{\bf Proof}.
Straightforward from
$$
\displaylines{
\log \vert W(X,nD,\R X,[\R X\setminus F])\vert \le\log GW_0(X,nD) \cr
=-K_XD \cdot n \log n + O(n), \;\;\; n \to + \infty\ ,
}
$$
(see \cite[Theorem 1]{IKS5}).
\proofend

\subsection{Non-negativity of $w$-numbers}\label{sec-nn} In each
of the regular, elliptic, or hyperbolic ABV families
introduced in Section \ref{sec-deg}, the central fiber
coincides with or contains a real
del Pezzo surface $Y$ of
degree $2$ with a smooth real rational curve $E\subset Y$ whose
real part $\R E$ lies in some connected component $F$ of
$\R Y$ (recall that $\R E \ne
\emptyset $).

\begin{lemma}\label{non-negativity}
If $Y$ appears in
a regular or hyperbolic ABV family
introduced in Section
\ref{sec-deg}, then, for any
$D\in\PicPPR(Y,E)$ such that $-K_YD\ge1$ and, for any
$\alp,\bet\in\Z^{\infty,\odd}_+$ such that $I(\alp+\bet)=DE$, we
have
\begin{equation}W_{Y,E,[\R Y\setminus F]}(D,\alp,\bet,0)\ge0\
.\label{e2040}\end{equation} If $Y$ appears in
an elliptic ABV family
introduced in Section \ref{sec-deg}, then, for any component $F_+$ of
$F\setminus\R E$, for any
$D\in\PicPPR(Y,E)$ such that $-K_YD\ge2$ and $DE$ is even, and
for any
$\alp,\bet\in\Z^{\infty,\even}_+$ such that $I(\alp+\bet)=DE$, we
have \begin{equation}W^-_{Y,F_+,[\R Y\setminus
F]}(D,\alp,\bet,0)\ge0\ .\label{e2041}\end{equation}
\end{lemma}

{\bf Proof}.
Suppose, first, that $Y$ appears in a
regular, hyperbolic, or elliptic ABV family
of a del Pezzo surface $X$ whose type is different from
 $\langle2\rangle^+$.
We use induction on
$R_Y(D,\bet)$ to prove
(\ref{e2040}) via formula (\ref{e44}) and to prove
(\ref{e2041}) via formula (\ref{e44om}). Notice that the
coefficients in both these formulas are non-negative. Indeed,
the values of
$\eta(l)$ given in (\ref{ee6}) and related to the considered cases
are all non-negative, since, if $L',L''$ are real, then always $\R L',\R L''\subset\overline F_+$.
Thus, it remains to verify the non-negativity of the initial values given in Propositions \ref{ini1} and \ref{ini2}. In the case of $X$
of type $\PP^2_{a,b}$, $a+2b=7$,
this is so, since ${\mathcal E}(E)$ consists of $2b$ pairs of
disjoint complex conjugate lines and of $6-2b$ pairs of intersecting
real lines. In the case of $X$ of type $1\langle1\rangle^+$,
$\langle0\rangle^-$, $1\langle1\rangle^-$, or $\langle1\rangle^+$ this is so, since the
corresponding real nodal plane quartic curve $Q_Y$ has no real lines
passing through the node $z$ of $Q_Y$ and tangent to $Q_Y$ at a
point $z'\ne z$ ({\it cf.} Remark \ref{r4}).

Suppose now that $Y$ appears in
a hyperbolic ABV family
of a del Pezzo surface $X$ of type
$\langle2\rangle^+$. We again use induction on
$R_Y(D,\bet)$ and prove
(\ref{e2040}) via formula (\ref{e44}). The base of induction is provided by Proposition \ref{ini1}(1), where all values equal $1$.
However, some terms in the second sum of the left-hand side of (\ref{e44})
can be negative, and to proceed by induction, we will modify formula
(\ref{e44}) in order to cancel out the negative summands.
Consider the nodal quartic $Q_Y$ (see Figure \ref{f2}(b)). It
has $6$ tangents passing through the node $z$: four real (shown by dashed lines), $L_1,L_2$
(tangent to the domain $F_2$), $L_3,L_4$ (tangent to the oval), and two complex conjugate $L_5,L_6$.
Each tangent line $L_i$ lifts to a pair of $(-1)$-curves $E'_i,E''_i\subset Y$ intersecting at one point, and
we have
\begin{equation}E''_i=\conj(E'_i),\ i=1,2,3,4,\quad E'_6=\conj(E'_5),\ E''_6=\conj(E''_5)\ ,
\label{e-new3}\end{equation}
and $E'_i+E''_i\in|-K_Y-E|$, $i=1,...,6$. By Proposition \ref{ini1}(2i,ii),
\begin{equation}W_{Y,E,[\R Y\setminus F]}(\{E'_i,E''_i\},0,0,e_1)=\begin{cases}-1,\quad & i=1,2,\\
1,\quad & i=3,4,\end{cases}\label{e-new1}\end{equation} \begin{equation}W_{Y,E,[\R Y\setminus F]}(\{E'_5,E'_6\},0,0,e_1)=
W_{Y,E,[\R Y\setminus F]}(\{E''_5,E''_6\},0,0,e_1)=1\label{e-new2}\end{equation}
(here we denote by $F$ the connected component of $\R Y$, which deforms into
the considered component $\R X$; no confusion will arise). Now notice that
the coefficients in the second sum of LHS of (\ref{e44}) are positive, since $E_{1/2}\circ[\R Y\setminus F]=L_{1/2}\circ[\R Y\setminus F]=0$. Each summand of that sum can be
written either as $(l+1)AB_m$, or $(l+1)A'B'_m$, or $(l+1)A''B''_m$, where all the factors of type (\ref{e-new1}) and (\ref{e-new2}) are separated to $A$, and the sum of the divisor classes in the factors, participating in $B_m$, $B'_m$, or $B''_m$ equals either $D-E+m(K_Y+E)$, $D-E+m(K_Y+E)-(E'_5+E'_6)$,
or $D-E+m(K_Y+E)-(E''_5+E''_6)$, respectively. Observe that, by Theorem
\ref{t1}(1g) the factors (\ref{e-new1}) and (\ref{e-new2}) appear in $A$ at most once, and that
$A'$ (resp. $A''$) necessarily contains the factor $W_{Y,E,[\R Y\setminus F]}(\{E'_5,E'_6\},0,0,e_1)$
(resp. $W_{Y,E,[\R Y\setminus F]}(\{E''_5,E''_6\},0,0,e_1)$), and, for a given $m\ge0$, all combinations of $l$ and $A$ (resp.,
$A'$ or $A''$) subject to the above restrictions are allowed. Thus, an easy computation gives that, combining together all summands with the same $B_m$, $B'_m$ or $B''_m$, $m\ge0$, we finally reduce formula (\ref{e44}) to the form
$$W_{Y,E,[\R Y\setminus F]}(D,\alp,\bet,0)=\sum_{j\ge 1,\
\bet_j>0}W_{Y,E,[\R Y\setminus F]}(D,\alp+e_j,\bet-e_j,0)$$ \begin{equation}\qquad\qquad+B_0+B_2+B'_0+B''_0\ ,\label{e44-2}
 \end{equation} which completes the proof in view of $B_0,B_2,B'_0,B''_0\ge0$ (by the induction assumption).

\proofend

\subsection{Proof of Theorem \ref{t9}}
It
is clear that the case of any real del Pezzo surface of degree
$\ge3$ can be reduced
to
that of degree $2$
by blowing up suitable real points.
To treat the degree $2$ case,
we include
$X$
into
a suitable ABV family, defined in Section \ref{sec-deg}
and apply
Theorems \ref{t2} and \ref{t4} expressing Welschinger invariants in
terms of $w$-numbers.
The nodal del Pezzo pairs in the central fibers of the considered ABV families
are monic log-del Pezzo, which
allows us to compute and estimate the above $w$-numbers using
Theorems \ref{t1} and \ref{t3}. The key observation is the non-negativity of the
considered $w$-numbers (see Lemma \ref{non-negativity}).

\subsubsection{Positivity and asymptotics statements for surfaces of
types $\PP^2_{a,b}$, $a+2b=7$,
$1\langle1\rangle^+$, and $\langle q\rangle^+$, $q=2,3,4$}\label{sec-pos1}

Consider
first surfaces $X$ of types $\PP^2_{a,b}$, $a+2b=7$,
$1\langle1\rangle^+$, and $\langle 2\rangle^+$ and
the
regular or hyperbolic ABV families for them introduced in Section \ref{sec-deg}.
Following
Proposition \ref{new-proposition},
we identify $\Pic^\R(X)$ and $\Pic^\R(Y)$,
where $Y$
is
the central fiber of the corresponding ABV family.
Furthermore, we
restrict
ourselves
to the case $DE\ge0$ ({\it cf.} Theorem \ref{t2}(1)). Fix a connected component $F\subset\R X$
and put $\varphi=[\R X\setminus F]\in H_2(Y\setminus F,\Z/2)$.

\noindent {\bf Positivity
for surfaces of
types $\PP^2_{a,b}$, $a+2b=7$,
$1\langle1\rangle^+$, and $\langle 2\rangle^+$}.
From formula (\ref{ee2}), Proposition \ref{u=w}, and inequality (\ref{e2040}) it follows
that
\begin{equation}
W(X,D,F,\varphi)\ge W(X,D-mE,F,\varphi)\quad\text{for all}\ m\ge0\
.\label{e2057}\end{equation} Indeed,
the both terms are sums of non-negative $w$-numbers,
and all the $w$-numbers occurring
in the development of the right-hand side
appear in the development of
the left-hand side
with non-smaller coefficients:
$$\binom{DE+2m+2k}{m+k}\ge\binom{DE+2m+2k}{k}\quad \text{for all}\ k\ge0\
.$$

We will prove inequality (\ref{ee20}) for all real big and $X$-nef
divisors $D\in\Pic(X)$ such that $DE \geq 0$, using induction on
$\rho(D)=-(K_X+E)D$.

Observe that $\rho(D)>0$, since
$D$ is big and $|-K_X-E|$ defines a conic bundle. Suppose that
$\rho(D)=1$, or, equivalently, $R_Y(D,0)=0$. From formula
(\ref{ee2}), Proposition \ref{u=w}, and inequality (\ref{e2040}) we get
\begin{equation}
W(X, D,F,\varphi)\ge W_{Y,E,\varphi}(D,0,(DE)e_1,0)\ ,\label{e2061}
\end{equation}
and then, applying formula (\ref{e44}) (resp. (\ref{e44-2}) for $X$ of type
$\langle2\rangle^+$) $DE$ times, we end up with
$$W_{Y,E,\varphi}(D,0,(DE)e_1,0)\ge W_{Y,E,\varphi}(D,(DE)e_1,0,0)=1\ ,$$
the latter equality coming from Proposition \ref{ini1}(1iii) and
inequality
$$DE=-K_YD-\rho(D)=-K_XD-1\overset{\text{Lemma \ref{ll1}(i)}}{>}0\ .$$

Assume that $\rho(D)>1$, or, equivalently,
\begin{equation}R_Y(D,0)>0\ .\label{e2048}\end{equation}
Suppose that ${\mathcal E}(E)^{\perp D}=\emptyset$.
Then, $D-E$ is $Y$-nef and satisfies
$(D-E)^2\ge0$ (see Lemma \ref{ll1}(ii)).

If $(D-E)^2>0$, then $D-E$ is big and $X$-nef, and
$\rho(D-E)=\rho(D)-2<\rho(D)$. Thus, by the induction assumption and
(\ref{e2057})
$$W(X,D,F,\varphi)\ge W(X,D-E,F,\varphi)>0\ .$$

If $(D-E)^2=0$, then by Lemma \ref{ll1}(iii), $D-E=kD''$ with
$k\ge1$ and a primitive $D''\in\PicPPR(Y,E)$ such that $D''E>0$,
$(D'')^2=0$, $\dim|D''|_Y=-K_YD''-1=1$, and the linear system
$|D''|_Y$ contains a real, rational, smooth curve. If $k=1$, then
$W(X,D-E,F,\varphi)=1$, and again (\ref{ee20}) follows. If $k\ge2$, then
$R_Y(D,(DE)e_1)=-K_YD-1=-kK_YD''-1=2k-1$, and we get
\begin{eqnarray}
W(X,D,F,\varphi)&\overset{\text{(\ref{ee2})\ \&\ (\ref{e2040})}}{\ge}&W_{Y,E,\varphi}(D,0,(DE)e_1,0)\nonumber\\
&\overset{\text{(\ref{e44}),(\ref{e44-2})\ \&\ (\ref{e2040})}}{\ge}&W_{Y,E,\varphi}(D,(k-2)e_1,k((D''E)-1)e_1,0)\nonumber\\
&\overset{\text{(\ref{e44}),(\ref{e44-2})\ \&\
(\ref{e2040})}}{\ge}&\left(W_{Y,E,\varphi}(D'',0,(D''E)e_1,0)\right)^k\overset{\text{Lemma
\ref{ll1}(iii)}}{=}1\ .\label{e2044}\end{eqnarray}

Suppose now that ${\mathcal E}(E)^{\perp D}\ne\emptyset$. By
Lemma~\ref{ll1}(v) and Theorem~\ref{t2}(2),
we have
\begin{equation}
W(X,D,F,\varphi)=W_{Y,E,\varphi}(D,0,(DE)e_1,0)\ .\label{e2045}
\end{equation}
By Lemma~\ref{ll1}(iv), the divisor class
$D'=D-E-\sum_{E'\in\;{\mathcal E}(E)^{\perp D}}E'$ is $Y$-nef, and hence, by
Lemma \ref{ll1}(ii), is also $X$-nef.

Assume that $(D')^2>0$. Since $\rho(D')=\rho(D)-2$, we have
$W(X,D',F,\varphi)>0$, which due to ${\mathcal E}(E)^{\perp D'}\supset{\mathcal
E}(E)^{\perp D}\ne\emptyset$ yields ({\it cf.} (\ref{e2045}))
$$W_{Y,E,\varphi}(D',0,(D'E)e_1,0)=W(X,D',F,\varphi)>0\ .$$
Appropriately applying formula (\ref{e44}) (resp. (\ref{e44-2}) for $X$ of type $\langle2\rangle^+$)
and using (\ref{e2040}),
we obtain
\begin{eqnarray}W(X,D,F,\varphi)&=&W_{Y,E,\varphi}(D,0,(DE)e_1,0)\nonumber\\
&\ge& W_{Y,E,\varphi}(D,(s-1)e_1,(DE-s+1)e_1,0)
\ ,\label{e2049}\end{eqnarray}
where $s=
\card({\mathcal E}(E)^{\perp D})$.
In addition, \begin{eqnarray}
W_{Y,E,\varphi}(D,(s-1)e_1,(DE-s+1)e_1,0) \ge W_{Y,E,\varphi}(D',0,(D'E)e_1,0)>0
\ .\label{e2049-new}
\end{eqnarray}
Indeed,
\begin{itemize}\item $D'E>0$ ({\it cf.} Lemma \ref{ll1}(1iv)), $s=-K_YD+K_YD'\le-K_YD-1$,
$s-1=(D-E-D')E-1=DE+1-D'E\le DE$, and $DE-s+1=D'E-1$;
\item if $X$ is of type $\langle2\rangle^+$, then the non-empty set ${\mathcal E}(E)^{\perp D}$
must be either $\{E'_5,E'_6\}$, or $\{E''_5,E''_6\}$ ({\it cf.} (\ref{e-new3})), and hence
$W_{Y,E,\varphi}(D',0,(D'E)e_1,0)$ appears in the summand $B'_0$, resp. $B''_0$ in
(\ref{e44-2}).
\end{itemize}

Assume that $(D')^2=0$. If $D'=0$, then the relations (\ref{e2049}) and (\ref{e2049-new})
transform to
$$W(X, D, F,\varphi) \geq \prod_{\D_1} W_{Y,E,\varphi}(\D_1,0,e_1,0)
\prod_{\D_2} W_{Y,E,\varphi}(\D_2,0,0,e_1)
> 0,$$
where $\D_1$ runs over the real elements of
${\mathcal E}(E)^{\perp D}$, and $\D_2$ runs over the pairs of complex
conjugate elements in ${\mathcal E}(E)^{\perp D}$. If $D'\ne0$, by Lemma
\ref{ll1}(iii), $D'=kD''$ with one-dimensional linear system
$|D''|_Y$ represented by a smooth real rational curve $C''$, and the
relations (\ref{e2049}) and (\ref{e2049}) transform to
$$W(X,D,F,\varphi)\ge W_{Y,E,\varphi}(D,le_1,(DE-l)e_1,0)$$
\begin{equation}\ge (W_{Y,E,\varphi}(D'',0,(D''E)e_1,0))^k=1\
,\label{e2050}\end{equation} where $l=-K_YD-2-k$. Note that
$$k=(-K_YD'')^{-1}(-K_YD-s)\le(-K_YD-1)/2\le-K_YD-2$$
and $l\le DE$, the latter inequality coming from the relations
$$\displaylines{
DE = D'E - 2 + E \cdot \sum_{E' \in\; {\mathcal E}(E)^{\perp D}}E' = D'E - 2
+ s \geq k - 2 + s, \cr l = -K_YD - 2 - k = -K_YD' - K_YE - K_Y
\cdot \sum_{E' \in\; {\mathcal E}(E)^{\perp D}}E' - 2 - k = k - 2 + s. }
$$

\noindent {\bf Asymptotics
for surfaces of
types $\PP^2_{a,b}$, $a+2b=7$,
$1\langle1\rangle^+$, and $\langle 2\rangle^+$}.
Let $D\in\Pic(X)$ be a real, big and $X$-nef divisor. By Theorem
\ref{t2}(1) we can suppose that $D$ is $Y$-nef. We prove the
asymptotic relation (\ref{ee21}) by induction on $\tau(D) =
\min\{DE'\ :\ E'\in{\mathcal E}(E)\}$.

Let $\tau(D)=0$, or, equivalently ${\mathcal E}(E)^{\perp D}\ne\emptyset$.
Put $s=\card({\mathcal E}(E)^{\perp D})$. By Lemma \ref{ll1}(iv), there
exists
an
integer $m_0\ge 1$ such that the divisors
$D'_m=mD-E-\sum_{E'\in\;{\mathcal E}(E)^{\perp D}}E'$ are big and $Y$-nef and
satisfy $D'_mE>0$ for all integers $m\ge m_0$. Note also that by
(\ref{ee20}) and (\ref{e2045}), one has
$$W_{Y,E,\varphi}(D'_m,0,(D'_mE)e_1,0)>0,\quad m\ge m_0\ .$$ Put
$\widetilde D=D'_{m_0}$ and $\widetilde s= \card({\mathcal E}(E)^{\perp\widetilde D})$.
Again there exists an integer $m_1\ge1$ such that $\widetilde
D'_{m_1}$ is big and $Y$-nef and satisfy
$$\widetilde
D'_mE>0,\quad W_{Y,E,\varphi}(\widetilde D'_m,0,(\widetilde
D'_mE)e_1,0)>0\quad\text{for all integers}\ m\ge m_1\ .$$ For any
integer $n \ge 2$, we have decompositions
$$\begin{cases}&\widetilde D'_{nm_1}-E=
\widetilde D'_{m_1i(n)}+\widetilde
D'_{(n-i(n))m_1}+\sum_{E'\in{\mathcal
E}(E)^{\perp D}}E',\\
&-K_Y\widetilde D'_{nm_1}-2=(-K_Y\widetilde
D'_{m_1i(n)}-1)+(-K_Y\widetilde
D'_{(n-i(n))m_1}-1)+s,\end{cases}\quad
i(n)=\left[\frac{n}{2}\right],$$ and inequality $\widetilde
D'_{nm_1}E\ge \tilde s$ coming from the first relation. (Note that
$\tilde s\ge s$ and, moreover, if we choose $m_0\ge 3$ then
${\mathcal E}(E)^{\perp\widetilde D}={\mathcal E}(E)^{\perp D}$ and $\tilde s=s$.)
By Proposition \ref{ini1}, the product of all the terms
$W_{Y,E,0}(\D,0,\bet^\re,\bet^\ima)$ with $\D$ combined from
$E'\in{\mathcal E}(E)^{\perp D}$ equals $1$, and hence by formula
(\ref{e44}) (resp. (\ref{e44-2}) for $X$ of type $\langle2\rangle^+$) and inequality (\ref{e2040}), one has
$$W_{Y,E,\varphi}(\widetilde D'_{nm_1},0,(\widetilde D'_{nm_1}E)e_1,0)\ge W_{Y,E,\varphi}
(\widetilde D'_{nm_1},\tilde se_1,(\widetilde D'_{nm_1}E-\tilde
s)e_1,0)\qquad\qquad\qquad$$
$$\ge\frac{1}{2}(-K_Y\widetilde D'_{nm_1}-2)!(\widetilde D'_{m_1i(n)}E)(\widetilde
D'_{(n-i(n))m_1}E)\qquad\qquad\qquad\qquad\qquad\qquad\qquad$$
$$\qquad\times\frac{W_{Y,E,\varphi}(\widetilde D'_{m_1i(n)},0,(\widetilde
D'_{m_1i(n)}E)e_1,0)}{(-K_Y\widetilde
D'_{m_1i(n)}-1)!}\cdot\frac{W_{Y,E,\varphi}(\widetilde
D'_{(n-i(n))m_1},0,(\widetilde
D'_{(n-i(n))m_1}E)e_1,0)}{(-K_Y\widetilde D'_{(n-i(n))m_1}-1)!}.$$
There exists $\lam>0$ such that, for all integers $n\ge2$, one has
$$\frac{(\widetilde D'_{m_1i(n)}E)(\widetilde
D'_{(n-i(n))m_1}E)}{2(-K_Y\widetilde
D'_{nm_1}-1)}=\frac{(m_1i(n)\widetilde DE+2-\tilde
s)((n-i(n))m_1\widetilde DE+2-\tilde s)}{2(-nm_1K_Y\widetilde
D-1-\tilde s)}>\lam\ .$$ Hence, the sequence
$$a_n=\frac{W_{Y,E,\varphi}(\widetilde D'_{nm_1},0,(\widetilde
D'_{nm_1}E)e_1,0)}{(-K_Y\widetilde D'_{nm_1}-1)!},\quad n\ge1\ ,$$
satisfies the relation (\ref{e2062}).
Thus, by Lemma \ref{ll1a}, one has
$$\log W_{Y,E,\varphi}(\widetilde D'_{nm_1},0,
(\widetilde D'_{nm_1}E)e_1,0)\ge$$
\begin{equation}
-K_Y\widetilde D'_{m_1} \cdot n \log n + O(n) = -K_XD \cdot m_0m_1n
\log n + O(n), \;\;\; n \to +\infty\ .\label{e2059}\end{equation}
Observe that
$$(n+1)m_0m_1D-E=D'_{(n+1)m_0m_1}+\sum_{E'\in\;{\mathcal E}(E)^{\perp D}}E'\ ,$$
$$D'_{(n+1)m_0m_1+j}-E=D'_{m_0m_1+j}+D'_{nm_0m_1}+\sum_{E'\in{\mathcal
E}(E)^{\perp D}}E',\quad 0\le j<m_0m_1\ ,$$ $$D'_{nm_0m_1}=\widetilde
D'_{nm_1}+nm_1(E+\sum_{E'\in{\mathcal E}(E)^{\perp D}}E')\ ,$$ and hence,
applying formula (\ref{e44}) and inequality (\ref{e2040}) as above
and omitting positive integer coefficients, we obtain
$$W_{Y,E,\varphi}((
(n + 1)m_0m_1+j)D,0,(
(n + 1)m_0m_1+j)(DE)e_1,0)
\qquad\qquad\qquad\qquad\qquad\qquad\qquad\qquad\qquad\qquad$$
$$\ge W_{Y,E,\varphi}(D'_{
(n + 1)m_0m_1+j},0,(D'_{
(n + 1)m_0m_1+j}E)e_1,0) \quad\qquad\qquad\qquad\qquad\qquad\quad$$
$$\ge W_{Y,E,\varphi}(D'_{nm_0m_1},0,(D'_{nm_0m_1}E)e_1,0)
\qquad\qquad\qquad\qquad\qquad\qquad\qquad\qquad\quad$$
$$\ge W_{Y,E,\varphi}(D'_{nm_0m_1}-E-\sum_{E'\in{\mathcal
E}(E)^{\perp
D}}E',0,(D'_{nm_0m_1}E+2-s)e_1,0)\ge\ ...\qquad\quad$$ $$\ge
W_{Y,E,\varphi}(D'_{nm_0m_1}-nm_1(E+\sum_{E'\in{\mathcal
E}(E)^{\perp D}}E'),0,(D'_{nm_0m_1}E+nm_1(2-s))e_1,0)$$
$$\qquad=W_{Y,E,\varphi}(\widetilde D_{nm_1},0,(\widetilde
D_{nm_1}E)e_1,0)\quad\text{for all integers}\ n\ge2,\ 0\le j<m_0m_1\
.\qquad\qquad\quad$$
These inequalities, together with (\ref{e2045}) and (\ref{e2059}),
imply
$$
\displaylines{ \log W(X,nD,F,\varphi) =\log W_{Y,E,\varphi}(nD,0,n(DE)e_1,0)
\cr \ge-K_XD \cdot n \log n + O(n), \;\;\; n \to + \infty\
. }
$$
and hence (\ref{ee21}) by Lemma \ref{las}.

Now suppose that $\tau(D)>0$.
By Lemma \ref{ll1}(ii),  $D-E$ is $Y$-nef, $(D-E)^2\ge0$ and
$\tau(D-E)=\tau(D)-1$.

If $(D-E)^2>0$, then by (\ref{e2057}) and the induction assumption
$$
\displaylines{
\log W(X,nD,F,\varphi)
\ge
\log W(X,n(D-E),F,\varphi)
= \cr
 -K_X(D -E)\cdot n \log n + O(n) =
 -K_XD \cdot n \log n + O(n),}
 $$
which as above implies (\ref{ee21}).

If $(D-E)^2=0$, then by Lemma \ref{ll1}(iii), $D-E=kD''$, where
$k\ge 1$, $D''$ is a primitive $Y$-nef divisor represented by a real
smooth rational curve crossing $E$, and $\dim|D''|=1$. Consider the
divisor $D'_2=2D-E$. It is $Y$-nef and satisfies $D'_2E=2DE+2\ge2$.
It follows from formula (\ref{e44}) (resp. (\ref{e44-2}) for $X$ of type $\langle2\rangle^+$), inequality (\ref{e2040}),
decompositions $$D'_2-E=2kD'',\quad-K_YD'_2-2=2k(-K_YD''-1)+2(k-1)$$
and inequality $D'_2E=2kD''E-2\ge2k-2$ that
$$W_{Y,E,\varphi}(D'_2,0,(D'_2E)e_1,0)\ge
W_{Y,E,\varphi}(D'_2,2(k-1)e_1,(D'_2E-2(k-1))e_1,0)$$ $$\ge
W_{Y,E,\varphi}(D'',0,(D''E)e_1,0)^{2k}=1\ .$$ In the same way from
decompositions
\begin{equation}\begin{cases} &D+nD'_2-E=kD''+nD'_2,\\
&-K_Y(D+nD'_2)-2=k(-K_YD''-1)+(-nK_YD'_2-1)+(k-1)\end{cases}\end{equation}
and
\begin{equation}\begin{cases}&(n+1)D'_2-E=2kD''+[n/2]D'_2+[(n+1)/2]D'_2,\\ &-(n+1)K_YD'_2-2=
2k(-K_YD''-1)+(-[n/2]K_YD'_2-1)\\
&\qquad\qquad\qquad\qquad\qquad+
(-[(n+1)/2]K_YD'_2-1)+2k\end{cases}\end{equation} for all $n\ge2$,
we obtain
$$W_{Y,E,\varphi}(D+nD'_2,0,(D+nD'_2E)e_1,0)\ge
W_{Y,E,\varphi}(nD'_2,0,n(D'_2E)e_1,0)$$ and
$$W_{Y,E,\varphi}((n+1)D'_2,0,(n+1)(D'_2E)e_1,0)\ge\frac{1}{2}(-(n+1)K_YD'_2-2)!\cdot(D'_2E)^2$$
$$\times\left[\frac{n}{2}\right]\cdot\left[\frac{n+1}{2}\right]\cdot\frac{W_{Y,E,\varphi}
([\frac{n}{2}]D'_2,0,[\frac{n}{2}](D'_2E)e_1,0)}{(-[\frac{n}{2}]K_YD'_2-1)!}\cdot
\frac{W_{Y,E,\varphi}([\frac{n+1}{2}]D'_2,0,[\frac{n+1}{2}](D'_2E)e_1,0)}{(-[\frac{n+1}{2}]K_YD'_2-1)!}\,
.
$$
Since there exists $\lam>0$ such that
$$\frac{[n/2]\cdot[(n+1)/2](D'_2E)^2}{2(-(n+1)K_YD'_2-1)}\ge\lam>0\quad\text{for
all}\ n\ge2,\
$$
by Lemma \ref{ll1a} we get
$$
\displaylines{ \log W_{Y,E,\varphi}((2n+j+1)D,0,(2n+j+1)(DE)e_1,0) \ge\cr
\log W_{Y,E,\varphi}(D+nD'_2,0,((D+nD'_2)E)e_1,0)
\ge
 \log
W_{Y,E,\varphi}(nD'_2,0,n(D'_2E)e_1,0)
\ge \cr
-K_YD'_2\cdot  n\log n + O(n)=-K_XD\cdot 2n\log n +O(n)}
$$ which in view of (\ref{e2057}),
(\ref{e2061}), and hence (\ref{ee21}) by Lemma \ref{las}.

\noindent {\bf The case of surfaces of type $\langle q\rangle^+$, $q=3,4$}.
Consider the $(3,0)$-ABV family for a surface $X$ of type $\langle3\rangle^+$ introduced in Section
\ref{sec-deg}. The general fiber of its mirror $(0,3)$-ABV family is a surface of type $\langle2\rangle^+$.
Hence (\ref{ee20}) and (\ref{ee21}) follow from the same relations for surfaces of type
$\langle2\rangle^+$, Corollary \ref{c4},
Proposition \ref{t10}, and relation $W(X,-K_X,F,[\R X\setminus F])=4$
({\it cf.} the table in Section \ref{sec414}). In the same manner, (\ref{ee20}) and (\ref{ee21}) for surfaces of type $\langle4\rangle^+$ follow from these relations for surfaces of type $\langle3\rangle^+$.

\subsubsection{Positivity and asymptotics statements for surfaces of
types $\langle0\rangle^-$
and $1\langle1\rangle^-$
}\label{sec-pos2}
Let $X$ be a real del Pezzo surface as in the title.
Let $F$ be
a non-orientable connected
component
of $\R X$
(in the
case of $X$ of type $\langle0\rangle^-$
assume that $F$ is chosen as specified in Section \ref{F-compatible}). Consider
an elliptic
ABV family of $X$
(see Section \ref{sec-deg}).

The central fiber of
that
ABV family
contains a real nodal del Pezzo pair $(Y,E)$. Observe that, for all divisors
$D\in\PicPPR(Y,E)$, the intersection number $DE$ is even.
Denote by $\hat F$ the connected component of $\R Y$
containing $\R E$, and by $F^+$ the component of $\hat F\setminus\R
E$ to which merges the component $F$ of $\R X$.

Following
Proposition \ref{new-proposition},
we identify $\Pic^\R (X)$  with a subgroup of  $\Pic^\R (Y)$.

\noindent {\bf Positivity}.
Let $D\in\Pic^\R(X)$ be $X$-nef, big, and
$F$-compatible.
Then
$DE=0$, and by formula (\ref{ee14}),
one has
\begin{equation}W(X,D,F,[\R X\setminus
F])=W^-_{Y,F_+,[\R Y\setminus\hat F]}(D,0,0,0)\
.\label{e2051}\end{equation} Thus, to prove (\ref{ee20}) it is
sufficient to show that
for any big, $Y$-nef, and $F_+$-compatible
divisor class $D'\in\Pic^\R(Y)$,
one has
\begin{equation}W^-_{Y,F_+,[\R Y\setminus\hat F]}(D',0,(D'E/2)e_2,0)>0\ .\label{e2052}\end{equation}
In what follows, we prove this inequality by induction on
$\rho(D')=-(K_Y+E)D'$.

As we know, $\rho(D')=-(K_Y+E)D'>0$, since $D'$ is a big
$X$-nef divisor and $\dim\vert -K_X-E\vert=1$.
If
$\rho(D')=1$, then
due to (\ref{e44om}) and inequality (\ref{e2041}) one has
$$W_{Y,F_+,[\R Y\setminus\hat
F]}(D',0,(D'E/2)e_2,0)\ge2^{D'E/2}W_{Y,F_+,[\R Y\setminus\hat
F]}(D',(D'E/2)e_2,0,0)>0\ ,$$
where the latter (strict) inequality follows from the fact that
the real part of the unique rational curve $C\in|D'|_Y$
quadratically tangent to $E$ at $D'E/2$ generic real points
lies in $\overline F_+$.

Let $\rho(D')\ge2$, which
is equivalent to $R_Y(D',0)>0$.
By Lemma \ref{ll1}(iv), the divisor $D''=D'-E-\sum_{E'\in{\mathcal
E}(E)^{\perp D'}}E'$ is $Y$-nef, and it is
easy to see that $D''$ is $F_+$-compatible.
Clearly,
$\rho(D'')\le\rho(D')-2$.
Since $D''E$ is even,
we have $\card({\mathcal E}(E)^{\perp D'})=2s$, $1\le s\le 3$.

So, if $(D'')^2>0$, then formula (\ref{e44om}), inequality
(\ref{e2041}), the induction assumption, and the relations
$$2s=(D'-E-D'')E\le(D'-E)E-2=D'E,\quad D''E/2=D'E/2-s+1$$ result in
\begin{eqnarray}& W^-_{Y,F_+,[\R Y\setminus\hat
F]}(D',0,(D'E/2)e_2,0)\ge2^sW^-_{Y,F_+,[\R Y\setminus\hat F]}(D',se_2,(D'E/2-s)e_2,0)\nonumber\\
& \quad\ge2^s\prod_\D W^-_{Y,F_+,[\R Y\setminus\hat
Y]}(\D,0,0,e_1)\cdot(D''E)\cdot W^-_{Y,F_+,[\R Y\setminus\hat
F]}(D'',0,(D''E/2)e_2,0)>0\ ,\nonumber\end{eqnarray} where $\D$ runs
over all elements $\D\in\PicPPR(Y,E)$ combined from $E'\in{\mathcal
E}(E)^{\perp D'}$.

If $D''=0$ (which is relevant only
if ${\mathcal E}(E)^{\perp D'}\ne\emptyset$),
we get by the same arguments the same expression but
without $(D''E)\cdot W^-_{Y,F_+,[\R Y\setminus\hat F]}(D'',0,(D''E)/2,0)$
in the very end, which again leads to the required positivity.

Thus, it remains to treat the case
$D''\ne0$ and $(D'')^2=0$. Then, since
$D''E$ is positive and even,
Lemma \ref{ll1}(iii)
implies that $D''=-k(K_Y+E)$, $k\ge1$.

If $X$ is of type $1\langle1\rangle^-$ and $F=F^{\;no}$, then the
both $L',L''\in|-K_Y-E|_Y$ (see Section \ref{log-nodal})
are real, $\R L'\cup\R L''\subset\overline F_+$. Using
\begin{equation}R_Y(D',(D'E/2)e_2)=-(K_Y+E)D'+D'E/2-1=2+(k+s-1)-1=k+s\
,\label{e2055}\end{equation}
and applying formula (\ref{e44om}) and inequality (\ref{e2041}), we
derive that
\begin{eqnarray}W^-_{Y,F_+,[\R Y\setminus\hat
F]}(D',0,(D'E/2)e_2,0)&=& W^-_{Y,F_+,[\R Y\setminus\hat
F]}(D',0,(k+s-1)e_2,0,0)\nonumber\\
&\ge&2^{k+s-1} W^-_{Y,F_+,[\R Y\setminus\hat
F]}(D',(k+s-1)e_2,0,0)\nonumber\\
&\ge&2^{k+s-1}\eta
(k)\prod_\D W^-_{Y,F_+,[\R Y\setminus\hat
Y]}(\D,0,0,e_1)>0\ ,\nonumber\end{eqnarray} where the latter
expression corresponds to the summand in the second sum of the
right-hand side of (\ref{e44om}), in which the product runs over
pairs of divisors $\D\in\PicPPR(Y,E)$, combined out of the lines
$E'\in{\mathcal E}(E)^{\perp D'}$, the parameters in the condition (3e) of
Theorem \ref{t3} are chosen to be $\alp^{(0)}=\bet^{(0)}=0$,
$l=k=(D-E)E/2$, and the value of $\eta
(k)$ is given by formula
(\ref{ee6}):
$$\eta
(k)=\begin{cases}(k/2+1)(2-(-1)^{k/2}),\quad & k\ \text{is
even},\\ 2k+2(-1)^{(k-1)/2},\quad & k\ \text{is odd}\end{cases}$$ If
$X$ is
of type $\langle0\rangle^-$,
then $k$ is even by the $F_+$-compatibility condition
({\it cf.}
Proposition \ref{remark-prop}).
Now from
(\ref{e2055}), formula (\ref{e44om}), and inequality (\ref{e2041}),
we derive
$$W^-_{Y,F_+,[\R Y\setminus\hat F]}(D',0,(D'E/2)e_2,0)=
W^-_{Y,F_+,[\R Y\setminus\hat F]}(D',0,(k+s-1)e_2,0,0)$$
$$\qquad\qquad\ge2^{k/2+s-1}W^-_{Y,F_+,[\R Y\setminus\hat
F]}(D',(k/2+s-1)e_2,(k/2)e_2,0,0)$$
$$\ge2^{k/2+s-1}4^{k/2}\prod_\D W^-_{Y,F_+,[\R Y\setminus\hat
Y]}(\D,0,0,e_1)>0\ ,$$ where the latter expression corresponds to
the summand in the second sum in the right-hand side of
(\ref{e44om}), matching the parameter values $\alp^{(0)}=0$,
$\bet^{(0)}=(k/2)e_2$, and $l=0$ in the condition (3e) in Theorem
\ref{t3}.

\noindent {\bf Asymptotics}.
Let $D\in\Pic^\R(X)$ be
$X$-nef, big, and $F$-compatible. In
particular, $DE=0$, which by Lemma \ref{ll1}(ii) yields that $D$
is $Y$-nef, and hence $W^-_{Y,F_+,[\R Y\setminus\hat
Y]}(D,0,0,0)>0$ (see (\ref{e2052})). Since $-K_YD>1$ (see Lemma
\ref{ll1}(i)), formula (\ref{e44om}) for $W^-_{Y,F_+,[\R
Y\setminus\hat F]}(D,0,0,0)$ reads
\begin{equation}W^-_{Y,F_+,[\R
Y\setminus\hat F]}(D,0,0,0)=2W^-_{Y,F_+,[\R Y\setminus\hat
F]}(D-E,0,e_2,0)\ ,\label{e2064}
\end{equation}
thus,
$$W^-_{Y,F_+,[\R Y\setminus\hat F]}(D-E,0,e_2,0)>0\ .$$
Again, using formula (\ref{e44om}) and
non-negativity statement (\ref{e2041}), we obtain
$$W^-_{Y,F_+,[\R Y\setminus\hat F]}(nD-E,0,e_2,0)\ge2(-nK_YD-2)!$$
$$\times\frac{W^-_{Y,F_+,[\R Y\setminus
\hat F ]}([\frac{n}{2}]D-E,0,e_2,0)}{(-[\frac{n}{2}]K_YD-1)!}\cdot
\frac{W^-_{Y,F_+,[\R Y\setminus\hat F]}
(\frac{n+1}{2}]D-E,0,e_2,0)}{(-[\frac{n+1}{2}]K_YD-1)!}$$ for all
$n\ge2$, and hence the sequence $$a_n=\frac{W^-_{Y,F_+,[\R
Y\setminus\hat F]}(nD-E,0,e_2,0)}{(-nK_YD)!},\quad n\ge1\ ,$$
satisfies (\ref{e2063}) with
$$\lam=\inf_{n>3}\frac{[n/2]\cdot[(n+1)/2]K_YD}{n(-nK_YD-1)}>0\
.$$ So, from Lemma \ref{ll1a} we derive
$$\liminf_{n\to\infty}\frac{\log W^-_{Y,F_+,[\R Y\setminus\hat
Y]}(nD-E,0,e_2,0)}{n\log n}\ge-K_YD=-K_XD\ ,$$ and hence in view of
(\ref{e2057}), (\ref{e2064}) applied to $nD$, and of Lemma \ref{las}, we
obtain the desired relation (\ref{ee21}).

\subsection{Proof of Theorem \ref{t11}}\label{sec-s2}
Take a real quadric surface in $\PP^3$ with a spherical real point
set and blow
up
this surface at
three pairs of complex conjugate points. The
resulting del Pezzo surface~$X$ is of degree~$2$ and of type
$\langle 1 \rangle^+$. We have
a natural basis $L_1, L_2, E_1, \ldots, E_6$
in $\Pic(X)$,
where:
\begin{itemize}
\item
$L_1$ and $L_2$ are complex conjugate, $L_1^2 = L_2^2 = 0$, and $L_1
L_2 = 1$,
\item $E_i^2=-1$, and $L_1E_i=L_2E_i=0$ for $1\le i\le 6$,
\item $E_iE_j=0$ for $1\le i<j\le6$,
\item
$E_{2i-1},E_{2i}$ are complex conjugate for $i = 1, 2, 3$.
\end{itemize}
Any real big effective divisors~$D$ in~$X$ can be represented as
$$D=d(L_1+L_2)-d_1(E_1+E_2)-d_2(E_3+E_4)-d_3(E_5+E_6),\quad
d>0,\quad d_1,d_2,d_3\ge0\ ,$$ in particular, $D^2$ is
even.

Consider
an elliptic ABV family
of $X$ (see Section \ref{sec-deg}),
and denote by $(Y,E)$ the real nodal del Pezzo pair in the central fiber of this family.
From
Proposition \ref{new-proposition}
and relation $E^2=-2$,
we immediately derive that the divisor class $E$ is either
$\pm(L_1-L_2)$, or $\pm(E_{2i-1}-E_{2i})$, $i=1,2,3$. Since the corresponding nodal
quartic $Q_Y$
does not have real tangents passing through the node (except for
tangent lines at the node), there
is no real
$(-1)$-curve on $Y$ crossing $E$,
and
we are left with the only option
$E = \pm(L_1 - L_2)$
(we can assume that $E = L_1 - L_2$).

Let $D\in\Pic(X)$ be an $X$-nef and big real divisor. By Theorem
\ref{t4}(1) and
non-negativity statement (\ref{e2040}), one obtains
$$W(X,D,\R X,0)=W_{Y,E,0}(D,0,0,0)\ge0\ ,$$ proving statement (i).
If $W(X,D,\R X,0)=W_{Y,E,0}(D,0,0,0)>0$, then formula (\ref{e44})
applied to $W_{Y,E,0}(D,0,0,0)$ must contain in the right-hand side
a summand $$c\cdot W_{Y,E,0}(D^{(1)},0,e_1,0)\cdot
W_{Y,E,0}(D^{(2)},0,e_1,0)$$ with a positive integer $c$, and
\begin{equation}D^{(1)},D^{(2)}\in\PicPPR(Y,E),\quad D^{(1)}+D^{(2)}=D-E,\quad
W_{Y,E,0}(D^{(i)},0,e_1,0)>0, \ i=1,2\ .\label{e2076}\end{equation}
Here $D^{(1)},D^{(2)}$ cannot be represented by $(-1)$-curves, since
$(D-E)^2=D^2-2\ge0$, and complex conjugate $(-1)$-curves, crossing
$E$, are disjoint.

Assume now that $D^2\le2$. Then $(D-E)^2=D^2-2\le0$, which in the
case of $W(X,D,\R X,0)>0$ leaves the only option
$D^{(1)}=D^{(2)}=D'$, where $(D')^2=0$ and $\dim|D'|_Y=1$; hence
$-K_YD'=2$ and $-K_XD=-2K_YD'=4$ as asserted in statement (iv).

Assume now that $D_1,D_2\in\Pic^\R(X)$ are $X$-nef and big, and
satisfy $W(X,D_i,\R X,0)>0$, $i=1,2$. Show that $W(X,D_1+D_2,\R
X,0)>0$ in agreement with statement (ii). As we have seen above,
there are $D^{(j)}_i\in\PicPPR(Y,E)$, $i,j=1,2$, such that
$$D^{(j)}_iE=1,\quad W_{Y,E,0}(D^{(j)}_i,0,e_1,0)>0,\quad i,j=1,2\ ,$$ $$D^{(1)}_i+D^{(2)}_i=D_i-E,\quad i=1,2\ .$$
Appropriately applying formula (\ref{e44}), we obtain the required
inequality from
$$W(X,D_1+D_2,\R X,0)\ge
c_1W_{Y,E,0}(D^{(1)}_1+D_2,0,e_1,0)W_{Y,E,0}(D^{(2)}_1,0,e_1,0)\ ,$$
and from $$W_{Y,E,0}(D^{(1)}_1+D_2,0,e_1,0)\ge
W_{Y,E,0}(D^{(1)}_1+D_2,e_1,0,0)\qquad\qquad\qquad\qquad\qquad\qquad$$
\begin{equation}\ge
c_{2j}W_{Y,E,0}(D^{(1)}_1,0,e_1,0)W_{Y,E,0}(D^{(1)}_2,0,e_1,0)W_{Y,E,0}(D^{(2)}_2,0,e_1,0)>0\
,\label{e2077a}\end{equation} where $c_1,c_{21}, c_{22}$ are some
positive integers. If $D\in\Pic^\R(X)$ is $X$-nef and big, and is
disjoint
from
a pair of complex conjugate $(-1)$-curves on $X$,
then
we can regard $D$ as a real
divisor on a surface of type $\PP^2_{1,3}$, whose Welschinger
invariants are positive by Theorem \ref{t9}. If $D\in\Pic^\R(X)$ is
$X$-nef and big with $W(X,D,\R X,0)>0$, then
$$W(X,D-K_X,\R X,0)=W_{Y,E,0}(D-K_X,0,0,0)$$ $$\ge
c_0\cdot W_{Y,E,0}(D^{(1)},0,e_1,0)\cdot
W_{Y,E,0}(D^{(2)}-K_X,0,e_1,0)$$
$$\ge2c_0\cdot W_{Y,E,0}(D^{(1)},0,e_1,0)\cdot
W_{Y,E,0}(D^{(2)},0,e_1,0)\cdot W_{Y,E,0}(-K_X-E,0,2e_1,0)>0\ .$$

To complete the proof of statement (ii), we have to show the
positivity of $W(X,-2K_X,\R X,0)$: a direct application of Theorem
\ref{t4}(1) and formula (\ref{e44}) gives $W(X,-2K_X,\R X,0)=8$.

Let $D\in\Pic^\R(X)$ be an $X$-nef and big divisor such that
$W(X,D,\R X,0)>0$, and let $D^{(1)},D^{(2)}\in\PicPPR(Y,E)$ be as in
(\ref{e2076}). By (\ref{e2077a}), $W_{Y,E,0}(D^{(1)}+mD,0,e_1,0)>0$
for all $m\ge0$. Hence, first, for any $n\ge1$,
\begin{eqnarray}W(X,(n+2)D,\R X,0)&=&W_{Y,E,0}((n+2)D,0,0,0)\nonumber\\
&\ge&
c_0W_{Y,E,0}(D^{(2)},0,e_1,0)W_{Y,E,0}(D^{(1)}+(n+1)D,0,e_1,0)\nonumber\\
&\ge& W_{Y,E,0}(D^{(1)}+(n+1)D,0,e_1,0)\ ,\nonumber\end{eqnarray}
and further on, $$W_{Y,E,0}(D^{(1)}+(n+1)D,0,e_1,0)\ge
W_{Y,E,0}(D^{(1)}+(n+1)D,e-1,0,0)$$
$$\frac{(-K_Y(D^{(1)}+(n+1)D)-2)!}{2}\cdot W_{Y,E,0}(D^{(2)},0,e_1,0)$$ $$\times
\frac{W_{Y,E,0}(D^{(1)}+[\frac{n}{2}]D,0,e_1,0)}{(-K_Y(D^{(1)}+[\frac{n}{2}]D)-1)!}
\cdot\frac{W_{Y,E,0}(D^{(1)}+[\frac{n+1}{2}]
D,0,e_1,0)}{(-K_Y(D^{(1)}+[\frac{n+1}{2}]D)-1)!}\ .$$ Then, the
sequence
$$a_n=\frac{W_{Y,E,0}(D^{(1)}+nD,0,e_1,0)}{(-K_Y(D^{(1)}+nD))!},\quad
n\ge0\ ,$$ satisfies (\ref{e2062}) in Lemma \ref{ll1a} with
$$\lam=\inf_{n\ge3}\frac{W_{Y,E,0}(D^{(2)},0,e_1,0)\cdot(-K_Y(D^{(1)}+[\frac{n}{2}]D))
\cdot(-K_Y(D^{(1)}+[\frac{n+1}{2}]D))}{2
(-K_Y(D^{(1)}+(n+1)D))\cdot(-K_Y(D^{(1)}+(n+1)D)-1)}>0\ .$$ Hence
(\ref{e2078}) follows. \proofend

\subsection{Case $\R X = S^2\sqcup \R P^2\#\R P^2$ and $F = S^2$}\label{prooft12}

Let $X$ be a real del Pezzo surface of degree $2$ and of type $1\langle1\rangle^-$, and let $F = F^{\;o}$.
Including such a surface in an ABV-family does not lead to non-negative $w$-numbers,
so the approach used in the proof of Theorems \ref{t9} and \ref{t11} cannot be performed in this case
without any modification.

\begin{proposition}\label{t12}
Let $X$ be a real del Pezzo surface of degree $2$ and of type $1\langle1\rangle^-$.
Then,
\begin{enumerate}\item[{\rm (i)}] for any real $F^{\;o}$-compatible big and nef divisor class $D\in\Pic(X)$, we have for each $i=1,2,3$,
$$\sum_{m\in\Z}W(X,D+2mE^{(i)},F^{\;o},[F^{\;no}])>0\ ,$$
$$\log \sum_{m\in\Z}W(X,nD+2mE^{(i)},F^{\;o},[F^{\;no}]) = -K_X D \cdot n \log n + O(n),
\quad n\to\infty\ ,$$ where $E^{(1)}$, $E^{(2)}$, $E^{(3)}$ are the
real divisors classes introduced in
Section \ref{F-compatible};
\item[{\rm (ii)}] if, in addition, $DE'=0$ for some real $(-1)$-curve $E'\subset X$, then the relations
(\ref{ee20}) and (\ref{ee21}) with $F=F^{\;o}$ hold true.
\end{enumerate}
\end{proposition}

{\bf Proof.}
For statement (i), consider
an external
ABV family
of $X$
(see Section \ref{sec-deg}). The mirror unscrew has a surface $X'$ of type $\langle2\rangle^+$ as a general fiber, and hence both formulas follow from relations (\ref{ee20}) and (\ref{ee21}) for $X'$, claimed in Theorem \ref{t9}, and from Corollary \ref{c4}. 

For statement (ii), we blow down the curve $E'$ and reduce the problem to the case of a real
del Pezzo surface of degree $3$, covered by Theorem \ref{t9}.
\proofend

Using
a slightly different approach, E. Brugall\'e \cite[Corollary 6.10]{B-new} proved the non-negativity of
the invariants $W(X,D,F^{\;o},[F^{\;no}])$ for all real effective divisor classes $D\in\Pic(X)$.

\section{Monotonicity}\label{sec-mon}

\begin{lemma}\label{ln2} ({\it cf.} \cite[Lemma 7.6]{IKS6})
Let $D_1,D_2$ be $X$-nef and big real divisor classes on a del Pezzo
surface $X$ of type $\PP^2_{a,b}$, $a+2b=7$. If $D_2-D_1$ is
effective, then $D_2-D_1$ can be decomposed into a sum
$E^{(1)}+...+E^{(k)}$, where $E^{(i)}$ is either a real
$(-1)$-curve, or a pair of disjoint complex conjugate $(-1)$-curves,
$i=1,...,k$, and, moreover, each real divisor
$D^{(i)}=D_1+\sum_{j\le i}E^{(j)}$ is $X$-nef and big, and satisfies
$D^{(i)}E^{(i+1)}>0$, $i=0,...,k-1$.
\end{lemma}

{\bf Proof}. It is well known that
the effective cone in $\Pic(X)$ is generated by $(-1)$-curves. It is
easy to verify that two complex conjugate $(-1)$-curves in $X$
intersect in at most one point, and
if they intersect, then their sum is
linearly
equivalent to a pair of real $(-1)$-curves, thus, $D_2-D_1$ can be
decomposed into a sum $E^{(1)}+...+E^{(k)}$, where $E^{(i)}$ is
either a real $(-1)$-curve, or a pair of disjoint complex conjugate
$(-1)$-curves, $i=1,...,k$.
We show that a suitable reordering of $E^{(1)},...,E^{(k)}$ ensures
the $X$-nefness and bigness of $D^{(i)}$ together with
$D^{(i)}E^{(i+1)}>0$, $i=0,...,k-1$. The divisor $D^{(0)}=D_1$ is
$X$-nef and big. Suppose now that $D^{(i)}$ is $X$-nef and big for
some $0\le i<k$. If $i=k-1$, then
$D^{(k)}=D_2$ is $X$-nef and big, and furthermore
$$D^{(k-1)}E^{(k)}=(D_2-E^{(k)})E^{(k)}=D_2E^{(k)}-(E^{(k)})^2>0\
.$$ If $i\le k-2$, then there exists $i<j\le k$ such that
$D^{(i)}E^{(j)}>0$. Indeed, otherwise, we would have
\begin{itemize}\item either all $E^{(i+1)},...,E^{(k)}$ orthogonal
to each other and to $D^{(i)}$, and thus, $D_2E^{(j)}<0$, $i<j\le k$
contrary to the $X$-nefness of $D_2$, \item or we would have some
$i<j<j'\le k$ such that $E^{(j)}E^{(j')}>0$, but then
$\dim|E^{(j)}+E^{(j')}|>0$, contradicting to the bigness of
$D^{(i)}$.
\end{itemize} So, we may assume that $D^{(i)}E^{(i+1)}>0$. Then
$D^{(i+1)}=D^{(i)}+E^{(i+1)}$ is $X$-nef and big:
\begin{eqnarray}
D^{(i+1)}E^{(i+1)}&=&D^{(i)}E^{(i+1)}+(E^{(i+1)})^2\ge0\
,\nonumber\\
(D^{(i+1)})^2&=&(D^{(i)})^2
+2D^{(i)}E^{(i+1)}+(E^{(i+1)})^2\ge(D^{(i)})^2>0. \quad
\text{\proofend}\nonumber
\end{eqnarray}

\begin{theorem}\label{tn7}
Let $D_1,D_2$ be $X$-nef and big divisor classes on a real del Pezzo
surface $X$ of type $\PP^2_{a,b}$, $a+2b=7$, $0\le b\le2$, such that
$D_2-D_1$ is effective. Then \begin{equation}W(X,D_2,\R X,0)\ge
W(X,D_1,\R X,0)\ .\label{e2066}\end{equation}
\end{theorem}

{\bf Proof}. By Lemma \ref{ln2}, we should only consider the case of
$E^*=D_2-D_1$ either a real $(-1)$-curve, or a pair of disjoint
complex conjugate $(-1)$-curves.
Let $\lam$ be the number of irreducible components of~$E^*$. We can
assume that $E^*$ consists of $\lam$ exceptional divisors of the
blow up $X\to\PP^2$. Specializing $6-\lam$ other exceptional
divisors so that their blow-downs and the blow-downs of the
components of~$E^*$
appear on a real conic, we degenerate $X$ in
a regular ABV family
into a real nodal del Pezzo pair $(Y,E)$
where $E$
is the strict transform of the above plane conic.
Since $D_2E=D_1E+\lam>D_1E$, we have
$\binom{D_2E+2m}{m}\ge\binom{D_1E+2m}{m}$ for all $m\ge0$, and hence
using the non-negativity statement (\ref{e2040}) and formula
(\ref{ee2}) for the both sides of (\ref{e2066}), we reduce the
problem to establishing inequality
\begin{equation}W_{Y,E,0}(D,\alp,\bet+\lam e_1,0)\ge
W_{Y,E,0}(D-E^*,\alp,\bet,0)\label{e2067}\end{equation} for all
divisors $D\in\PicPPR(Y,E)$ such that $DE\ge\lam$, and for all
vectors $\alp,\bet\in\Z_+^{\infty,\odd}$ such that
$I(\alp+\bet)=DE-\lam$. We prove (\ref{e2067}) by induction on
$R_Y(D,\bet+\lam e_1)$. The case
$R_Y(D,\bet+\lam e_1)<\lam$ is trivial, since then
$R_Y(D-E^*,\bet)=R_Y(D,\bet+\lam e_1)-\lam<0$. If $R_Y(D,\bet+\lam
e_1)=\lam$ and, respectively, $R_Y(D-E^*,\bet)=0$, the only relevant
case is that of Proposition \ref{ini1}(1iii) with $D-E^*$ playing
the role of $D$ and $\bet=0$, in which case by (\ref{e2040}),
formula (\ref{e44}), and Proposition \ref{ini1}(1iii) we have
$$W_{Y,E,0}(D,\alp,\lam e_1,0)\ge W_{Y,E,0}(D,\alp+\lam e_1,0,0)=1=W_{Y,E,0}(D-E^*,\alp,0,0)\ .$$ If
$R_Y(D,\bet+\lam e_1)=R_Y(D-E^*,\bet)+\lam>\lam$, we compute both
sides of (\ref{e2067}) by formula (\ref{e44}) and compare them using
the induction assumption (in the sequel we shortly write
RHS(\ref{e44})$_l$ and RHS(\ref{e44})$_r$ for the right-hand side of
(\ref{e44}) expressing the left and the right terms of (\ref{e2067})
respectively). So, for the summands in the first sum in
RHS(\ref{e44})$_l$ and RHS(\ref{e44})$_r$ the induction assumption
yields
$$W_{Y,E,0}(D,\alp+e_j,\bet-e_j+\lam e_1,0)\ge
W_{Y,E,0}(D-E^*,\alp+e_j,\bet-e_j,0)\ .$$ For the second sum in
RHS(\ref{e44})$_l$ and RHS(\ref{e44})$_r$ we perform the following
comparison. Let
$$
S_r=c\cdot\frac{2^{\|\bet^{(0)}\|}}{\bet^{(0)}!}\cdot\frac{(n-1-\lam)!}
{n_1!...n_m!}\cdot \prod_{i=1}^m\binom{(\bet^\re)^{(i)}}{\gam^{(i)}}
W_{Y,E,0}(\D^{(i)},\alp^{(i)},(\bet^\re)^{(i)},(\bet^{\ima})^{(i)})
$$
be a summand in the second sum of RHS(\ref{e44})$_r$,
where $n=R_Y(D,\bet+\lam e_1)$ and
$n_i=R_Y(\D^{(i)},(\bet^\re)^{(i)}+2(\bet^{\ima})^{(i)})$, $1\le
i\le m$. Notice that $m\ge1$, since $DE^*>0$, and hence there is
$\D^{(i)}$ such that $[\D^{(i)}]E^*>0$. Pick $\D^{(j)}$ with
$[\D^{(j)}]E^*>0$ and associate with $S_r$ the following summand
$S_l$ in the second sum of RHS(\ref{e44})$_l$:
\begin{itemize}
\item if $[\D^{(j)}]\ne-\lam(K_Y+E)-E^*$ (in which case $\D^{(j)}$ is a real divisor), then
\begin{equation}S_l=c\cdot\frac{s_j^r}{s_j^l}
\cdot\frac{2^{\|\hat\bet^{(0)}\|}}{\hat\bet^{(0)}!}\cdot\frac{(n-1)!}
{\hat n_1!...\hat
n_m!}\cdot\prod_{i=1}^m\binom{(\hat\bet^\re)^{(i)}}{\gam^{(i)}}W_{Y,E,0}
(\hat \D^{(i)},\alp^{(i)},
(\hat\bet^\re)^{(i)},(\bet^{\ima})^{(i)}),\label{e2069}\end{equation}
where $\hat\bet^{(0)}=\bet^{(0)}$, $\hat n_i=R_Y(\hat \D^{(i)},
(\hat\bet^\re)^{(i)}+2(\bet^{\ima})^{(i)})$, $$\hat \D^{(i)}
=\begin{cases}\D^{(i)},\quad&i\ne j,\\
\D^{(j)}+E^*,\quad&i=j,\end{cases}\quad\text{and}\quad(\hat\bet^\re)^{(i)}
=\begin{cases}(\bet^\re)^{(i)},\quad&i\ne j,\\ (\bet^\re)^{(j)}+\lam
e_1,\quad&i=j,\end{cases}$$ $s_j^r$ counts how many times the tuple
$(\D^{j},\alp^{(j)},(\bet^\re)^{(j)},\gam^{(j)})$ occurs in the list
$\{(\D^{(i)},\alp^{(i)},(\bet^\re)^{(i)},\gam^{(i)})\}_{i=1}^m$, and
$s_j^l$ counts how many times the tuple $(\hat
\D^{j},\alp^{(j)},(\hat \bet^\re)^{(j)},\gam^{(j)})$ occurs in the
list $\{(\hat
\D^{(i)},\alp^{(i)},(\hat\bet^\re)^{(i)},\gam^{(i)})\}_{i=1}^m$;
\item if $D^{(j)}=-\lam(K_Y+E)-E^*$, then
\begin{equation}S_l=c\cdot\frac{2^{\|\hat\bet^{(0)}\|}}{\hat\bet^{(0)}!}
\cdot\frac{(n-1)!}{\hat n_1!...\hat n_m!}
\cdot\prod_{\renewcommand{\arraystretch}{0.6}
\begin{array}{c}
\scriptstyle{1\le i\le m}\\
\scriptstyle{i\ne j}
\end{array}}\binom{(\hat\bet^\re)^{(i)}}{\gam^{(i)}}
W_{Y,E,0}(\hat
\D^{(i)},\alp^{(i)},(\hat\bet^\re)^{(i)},(\bet^{\ima})^{(i)})\
,\label{e2070}\end{equation} where $\hat\bet^{(0)}=\bet^{(0)}+\lam
e_1$, $(\hat\bet^\re)^{(i)}=(\bet^\re)^{(i)}$, $\hat n_i=n_i$, and
$\hat \D^{(i)}=\D^{(i)}$ for $1\le i\le m$, $i\ne j$, and $\hat
n_j=0$.
\end{itemize}
It is easy to verify (again using the induction assumption) that
$$S_r\le \begin{cases}S_l\cdot\frac{s_j^l\hat n_j}{s_j^r(n-1)},
\quad &\text{if}\ \lam=1,\\
S_l\cdot\frac{s_j^l\hat n_j(\hat n_j-1)}{s_j^r(n-1)(n-2)},\quad
&\text{if}\ \lam=2,\end{cases}\quad\quad\text{in (\ref{e2069})}\ ,$$
and
$$S_r\le \begin{cases}S_r\cdot\frac{\hat\bet^{(0)}_1}{2(n-1)},
\quad&\text{if}\ \lam=1,\\
S_r\cdot\frac{\hat\bet^{(0)}_1(\hat\bet^{(0)}_1-1)}{4(n-1)(n-2)},
\quad&\text{if}\ \lam=2,\end{cases}\quad\quad \text{in
(\ref{e2070})}\ .$$ Since $n-1=\sum_js^l_j\hat
n_j+\|\hat\bet^{(0)}\|$ ({\it cf.} Remark \ref{r2}), we conclude
that the total value of the terms in the second sum in
RHS(\ref{e44})$_r$ associated with a given summand $S_l$ in the
second sum in RHS(\ref{e44})$_l$ does not exceed $S_l$, which
completes the proof. \proofend

\section{Mikhalkin's congruence}\label{secn22}

\begin{theorem}\label{t8}
For any $X$-nef and big divisor class $D$ on a surface $X$ of type
$\PP^2_{7,0}$, one has
$$W(X,D,\R X,0)
=
GW_0(X,D)\mod 4\ .$$
\end{theorem}

{\bf Proof}. Using
a regular ABV family
of $X$
and
formulas (\ref{ee1}),
(\ref{ee2}), we reduce the question to the congruence
$$W_{Y,E,0}(D,0,(DE)e_1,0)
=
N_Y(D,0,(DE)e_1)\mod 4$$ for all
divisors
$D\in\PicPPR(Y, E)$.

In fact, a more general statement holds: for any divisor class
$D \in \PicPPR(Y, E)$
and any vectors $\alp,\bet\in\Z^\infty_+$ such that
$I(\alp+\bet)=DE$, one has
\begin{equation}
W_{Y,E,0}(D,\alp,\bet,0)
=
I^\bet N_Y(D,\alp,\bet)\mod 4 \quad
\text{\rm if} \ \beta \in Z^{\infty, \; \odd}_+ ,\label{e112}
\end{equation}
and
\begin{equation}
I^{\bet}\cdot N_Y(D,\alp,\bet)
=
0\mod4\quad\text{if} \ \beta
\not\in\Z_+^{\infty,\;\odd} \ ,\label{eCO}
\end{equation}
where the numbers $N_Y(D,\alp,\bet)$
are the degrees of varieties $V_Y(D, \alp, \bet, \bp^\flat)$
(here $\bp^\flat=\{p_{i,j}\ :\ i\ge1,\ 1\le j\le\alp_i\}$
are sequences
of $\|\alp\|$ distinct generic points on $E$)
introduced in Section \ref{Welsch-inv}.
The proof literally coincides with that of \cite[Theorem 5]{IKS7}
and uses induction on $R_Y(D,\bet)$ and the
recursive formulas \cite[Formula (66)]{MS2} and (\ref{e44}).
\proofend

\newpage

\centerline{\bf Index of notations}

\vskip10pt

$$
\begin{array}{lllll}
\|\alpha\| & \mbox{Section \ref{some_notations}}& \qquad & \Pic^{\R}(\Sig)& \mbox{Section \ref{sec414}} \\
I\alpha & \mbox{Section \ref{some_notations}}& \qquad & \Pic^{\R}_+(\Sig)& \mbox{Section \ref{div-cl}} \\
I^{\alpha }& \mbox{Section \ref{some_notations}}& \qquad & \PicPP(\Sig, E)& \mbox{Section \ref{div-cl}} \\
e_i & \mbox{Section \ref{some_notations}}& \qquad & \PicPPR(\Sig)& \mbox{Section \ref{div-cl}} \\
\Z^{\infty}_+ & \mbox{Section \ref{some_notations}}& \qquad & [\D]&  \mbox{Section \ref{div-cl}} \\
\Z^{\infty, \; \odd}_+ & \mbox{Section \ref{some_notations}}& \qquad & R_{\Sig}(\D, \beta) & \mbox{Section \ref{div-cl}}\\
\Z^{\infty, \; \even}_+ & \mbox{Section \ref{some_notations}}& \qquad & \PP^2_{a, b} & \mbox{Section \ref{sec12}}\\
\Z^{\infty, \; \odd \cdot \even}_+ & \mbox{Section \ref{some_notations}}& \qquad & \langle q \rangle^{\pm} & \mbox{Section \ref{sec12}}\\
W(X, D, F, \varphi) & \mbox{Section \ref{intro1}}& \qquad & \langle1\langle1\rangle\rangle^{\pm} & \mbox{Section \ref{sec12}} \\
W_{Y,E,\varphi}(\D,\alp,\bet^{\re},\bet^{\ima}) & \mbox{Section \ref{new-ordinary}}&
\qquad & L', L'' & \mbox{Section \ref{log-nodal}} \\
W_{Y,F_+,\varphi}^\pm(\D,\alp,\bet^{\re},\bet^{\ima}) & \mbox{Section \ref{sec6}}&
\qquad & {\mathcal E}(E) & \mbox{Section \ref{log-nodal}} \\
U_{Y,F,\varphi}(D, ke_1, le_1, \bp^\sharp) & \mbox{Section \ref{unum}} &
\qquad & {\mathcal E}(E)^{\perp D} & \mbox{Section \ref{log-nodal}} \\
U_{Y,F_+,\varphi}^{\pm}(D, 0, (DE/2)e_1, \bp^\sharp) & \mbox{Section \ref{unum}} &
\qquad & Q_X & \mbox{Section \ref{sec12}} \\
U_{Y,E,\varphi'}(D,\bp^\sharp) & \mbox{Section \ref{eun}}&
\qquad & F^{\; o} & \mbox{Section \ref{sec12}} \\
V^\R_Y (\D,\alp,\bet^{\re},\bet^{\ima},\bp^\flat) & \mbox{Section \ref{Welsch-inv}} &
\qquad & F^{\; no} & \mbox{Section \ref{sec12}} \\
V^\R_Y (\D,\alp,\bet^{\re},\bet^{\ima},\bp^\flat,\bp^\sharp) & \mbox{Section \ref{Welsch-inv}} &
\qquad & F^{\; no} & \mbox{Section \ref{sec12}} \\
V^\R_{Y,F_+}(\D,\alpha,\beta^{\re},\beta^{\ima},\bp^\flat,\bp^\sharp) & \mbox{Section \ref{sec6}} &
\qquad & V_Y(D,\alp,\bet,\bp^\flat) & \mbox{Section \ref{Welsch-inv}} \\
 \mu_\varphi(C) & \mbox{Section \ref{new-ordinary}}  &
\qquad & GW_0(Y,D) & \mbox{Section \ref{sec-vakil}} \\
\mu^\pm_\varphi(C) & \mbox{Section \ref{sec6}}  &
\qquad & N_Y(D,0,(DE)e_1) & \mbox{Section \ref{sec-vakil}}
\end{array}
$$

\end{document}